\documentstyle[11pt,twoside,amssymb,amscd]{article}

\def\Ann{{\rm Ann}}
\def\Perv{{\rm Perv}}
\def\Comp{{\rm Comp}}
\def\Ind{{\rm Ind}}

\newcommand{\C}{{\Bbb C}}
\newcommand{\Z}{{\Bbb Z}}
\newcommand{\R}{{\Bbb R}}

\let\ccal=\cal
\renewcommand{\cal}[1]{{\ccal #1}}
\title{\bfseries Geometric methods in the representation theory\\
 of Hecke algebras and quantum groups}
\author{Victor GINZBURG\\[5pt]
\itshape Department of Mathematics \\
\itshape University of Chicago \\
\itshape Chicago, IL 60637 \\
\itshape USA\\[12pt]
Notes by \\ Vladimir BARANOVSKY}
\date{}

\begin{document}
\maketitle

\begin{abstract}
\noindent
These  lectures 
given in Montreal in Summer 1997 
are mainly based on, and form a condensed survey of the book by {\rm N.
  Chriss and V. Ginzburg,} {\it Representation Theory and Complex
  Geometry} , {\rm  Birkh\"auser 1997}.  Various algebras arising
naturally in Representation Theory such as the group algebra of a Weyl
group, the universal enveloping algebra of a complex semisimple Lie
algebra, a quantum group or the Iwahori-Hecke algebra of bi-invariant
functions (under convolution) on a $p$-adic group, are considered.  We
give a uniform geometric construction of these algebras in terms of
homology of an appropriate ``Steinberg-type" variety $Z$ (or its
modification, such as K-theory or elliptic cohomology of $Z$, or an
equivariant version thereof). We then explain how to obtain a complete
classification of finite dimensional irreducible representations of
the algebras in question, using our geometric construction and
perverse sheaves methods.

Similar techniques can be applied to other algebras, e.g.
the double-affine Hecke algebras, elliptic algebras,
quantum toroidal algebras.
\end{abstract}

\section*{Introduction}

A new branch has emerged during the last decade within the part of
mathematics dealing with Lie groups.  That new branch may be called
{\it Geometric Representation Theory}.  As Beilinson-Bernstein put it
in their seminal paper [BeBe], the discovery of $\cal D$-modules and
of {\it perverse sheaves} has made Representation Theory, to a large
extent, part of Algebraic Geometry.  Among applications of perverse
sheaf methods to Representation Theory that have already proved to be
of primary importance we would like to mention here the proof of the
Kazhdan-Lusztig conjecture by Beilinson-Bernstein [BeBe] and
Brylinski-Kashiwara [BrKa], Lusztig's construction of canonical bases
in quantum groups [L1], and the work of Beilinson-Drinfeld on the
Geometric Langlands conjecture, cf. [Gi3].  We refer the reader to
[L2] for further applications.

In these   notes we discuss  another (not  completely unrelated to the
above) kind of   applications of equivariant  $K$-theory and  perverse
sheaves to representations  of Hecke algebras and  quantum groups.  We
study  various   associative   algebras   that   arise  naturally   in
Representation Theory.   These may be, for   example, either the group
algebra of   a Weyl group, or  the  universal enveloping algebra  of a
complex  semisimple Lie  algebra, or a    quantum group, or the  Hecke
algebra of  bi-invariant functions (under   convolution) on a $p$-adic
group.  Further  examples, such as the  double-affine Hecke algebra of
Cherednik, elliptic  algebras,  quantum  toroidal algebras, etc.,  fit
into the same  scheme, but will not  be  considered here; see  [GaGr],
[GKV1], [GKV2], [Na] for  more details.  In spite  of the diversity of
all the examples  above,  our strategy  will  always follow  the  same
pattern that we now outline.

The first step consists of giving an ``abstract-algebraic" presentation
of our  algebra $A$ in   terms of a  convenient  set of generators and
relations.  The  second  step is  to find  a geometric construction of
$A$.  More specifically, we are looking for a complex manifold $M$ and
a ``correspondence" $Z\subset  M\times M$ such that  the algebra $A$ is
isomorphic to the homology $H_*(Z, {\Bbb C})$  or its modification, such
as $K$-theory or elliptic cohomology of $Z$, or an equivariant version
thereof. Here the subvariety $Z$ that we are seeking should be thought
of as the graph of a  multivalued map $f : M  \to M$, and this map $f$
should  satisfy the {\it  idempotency equation}: $f\circ f=f$. Such an
equation rarely holds  for genuine maps, but  becomes not so  rare for
multivalued  maps. We will see that  the idempotency  equation for $f$
gives rise to a multiplication-map on homology: $H_*(Z, {\Bbb C}) \times
H_*(Z,  {\Bbb C}) \to H_*(Z, {\Bbb  C})$, called {\it convolution}. Such a
convolution  makes   $H_*(Z,  {\Bbb   C})$  an    associative, typically
noncommutative, ${{\Bbb  C}}$-algebra, and it  is this algebra structure
on homology (or $K$-theory or  elliptic cohomology) of $Z$ that should
be isomorphic to the one on $A$.

It should be mentioned that, in all examples above, the only known way
of proving an isomorphism $A \simeq H_*(Z, {\Bbb C})$ is by showing
that the algebra on the right hand side has the same set of generators
and relations as were found for $A$ in Step 1.  Of course, given an
algebra $A$, there is no {\it a-priori} recipe helping to find a
relevant geometric data $(M,Z)$; in each case this is a matter of good
luck.  Sometimes a partial indication towards finding (in a
conceptual way) a geometric realization of our algebra $A$ as the
convolution algebra $H_*(Z, {\Bbb C})$ comes from $\cal D$-modules;
more precisely, from the notion of the {\it characteristic cycle} of a
holonomic $\cal D$-module (this links our subject to that discussed in
the first paragraph, see [Gi1]).  It is fair to say, however,
that it is still quite a mystery, why a geometric realization of the
algebras $A$ that we are interested in is possible at all.  But once a
geometric realization is found, a complete classification of finite
dimensional irreducible representations of $A$ can be obtained in a
straightforward manner.  This constitutes the last step of our
approach, which we now outline and which is to be explained in more
detail in Section 5.

The geometric realization of the algebra reduces the problem to the
classification of finite dimensional irreducible representations of a
convolution algebra, like $H_*(Z, {\Bbb C})$.  This problem is solved
as follows.  First we show, using the techniques of sheaf theory (see
Sections 4 and 5), that the convolution algebra is isomorphic to the
Ext-algebra, ${\rm Ext}^\bullet(\cal L, \cal L)$, equipped with the Yoneda
product, where $\cal L$ is a certain constructible complex on an
appropriate complex variety.  The structure of $\cal L$ is then
analyzed using the very deep {\it Decomposition Theorem} [BBD].  The
theorem yields an explicit decomposition of the Ext-algebra as the sum
of a nilpotent ideal and a direct sum of finitely many matrix
algebras. Hence the nilpotent ideal is the {\it radical}, and each
matrix algebra occurring in the direct sum gives an irreducible
representation of the Ext-algebra.  Therefore, the non-isomorphic
irreducible representations are parametrized by the matrix algebras
that occur in the decomposition above.  Thus, the classification of
finite dimensional irreducible representations of the original algebra
$A$ can be read off from the decomposition of the constructible
complex $\cal L$ provided by the Decomposition Theorem.

These  notes are  mainly based on, and  form a condensed survey
of, the book [CG].  The reader is referred to the introduction to [CG]
for more  motivation  and   historical background.  We    have  tried,
however, to  make our present exposition  as complementary  to [CG] as
possible. For example, the geometric  construction of Weyl groups  and
enveloping algebras given here  is based on Fourier transform, whereas
in [CG] another   approach  has  been  used.  We  also  discuss   here
degenerate  affine Hecke  algebras and quantum  affine algebras, which
were not present in [CG].

\vskip 0.5cm
\noindent{ TABLE OF CONTENTS}

\medskip\noindent
\begin{tabular}{rl}
\small
1& Borel--Moore homology.\\
2& Convolution in Borel--Moore homology.\\
3& Constructible complexes.\\
4& Perverse sheaves and the Decomposition Theorem.\\
5& Sheaf--theoretic analysis of the convolution algebra.\\
6& Representations of Weyl groups.\\
7& Springer theory for $\cal U( {\frak s} {\frak l}_n)$.\\
8& Fourier transform.\\
9& Proof of the geometric construction of $W$.\\
10& Proof of the geometric construction  of $\cal U( {\frak s} {\frak l}_n)$.\\
11& $q$-deformations: Hecke algebras and a quantum group.\\
12& Equivariant cohomology and degenerate versions.
\end{tabular}

\section{Borel-Moore homology}

Borel-Moore homology will   be the principal functor  we  use in these
lectures  for constructing representations  of Weyl groups, enveloping
algebras, and   Hecke algebras. We    review here the   most essential
properties of the Borel-Moore homology theory and  refer the reader to
the monographs [Bre] and  [Iv] for  a more  detailed treatment  of the
subject.  We have to say a few words about the kind  of spaces we will
be dealing with.

By a ``space" (in the topological sense) we will mean a locally compact
topological  space  $X$  that  has  the  homotopy   type of  a  finite
$CW$-complex; in particular,  has  finitely many  connected components
and   finitely    generated homotopy  and  homology  groups   (with
$\Z$-coefficients).  Furthermore, our space $X$  is assumed to admit a
closed embedding into  a countable at infinity $C^\infty$-manifold $M$
(in particular, $X$ is paracompact).  We assume also that there exists
an open neighborhood $U \supset X$ in $M$ such that  $X$ is a homotopy
retract of $U$.  It is known, cf.  [GM],  [RoSa], that any complex  or
real algebraic variety satisfies  the above conditions.  These are the
spaces we will mainly use in applications.

Similarly, by   a     closed  ``subset"  of  a
$C^\infty$-manifold we   always mean a subset  $X$   which has an open
neighborhood $U \supset X$ such that $X$ is a homotopy retract of $U$.
In that case one can also find a  smaller {\it closed} neighborhood $V
\subset U$  such that $X$  is a  {\it proper}  homotopy retract of $V$
(recall that a  continuous map $f: X \to  Y$ is called {\it proper} if
the  inverse image of  any compact set is compact).

We now give a list of the various equivalent definitions of
Borel-Moore homology of a space $X$, see [BoMo], [Bre].  In what
follows, all homology and cohomology is taken with complex
coefficients, which may be replaced by any field of characteristic
zero.

\smallskip
$(1)\enspace$ Let  $\hat{X}=X\cup  \{\infty\}$ be  the one-point
compactification   of $X$.  Define $H_*^{BM}(X)= H_*(\hat{X},\infty)$,
where    $H_*$  is   ordinary    relative  homology    of  the    pair
$(\hat{X},\infty)$.

\smallskip
$(2)\enspace$ Let $\overline{X}$ be an arbitrary compactification of $X$
such that $(\overline{X},\overline{X}\backslash X)$ is a 
$CW$-pair. Then,  
$H_*^{BM}(X)\simeq H_*(\overline{X},\overline{X}\backslash X)$;
see [Sp].
The fact that this definition agrees with (1) is proved in [Bre].

\smallskip
$(3)\enspace$ Let $C_*^{BM}(X)$ be the chain complex of 
{\it infinite} singular
chains
$
\sum_{i=0}^{\infty} a_i\sigma_i,
$
where $\sigma_i$ is  a singular simplex, $a_i\in {\Bbb  C}$, and the sum
is locally  finite in the  following  sense:  for any  compact set  $D
\subset  X$ there are only  finitely  many non-zero coefficients $a_i$
such that $D\cap  \mbox{supp}\,\sigma_i \not = \emptyset\,.$ The usual
boundary  map    $\partial$ on singular  chains  is   well  defined on
$C_*^{BM}(X)$ because taking boundaries  cannot destroy the finiteness
condition.  We then have
$$
H_{_\bullet}^{BM}(X)= H_{_\bullet}(C_*^{BM}(X)\,,\,\partial).
$$

\smallskip
$(4)\enspace$
Poincar\'e duality: let $M$ be a smooth, oriented manifold, and
$\dim_\R M=m$.
Let
$X$ be a closed subset of $M$ which has a closed neighborhood
$U \subset M$ such that $X$ is a proper deformation retract of $U$.
 Then there is a canonical isomorphism ([Iv], [Bre]):
\begin{equation}
\label{pd_defn_bm}
H_i^{BM}(X) \simeq H^{m-i}(M,M\setminus X),
\end{equation}
where each side of the equality is understood to be with complex
coefficients. In particular, setting $X=M$  we obtain, for any {\it smooth}
not necessarily compact variety $M$, a canonical isomorphism
(depending on the orientation of~$M$)
\begin{equation}
\label{standard_pd}
H^{BM}_i(M) \simeq H^{m-i}(M)\,.
\end{equation}
We will often use the ``Poincar\'e duality'' definition (formula
(\ref{pd_defn_bm}) above)
to prove many of the basic
theorems about Borel-Moore homology by appealing to the
same theorems for singular cohomology.  In these instances we will
refer the reader to [Bre], [Sp] for the proofs in singular cohomology,
despite the fact that Borel-Moore homology is not explicitly developed
there.

\medskip

\noindent
{\bf Notation } From now on $H_*$ will stand for
$H_*^{BM}$ (since Borel-Moore homology is the main functor used in
these notes).

\medskip

We now study the functorial properties of Borel-Moore homology.

\medskip

\noindent
{\bf Proper pushforward}\enspace
\label{bm_covariance}
Borel-Moore homology is a covariant functor with respect to
proper maps.  
If ${f}:X\to Y$ is a proper map, then we may define the direct image 
(or proper push-forward) map 
$$
{f}_*:H_*(X)\to H_*(Y)
$$
by extending $f$ to a map
${\bar{f}}:\bar{X}\to \bar{Y}$ where $\bar{X}=X\cup \{\infty\}\,,$
 resp. $\bar{Y}= Y\cup \{\infty\}\,,$ and $\bar{f}(\infty)=\infty$
(observe that $f$ being proper ensures that $\bar{f}$ is continuous).

\medskip

\noindent
{\bf Long exact sequence of Borel-Moore homology}\enspace
\label{long_exac}
Given an open subset $U  \subset  X$ there is a natural restriction
morphism $H_*(X) \to H_*(U)$ induced by the composition of maps:
\[ H_*(X)= H^{ord}_*(\overline{X}, \overline{X} \setminus X)
 \to H^{ord}_*(\overline{X}, \overline{X} \setminus
U) = H_*(U)\,,\]
where $\overline{X}$ stands for a compactification of $X$,
cf. definition (2) of Borel-Moore homology, and the map in the middle
is induced by the natural morphism of pairs $(\overline{X}, \overline{X}
\setminus X)
\to (\overline{X}, \overline{X} \setminus
U)\,.$
For an alternative ad hoc definition of the restriction to
an open subset see [Iv].

Suppose that next $F$ is a closed subset of $X$.  Write $i: F
\hookrightarrow X$ for the (closed) embedding, set $U= X \setminus F$,
and consider the diagram
\[
F\stackrel{i}{\hookrightarrow} X \stackrel{j}{\hookleftarrow} U.
\] 
 Since $i$ is
proper and $j$ is an open embedding, the functors $i_*$ and $j^*$
are  defined.  Then there is a natural
 long exact sequence in Borel-Moore homology (see [Bre],
[Sp] for more details):
\begin{equation}
\label{longexact_homology}
\cdots \to H_p(F) \to H_p(X) \to H_p(U) \to H_{p-1}(F)\to \cdots 
\end{equation}
To construct this long exact sequence, choose an embedding of $X$ as a
closed subset  of a smooth manifold $M$.   Then the Poincar\'e duality
isomorphism (1) gives:
\[H^{m-p}(M, M \setminus X) \simeq H_p(X)\quad\mbox{and}\quad
H^{m-p}(M, M \setminus F) \simeq H_p(F).\]
Further, the set $U$ being locally closed in $M$, we may find an open
subset $M' \subset M$ such that $U$ is a closed subset of $M'$.
Then, the excision axiom, see [Sp], combined with
Poincar\'e duality 
yields
\[H^{m-p}(M, M \setminus U) \simeq H^{m-p}(M', M' \setminus U) \simeq
H_p(U).\] 
Thus, we see that terms of the standard relative
cohomology long exact sequence, cf. [Sp]:
\begin{equation}
\label{longexact_cohomology}
\ldots \to H^k(M, M \setminus F) 
\to H^k(M, M \setminus X) 
\to H^k(M, M \setminus U) \to
H^{k+1}(M, M \setminus F) \to \ldots 
\end{equation}

\medskip
\noindent
get identified via the above isomorphisms with the corresponding terms
of      (\ref{longexact_homology}).   In this  way   we      define
(\ref{longexact_homology}) to  be  the exact  sequence induced  by the
cohomology exact sequence (\ref{longexact_cohomology}).

\medskip

\noindent
{\bf Fundamental class}\enspace
\label{remark_lecture7}
The main reason we are using Borel-Moore homology is the existence of
fundamental classes. Recall that 
any smooth  oriented manifold $X$ has a  well defined {\it fundamental
  class} in Borel-Moore homology:
$$
[X]\in H_m(X),\qquad m=\dim_{{\Bbb R}} X.
$$ 
Note that  there is no fundamental class   in ordinary homology unless
$X$ is  compact.

The essential feature of Borel-Moore homology is the 
existence of a fundamental
class, $[X]$, of any (not necessarily smooth or compact) 
complex
algebraic variety $X$.
If $X$ is irreducible of real dimension $m$, then $[X]$ 
is the unique class in $H_m(X)$ that
restricts to the fundamental class of the non-singular part of $X$.
More precisely, write $X^{reg}$ for the Zariski open dense subset 
consisting of the non-singular points  of $X$. Being a smooth complex
manifold, $X^{reg}$ has a canonical orientation coming from the complex
structure, and hence a fundamental class $[X^{reg}] \in H_m(X^{reg})$.
The inequality $\dim_{{\Bbb R}}(X \backslash X^{reg})\le m-2 $ 
yields (say by definition (1) of Borel-Moore homology)
\[
H_k(X \backslash X^{reg})=0\quad\mbox{for any}\quad k> m-2.
\] 
The long exact sequence of Borel-Moore homology (\ref{longexact_homology}) 
shows that the restriction $H_{m}(X)\to H_{m}(X^{reg})$
is an isomorphism. We define $[X]$ to be the preimage of $[X^{reg}]$
under this isomorphism.
If $X$ is an arbitrary complex algebraic 
variety with irreducible components $X_1, X_2,\dots, X_n$, 
then $[X]$ is set to be a non-homogeneous class equal to  $\sum [X_i]$.

The top Borel-Moore homology of a complex algebraic variety is
particularly easy to understand in the light of the following proposition.

\medskip

\noindent
{\bf Proposition 1.1}  {\em Let $X$ be a complex algebraic variety of
complex dimension $n$ and let  $X_1,\dots,X_m$ be the  $n$-dimensional
irreducible components of $X$.   Then the fundamental classes $[X_1]$,
$\dots$,    $[X_m]$    form  a    basis    for  the     vector   space
$H_{top}(X)=H_{2n}(X)$.$\quad\square$ }

\bigskip

\noindent
{\bf Intersection Pairing} \enspace
\label{int_pairing}
Let $M$ be a smooth, oriented manifold
and $Z,\tilde{Z}$ two
closed subsets (in the sense explained at the beginning of this section) in $M$.  We define
a bilinear pairing
\begin{equation}
\label{qs_pairing}
{\cap}:H_i(Z)\times H_j(\tilde{Z})\to H_{i+j-m}(Z\cap \tilde{Z}),
\quad m=\dim_{{\Bbb R}} M
\end{equation}
which refines the standard intersection of cycles in a smooth variety.
The only new feature is that instead of regarding cycles as homology
classes in the ambient manifold $M$ we take their supports into
account.
So, given two singular chains with supports in the subsets
 $Z$ and $\tilde{Z}$, respectively,
we would like to define their intersection
to be a class in the homology of the set-theoretic intersection, 
$Z\cap \tilde{Z}$. To that end we
use the
standard $\cup$-product in relative cohomology
(cf. [Sp]): 
$$
{\cup}:H^{m-i}(M,M\setminus Z)\times H^{m-j}(M, M\setminus
\tilde{Z}) \to H^{2m-j-i}(M,\,(M\setminus Z)\cup (M\setminus
\tilde{Z}))\,.
$$
Applying 
Poincar\'e  duality (1) to each term of this $\cup$-product
we get the intersection pairing (\ref{qs_pairing}).

The intersection pairing introduced above has an especially clear
geometric meaning in the case when $M$ is a real analytic manifold and
$Z,\,\tilde{Z}$ are closed analytic subsets in $M$. One can then use
the definition of Borel-Moore homology as the homology of the complex
formed by subanalytic chains, cf.  e.g. [GM2] or [KS].  It is known
further, see [RoSa], that the set $Z\cap\tilde{Z}$ has an open
neighborhood $U$ in $M$ such that $Z\cap\tilde{Z}$ is a proper
homotopy retract of $\overline{U}$, the closure of $U$ (this is a
general property of analytic sets).  Now, given two subanalytic cycles
$ c\in H_*(Z)$ and $\tilde{c} \in H_*(\tilde{Z}) $, one can give the
following geometric construction of the class $c\cap \tilde{c} \in
H_*(Z\cap \tilde{Z})$.

First choose $V$, an open neighborhood of $Z$ in  $M$, such that $Z$ is
 a proper homotopy  retract of $\overline{V}$, and  $\overline{V} \cap
 \tilde{Z} \subset U$. Second,  since $V$  is  smooth, one can  find a
 subanalytic  cycle $c'$ in $\overline{V}$  which is homologous to $c$
 in $\overline{V}$  and  such that  the set-theoretic intersection  of
 $c'$  with $\tilde{c}$  is  contained   in  $V$ and, moreover,   $c'$
 intersects $\tilde{c}$ transversely at smooth points of both $c'$ and
 $\tilde{c}$. Hence,  the set-theoretic intersection $c'\cap\tilde{c}$
 gives a  well-defined    subanalytic cycle in   $H_*(\overline{V}\cap
 \tilde{Z})$,  and therefore   in  $H_*(\overline{U})$. Finally,   one
 defines $c\cap\tilde{c} \in H_*(Z\cap \tilde{Z})$ as the direct image
 of $c'\cap\tilde{c}$  under   a proper contraction   $\overline{U}\to
 Z\cap  \tilde{Z}$  which    exists  by  assumption.   It  is   fairly
 straightforward to check that this way one obtains  the same class as
 the  one  defined in  (\ref{qs_pairing})  via  the $\cup$-product  in
 cohomology. It follows in particular that the result of the geometric
 construction above  does not  depend  on the choices  involved in the
 construction.

\section{ Convolution in Borel-Moore homology}
\label{convolution_algebra}

In this section we  give a general  construction of a convolution-type
product in   Borel-Moore homology. Though   looking technically  quite
involved,   the  construction    is   essentially    nothing  but    a
``homology-valued"   version   of   the  standard  definition   of  the
composition of multi-valued maps.
\medskip

\noindent
{\bf Toy example}
\label{toy}

\noindent
We begin with the trivial case of the convolution product.
 We write $\C(M)$ for the finite dimensional
vector space of $\C$-valued
functions on a finite set $M$. Given
finite
sets  $M_1,\,M_2,\,M_3$,  define a
convolution product:
\[
\C(M_1 \times M_2) \otimes \C(M_2 \times M_3) \to \C(M_1 \times
M_3)
\]
by the formula
\begin{equation}
\label{toy1}
f_{12} * f_{23} : (m_1,m_3) \mapsto
\sum_{m_2\in M_2}
f_{12}(m_1,m_2)\cdot f_{23}(m_2,m_3)\,.
\end{equation}
Writing $d_i$ for the cardinality of the finite set $M_i$ we may
naturally identify $\C(M_i \times M_j)$ with the vector space of
$d_i\times d_j$-matrices with complex entries. Then, formula
(\ref{toy1}) turns into the standard formula for the
matrix multiplication.

As a  next step of our toy   example we would like   to find a similar
convolution  construction  assuming that  $M_1,\,M_2,\,M_3$ are smooth
compact manifolds  rather than finite sets  (note that the compactness
condition is a natural generalization of the finiteness condition. The
latter was needed in order to make the sum  in the RHS of (\ref{toy1})
finite).  As one   knows from elementary analysis,  it  is usually the
measures and not the functions that can be convoluted in a natural way.
In differential   geometry the role    of measures is  played  by  the
differential  forms. Thus,  given   a   smooth manifold $M$,  we   let
$\Omega^\bullet(M)$   denote   the    graded    vector     space    of
$C^\infty$-differential forms on $M$. This is the right substitute for
the vector space $\C(M)$ when a finite set is replaced by a manifold.

Let $M_1\,,\, M_2\,,\, M_3$ be smooth compact oriented manifolds, and
${p}_{ij}: M_1\times M_2\times M_3 \to M_i\times M_j$ the projection
to the $(i,j)$-factor. Put $d=\dim M_2$.  We now define a convolution
product:
\[
\Omega^{i}(M_1 \times M_2) \otimes \Omega^{j}(M_2 \times M_3) \to 
\Omega^{i+j-d}(M_1 \times
M_3)
\]
by the formula
\begin{equation}
\label{toy2}
f_{12} * f_{23}=\int_{M_2} p^*_{12}f_{12}\wedge p^*_{23}f_{23}.
\end{equation}
Here $ \int_{M_2} $ stands for operation of integrating over the fibers
of the projection $p_{13}: M_1 \times M_2 \times M_3 \to M_1 \times
M_3$ (see [BtTu]).

The standard properties of differential calculus on manifolds show
that the convolution (\ref{toy2}) is compatible with the De Rham
differential, i.e., we have
\[d(f_{12} * f_{23}) = (df_{12}) * f_{23} + (-1)^j f_{12}
 * (df_{23})\quad,\quad j=\deg f_{12}.\]
It follows that the convolution product of differential forms
induces a convolution product on the De Rham cohomology:
\begin{equation}
\label{toy3}
H^{i}(M_1 \times M_2) \otimes H^{j}(M_2 \times M_3) \to 
H^{i+j-d}(M_1 \times M_3).
\end{equation}
The latter can be transported,
via the
Poincar\'e duality, to a similar convolution in  homology. 

In what follows, we are going to give an
alternative ``abstract" definition of the convolution product
(\ref{toy3}) in terms of algebraic topology. One advantage
of such an ``abstract" definition is that it works for {\it any}
generalized homology theory, e.g., for K-theory. Such a
K-theoretic convolution will be discussed below and applied to
representation theory later. Another
advantage of the ``abstract" definition is that it enables
us to make a refined convolution construction ``with supports".

\medskip

\noindent
{\bf General case}
\label{gener}

\noindent
We proceed now to the ``abstract" construction of the convolution
product.  Let $M_1, M_2,M_3$ be connected, oriented
$C^\infty$-manifolds and let
$$
Z_{12}\,\subset\, M_1\times M_2,\qquad Z_{23}\,\subset\, M_2\times M_3
$$
be closed subsets.
Define the set-theoretic composition $Z_{12}\circ Z_{23}$ 
as follows

\begin{eqnarray}
\label{circ_comp}
Z_{12}\circ Z_{23}
&=&\{(m_1,m_3)\in M_1\times M_3\mid\; 
   \mbox{ there exists }\, m_2\in M_2\\
& &\mbox{ such that } (m_1,m_2)\in Z_{12} \mbox{ and } 
   (m_2,m_3)\in Z_{23}\}.\nonumber
\end{eqnarray}

\smallskip

If we think of $Z_{12}$ (resp. $Z_{23}$) as the graph of a multivalued map
from $M_1$ to $M_2$ (resp. from $M_2$ to $M_3$), then
$Z_{12}\circ Z_{23}$ may
be viewed as the graph of the composition of $Z_{12}$ and $Z_{23}$.

\medskip

\noindent
{\bf Example}\enspace
Let ${f}:M_1\to M_2$ and ${g}:M_2\to M_3$ be smooth maps. Then
$$
\mbox{Graph}(f)\circ \mbox{Graph}(g) = \mbox{Graph}(g\circ f). \quad \Box
$$

\smallskip

We will need another form of definition (\ref{circ_comp})  in the future.
Let ${p}_{ij}: M_1\times M_2\times M_3 \to M_i\times M_j$
be the projection to the $(i,j)$-factor.
From now on, we assume, in addition, that the map
\begin{equation}
\label{some_number}
{p}_{13}: p_{12}^{-1}(Z_{12}) \cap p_{23}^{-1}(Z_{23}) 
\to M_1\times M_3\quad\mbox{\it is proper.}
\end{equation}
We observe that 
\[ p_{12}^{-1}(Z_{12}) \cap p_{23}^{-1}(Z_{23}) 
= (Z_{12}\times M_3)\cap (M_1\times Z_{23})
= Z_{12} \times_{X_2}Z_{23} \,.\]
 Therefore the set
$Z_{12}\circ Z_{23}$ defined in (\ref{circ_comp}) is equal
to the image of the map (\ref{some_number}).
 In particular, this set
is a closed subset in $M_1\times M_3$, since the map in (\ref{some_number})
is proper.

\vspace{2pt}

Let $d=\dim_\R M_2$.
We define a {\it convolution in Borel-Moore homology}, cf. also [FM],
\begin{equation}
\label{bm_conv_map}
H_i(Z_{12})\times H_j(Z_{23}) \to H_{i+j-d}(Z_{12}\circ Z_{23})
\quad,\quad (c_{12},c_{23})\mapsto c_{12}*c_{23}
\end{equation}
by translating the set theoretic composition into composition of cycles.  
Specifically put (compare with (\ref{toy1})):
$$
c_{12}*c_{23}=
(p_{13})_*((p_{12}^*c_{12})\,\cap\,(p_{23}^*c_{23}))
\in H_*(Z_{12} \circ Z_{23}),
$$
where $p_{12}^*(c_{12}):=c_{12}\boxtimes [M_3]$, and
$p_{23}^*(c_{23}):=
[M_1]\boxtimes
c_{23}$ are given by the K\"unneth formula , and the
intersection pairing $\cap$ was defined in (\ref{qs_pairing}). Note that 
$$
((c_{12}\boxtimes [M_3])\cap ([M_1]\boxtimes 
c_{23}))\,\subset\, (Z_{12}\times M_3)\cap (M_1\times Z_{23}),
$$
so that the direct image is well defined due to the condition
that the map $p_{13}$ in (\ref{some_number}) is proper.  The reader
should be warned that although the ambient manifolds $M_i$ are not
explicitly present in (\ref{bm_conv_map}), the convolution map {\it
does} depend on these ambient spaces in an essential way.  Note that
changing the orientation would change the sign of the fundamental classes
$[M_i]$ in the formula, hence it would change the convolution product.

\medskip
\noindent
{\bf Associativity of convolution}\enspace
\label{associativity}
Given a fourth oriented manifold,
$M_4$, and a closed subset $Z_{34}\,\subset\, M_3\times M_4$,
the following associativity equation holds in Borel-Moore homology.
\begin{equation}
\label{assoc}
(c_{12}*c_{23})*c_{34} =
c_{12}*(c_{23}*c_{34}),
\end{equation}
where $c_{12}\in H_*(Z_{12}),\,c_{23}\in H_*(Z_{23}),\,
c_{34}\in H_*(Z_{34}).$ For the proof of the associativity equation see [CG].

\medskip
\noindent
{\bf Remark} \enspace 
The same definition applies in the disconnected case as
well, provided $[M_1]$, resp. $[M_3]$, is understood as the sum of the
fundamental classes of connected components of $M_1$, resp. $M_3$.

\bigskip
\noindent
{\bf Variant: Convolution in equivariant K-theory}

\smallskip
 A similar convolution  construction  works for any  {\it  generalized
homology  theory} that has pull-back  morphisms  for smooth maps, push
forward morphisms  for proper  maps  and an intersection pairing  with
supports.   This is the    case, e.g. for the  topological  K-homology
theory used in [KL1]  and also for  the algebraic  equivariant K-theory
(though the latter is not a generalized homology theory).

Given a complex linear algebraic group $G$ and a complex algebraic
$G$-variety $X$, let $\mbox{Coh}^G(X)$ denote the abelian category of
$G$-equivariant coherent sheaves on $X$.  Let $K^G(X)$ be the
Grothendieck group of $\mbox{Coh}^G(X)$.  Given $\cal
F\in\mbox{Coh}^G(X)$ let $[\cal F]$ denote its class in $K^G(X)$.  For
any $X$, the $K$-group has a natural ${\bold R}(G)$-module structure
where ${\bold R}(G)= K^G(pt)$ is the representation ring of $G$.  We
recall a few properties of equivariant K-theory (see [CG] for more
details).

\vskip2mm

\noindent (a) For any proper map $f\,:\,X\to Y$ between two $G$-varieties 
$X$ and $Y$ there is a direct image $f_*\,:\,K^G(X)\to K^G(Y)$. 
The map $f_*$ is a group homomorphism.

\vskip2mm

\noindent (b)  If $f\,:\,X\to  Y$    is flat  (for instance an    open
embedding) or is a closed embedding of a smooth $G$-variety and $Y$ is
smooth,   there       is an  inverse       image   homomorphism (of
groups) $f^*\,:\,K^G(Y)\to K^G(X)$.

\vskip2mm

Recall the general convolution setup.
Let $M_1$, $M_2$ and $M_3$ be smooth $G$-varieties. Let
$$p_{_{ij}}\,:\,M_1\times  M_2\times M_3   \to M_i\times M_j$$  be the
projection along the factor not  named.  The $G$-action on each factor
induces a natural  $G$-action on the Cartesian  product such  that the
projections    $p_{_{ij}}$ are $G$-equivariant.  Let $Z_{_{12}}\subset
M_1\times M_2$  and  $Z_{_{23}}\subset  M_2\times  M_3$  be $G$-stable
closed subvarieties such that (\ref{some_number}) holds. 
 Define a {\it convolution} map in K-theory
$$*\,:\,K^G(Z_{_{12}})\otimes K^G(Z_{_{23}})\to K^G(Z_{_{12}}\circ
Z_{_{23}})$$ as follows. Let $\cal F_{_{12}}$, $\cal F_{_{23}}$ be two
equivariant   coherent    sheaves on   $Z_{_{12}}$  and
$Z_{_{23}}$, respectively. Set
$$[\cal F_{_{12}}] *[\cal F_{_{23}}]= p_{_{13*}}\Big(p^*_{_{12}}[\cal
F_{_{12}}] \buildrel {L}\over\otimes p^*_{_{23}}[\cal
F_{_{23}}]\Big).$$
In this formula, the upper star stands for the
pullback morphism, well-defined on smooth maps, and $\buildrel
{L}\over\otimes $ is defined by choosing a finite locally free
$G$-equivariant resolution $F_{12}^{\bullet}$ of $p^*_{_{12}} \cal
F_{_{12}}$ (resp. $F_{23}^{\bullet}$ of $p^*_{_{23}} \cal F_{_{23}}$)
on the ambient smooth space $M_1 \times M_2 \times M_3$, and taking
the simple complex associated with the double-complex
$F_{12}^{\bullet} \otimes F_{23}^{\bullet}$.

\bigskip
\noindent
{\bf Examples}

{\rm (i)} Let $M_1=M_2=M_3=M$ be smooth, and 
$$
Z_{12},\; Z_{23}\,\subset\, M_\Delta \,\hookrightarrow\, M\times M,
$$
where $M_\Delta\,\hookrightarrow\, M\times M$ is the diagonal embedding.  If
$Z_{12}$ and $Z_{23}$ are closed then $p_{13}$ in (\ref{some_number})
is always proper, and moreover,
$$
Z_{12}\circ Z_{23}=Z_{12}\cap Z_{23}\,\subset\, M_{\Delta}\,\subset\, 
M\times M.
$$
In  this case we see  that the  $*$-convolution product in homology
reduces  to    the       intersection  $\cap$-product   defined     in
(\ref{qs_pairing}) above, and  $*$-convolution in K-theory reduces  to
the tensor product with supports (see [CG, Corollary 5.2.25]).

\vskip 2mm
{\rm (ii)}
Let $M_1$ be a point and ${f}:M_2\to M_3$ be a proper map of connected
varieties.
Set 
$\;Z_{12}= pt\times M_2= M_2,$ and
$\,Z_{23}=\mbox{Graph}(f)\,.\,$
Then $Z_{12}\circ Z_{23}=\mbox{Im} f\,\subset\, pt \times M_3=M_3$.
Let $c\in H_*(M_2)= H_*(Z_{12})$. Then we have
$c*[\mbox{Graph} f]=f_*(c)$. 

\vskip 2mm
{\rm (iii)} Let $M_3=pt$ and $Z_{23}=M_2\times \,pt$ .
Then the convolution 
\begin{equation}
\label{p2add1}
H_i(Z_{12})\otimes H_j(M_2)\to H_{i+j-d}(M_1)\quad,\quad d=\dim M_2\,.
\end{equation}
allows one to think of $H_*(Z_{12})$ as part of ${\rm Hom} (H_*(M_2),
H_*(M_1))$.

\medskip
\noindent
{\bf The convolution algebra [Gi1]}
\label{c_a}
\smallskip

\noindent
Let $M$ be a smooth complex manifold, let $N$ be a (possibly singular)
variety, and let
${\mu}:M\to N$ be a proper map.
Put $M_1=M_2=M_3=M$ and 
$Z=Z_{12}=Z_{23}=M\times_{_N} M$ in the general convolution setup.  
Explicitly, we have
$$
Z=\{(m_1,m_2)\in M\times M\mid \mu(m_1)=\mu(m_2)\}.
$$
It is obvious that $Z\circ Z=Z$.  Therefore we have the convolution
maps, cf. (\ref{bm_conv_map}),
$$
H_*(Z)\times H_*(Z)\to H_*(Z), \quad \mbox{resp. } K^G(Z) \times K^G(Z)
\to K^G(Z)
$$
in the $G$-equivariant setup.

The following corollary is an immediate consequence of 
(\ref{assoc}).

\medskip

\noindent
{\bf Corollary 2.1} (i) {\it
$H_*(Z)$
has a natural structure of an associative algebra with unit.
Similarly, in the $G$-equivariant setup, $K^G(Z)$
has a natural structure of an associative ${\bold R}(G)$-algebra with
unit.}

(ii) {\it
The unit in $H_*(Z)$, resp. in $K^G(Z)$,
is given by the fundamental class of $M_\Delta\,\subset\, Z$,
resp. by the structure sheaf of $M_\Delta$.
$\square$}

\medskip 

Choose $x\in N$ and set $M_x=\mu^{-1}(x)$.  Apply the convolution
construction for $M_1=M_2=M$ and $M_3$ a point.  Let
$Z=Z_{12}=M\times_{_N} M$ and $Z_{23}=M_x\,\subset\, M\times \{pt\}$.
We see immediately that $Z\circ M_x= M_x$.

\medskip

\noindent
{\bf Corollary 2.2} {\it
$H_*(M_x)$ has a natural structure of a left $H_*(Z)$-module
under the convolution map.} $\Box$

\medskip
\noindent
{\bf Examples}

{\rm (i)} Assume $N=M$ and $\mu: M \to N$ is the identity map. Then $Z
= M_{\Delta}$ is the diagonal in $M \times M$. Then the convolution
algebra $H_*(Z)$ is isomorphic to the cohomology algebra $H^*(M)$,
this is easily derived from (\ref{standard_pd}). In particular, the
convolution algebra $H_*(Z)$ is in this case a graded commutative
local ${\Bbb C}$-algebra. Moreover, for any $x \in N$, $M_x =
\mu^{-1}(x) = \{x\}$, so that $H_*(M_x) \simeq {\Bbb C}$ is the (only)
simple module over this local algebra.

{\rm (ii)} Assume $M$ is smooth and {\it compact}, and $N = pt$, so
that $\mu$ is a constant map. Then $Z = M \times M$ and $M_x = M$,
$\{x\}= N$. Furthermore, the convolution action $H_*(Z) \times
H_*(M_x) \to H_*(M_x)$ can be seen to give an algebra isomorphism
$$
H_*(Z) = {\rm End}_{\,{\Bbb C}}\, H_*(M).
$$
In particular, $H_*(Z)$ is a simple (matrix) algebra.

\medskip

In the general case of an arbitrary morphism $\mu: M \to N$ the
algebra $H_*(Z)$ is, in a sense, a combination of the special cases
(i) and (ii) considered above. In general, the variety $Z$ is the
union of the family $\{Z_x = M_x \times M_x, \, x \in N \}$. The algebra
$H_*(Z)$ is neither simple nor local, and is, in a sense, ``glued''
from the ``family'' of simple algebras $\{{\rm End}_{\,{\Bbb C}}\, H_*(M_x),\,
x \in N \}$. However, these simple algebras are ``glued together'' in a
rather complicated way depending on how far the map $\mu: M \to N$ is
from a locally trivial fibration.

\medskip
\noindent
{\bf The dimension property} 
\label{dimension_property}

\noindent
Let $M_1,M_2,M_3$ be smooth varieties of real dimensions $m_1$, $m_2$,
$m_3$, respectively.  Let $Z_{12}\,\subset\, M_1\times M_2$ and
$Z_{23}\,\subset\, M_2\times M_3$, and let
$$
p=\frac{m_1+m_2}{2},\quad q=\frac{m_2+m_3}{2}, \quad r=\frac{m_1+m_3}{2}\, .
$$
Then it is obvious from (\ref{bm_conv_map}) that convolution induces a map
(assuming that $p,q$ and $r$ are integers)
$$H_p(Z_{12})\times H_q(Z_{23})\to H_r(Z_{12}\circ Z_{23}).$$
We say that this is the property that
``the middle dimension part is always preserved.''

\medskip

Therefore, in our convolution-algebra setup $Z=M \times_N M$, 
and the dimension property yields:

\smallskip

\noindent
{\bf Corollary 2.3 [Gi1]}
\label{cor1_lecture10}
{\it $H(Z)$ is a subalgebra of $H_*(Z)$.}
\medskip

The last result is especially concrete in view of the following
\medskip

\noindent
{\bf Lemma 2.4}
\label{basis_cor}{\it
Let $\{Z_w\}_{w\in W}$ be the irreducible components of $Z$
indexed by a finite index set $W$.  If all the components have the same
dimension then the fundamental classes
$[Z_w]$ form a basis for the convolution algebra $H(Z)$.}
\medskip

\noindent
{\bf Proof} \ 
This follows from Proposition 1.1. $\Box$

\smallskip
In a similar way, one derives from formula (\ref{p2add1}):
\medskip

\noindent
{\bf Corollary 2.5}
\label{cor2_lecture10} {\it
The convolution action of the subalgebra $H(Z)\,\subset\,H_*(Z)$
on $H_*(M_x)$ is degree preserving, i.e,
for any
$i\ge 0$ we have $H(Z)\,*\, H_j(M_x)\,\subset\, H_j(M_x)$.} $\square$

\section{ Constructible complexes}

This  section contains definitions   and  theorems that will allow  us
later  to interpret   the  Borel-Moore homology  and   the convolution
product in  sheaf-theoretic  terms.  

For any topological space $X$ (subject to conditions described at the
beginning of Section 1), let $\cal Sh(X)$ be the abelian category of 
sheaves
of $\C$-vector spaces on $X$. Define the category $\Comp^b(\cal Sh(X))$ as
the category whose objects are {\it finite} complexes of sheaves on $X$
$$
A^{\bullet}=( 0\to A^{-m} \to A^{-m+1} \to\ldots \to A^{n-1} \to
A^{n}\to
 0)
,\qquad m,n \gg 0,
$$
and whose morphisms are morphisms of complexes 
$A^\bullet \to B^{\bullet}$  commuting with the differentials. 
Given a complex of sheaves $A^{\bullet}$ we let 
$$
\cal H^i(A^{\bullet})= \mbox{Ker}(A^i \to A^{i+1})/
\mbox{Im}(A^{i-1} \to A^{i})\,
$$ 
denote the $i$-th cohomology sheaf. A
morphism of complexes is called a {\it quasi-isomorphism} provided it
induces isomorphisms between cohomology sheaves.

The derived category, $D^b(\cal Sh(X))$, is by definition the category
with  the same objects as  \\  $\Comp^b(\cal Sh(X))$ and with morphisms
which  are obtained  from those in   $\Comp^b(\cal Sh(X))$ by  formally
inverting  all quasi-isomorphisms;  thus {quasi-}isomorphisms become
isomorphisms in the  derived category. For  example, we may (and will)
identify $D^b(\cal Sh(pt))$, the derived  category on $X=pt$, with the
derived category of  bounded complexes of  vector spaces.  
In general, the kernels and
cokernels of morphisms are  not  well-defined in $D^b(\cal  Sh(X))$ so
that this category is no longer abelian.  It has instead the structure
of  a {\it triangulated} category.   This structure involves, for each
$n \in \Z$, a  translation  functor ${[n]}:A  \mapsto A[n]$  such that
$\cal H^i(A[n]) =\cal H^{i+n}(A)$, for all $i\in  \Z\,$, and a class
of {\it   distinguished} triangles  that  come  from  all short  exact
sequences of complexes. The precise
definition of the  derived category is  a bit more involved  than this
oversimplified exposition leads  one  to  believe.   For more on   the
derived  category  see [KS], [Iv],  and  [Ha] [Ver2].

The reason for introducing derived categories is that most of the
natural functors on sheaves, like direct and inverse images, are not
generally exact, i.e.  do not take short exact sequences into short
exact sequences. The exactness is preserved, however, provided the
sheaves in the short exact sequences are injective.  Now, the point is
that any sheaf admits an injective resolution (possibly not unique)
and, more generally, any complex of sheaves is quasi-isomorphic to a
complex of injective sheaves. The notion of an ``isomorphism" in
$D^b(\cal Sh(X))$ is defined so as to ensure that any object of
$D^b(\cal Sh(X))$ can be represented by a complex of injective
sheaves.  In this way, all the above-mentioned natural functors become
exact, in a sense, when considered as functors on the derived
category.

>From now on we assume $X$ to be a  complex algebraic variety.  A sheaf
$\cal F$ on $X$ is said to be {\it constructible} if there is a finite
algebraic stratification $\,X= \bigsqcup  X_\alpha\,,$ such that for each
$\alpha$, the stratum    $X_{\alpha}$  is  a  locally closed    smooth
connected algebraic subvariety of $X$, and the restriction of $\cal F$
to the stratum   $\,X_\alpha$ is a  locally-constant  sheaf of finite
dimensional  vector  spaces   (such locally-constant  sheaves  will be
referred to as  {\it local systems} in the  future).  An object $A \in
D^b(\cal Sh(X))$  is said to be  a {\it constructible complex}  if all
the cohomology sheaves  $\cal H^i(A)$ are constructible.  Let $D^b(X)$
be the full subcategory  of $D^b(\cal Sh(X))$ formed  by constructible
complexes ({\it full}  means that the morphisms  remain the same as in
$D^b(\cal Sh(X))$).   The  category $D^b(X)$   is  called the  bounded
derived category of  constructible  complexes on  $X$ in  spite of the
fact that it  is {\it not}  the derived  category  of the  category of
constructible sheaves.

Our next objective is to give a definition of the dualizing
complex and the Verdier duality functor on $D^b(X)$.

Let ${i}:X \hookrightarrow M$
be a closed  embedding of a topological space $X$ into a {\it smooth} 
manifold $M$ (this always exists).
We define a functor 
$$
{i}^!:\cal Sh(M) \to\cal Sh(X),
$$
by taking germs of sections supported on $X$. Specifically, given a sheaf
$\cal F$ on $M$ and an open set $U \,\subset\, M$ set
$$
\Gamma_{[X]}(U,\cal F)=
\{f \in\Gamma(U, \cal F)\,|\, \mbox{supp}(f) \,\subset\, X\cap U\}.
$$ 
The stalk of the sheaf $i^!\cal F$ at a point  $x \in X$ is defined by
the formula
\[(i^!\cal F)_{|x}=\lim_{\rightarrow} \Gamma_{[X]}(U,\cal
  F), \] 
where the direct limit is taken over all open neighborhoods $U
\ni x$.  The functor $i^!$ is  left exact, and we let ${Ri}^!:D^b(\cal
Sh(M))\to D^b(\cal Sh(X))$  denote the  corresponding derived functor.
If $X$ and $M$ are algebraic varieties one  proves that ${Ri}^!$ sends
$D^b(M)$ to $D^b(X)$.

\vskip 4pt

Let $\C_X\in D^b(X)$ be the constant sheaf, regarded
as a complex concentrated in degree zero.
Define the ``dualizing complex'' of $X$, denoted
${\Bbb D}_X$, 
to be 
\begin{equation}
\label{dualizing_def_eqn}
{{\Bbb D}}_X=Ri^!(\C_M)\,[2\dim_\C M],
\end{equation}
where ${i}:X\hookrightarrow M$ as above.

The stalks of the cohomology sheaves of the
dualizing complex are given by the formula 
\begin{equation}
\label{stalk_D}
\cal H^j_x({\Bbb D}_X)= H^{^{j+2dim_\C M}}\!(U\,, U\setminus (U\!\cap\!X)\,)
=H^{BM}_{-j}(U\cap X)
\;  \mbox{ for all } x \in X,
\end{equation}
where $U \,\subset\, M$ is a small contractible open neighborhood of $x$ in $M$,
and the last isomorphism is due to Poincar\'e duality (\ref{pd_defn_bm}).

\medskip

\noindent
{\bf Proposition 3.1} (i) {\em
Let ${i}:N\hookrightarrow M$ be a closed embedding of a 
smooth complex variety $N$ into a smooth complex variety M. Then we have}
$$
Ri^!(\C_M)= \C_N[-2d], \mbox{ where } d = \dim_\C M - \dim_\C N.
$$

(ii) {\em The dualizing complex ${\Bbb D}_X$ does not depend on the
  choice of the embedding ${i}:X\hookrightarrow M$. Moreover, for a
  smooth variety $X$ we have}
$$
{\Bbb D}_X = \C_X[2\dim_\C X].
$$

\noindent
{\bf Proof} \,See Lemma 8.3.3 and Proposition 8.3.4 in [CG]. $\Box$

\medskip

>From now on we will never make use of the functor $i^!$ itself
and will only use the corresponding derived functor. Thus, to simplify
notation we write $i^!$ for $Ri^!$, starting from this moment.

To any object $\cal F \in D^b(X)$ and any  integer $i\in \Z$ we assign
the {\it  hyper-cohomology}  group $H^i(\cal  F)=H^i(X,\cal F)$.  This
is, by definition,  the $i$-th derived functor  to the global sections
functor $\Gamma:\cal Sh(X)\to\{\mbox{\it complex vector    spaces}\}$.
Explicitly, to   compute  the   derived   functors above,     find   a
representative  (up to quasi-isomorphism) of $\cal  F \in D^b(X)$ by a
complex of injective  sheaves $\cal I^\bullet\in \Comp^b(Sh(X))$.  Then
we have by definition of derived functors, see [Bo]:
$$
H^i(\cal F):=H^i(\Gamma(\cal I^\bullet))=
H^i(\mbox{Hom}_{Sh(X)}(\C_X,\cal I^\bullet)).
$$

We list the following basic isomorphisms, which we will use extensively:
\begin{equation}
\label{equ1}
\label{equ2}
H^i(X)= H^i(X,\C_X)\quad, \quad
H_i(X)= H^{-i}(X,{\Bbb D}_X).
\end{equation}
The second isomorphism is a global counterpart of (\ref{stalk_D}). This can
be seen as follows. The complex  ${\Bbb D}_X$ is obtained by applying the
functor $Ri^!$ to the constant sheaf on an ambient smooth variety $M$.
The hyper-cohomology is the derived functor of the functor of global
sections. Thus, $H^\bullet({\Bbb D}_X)$ is equal to the hyper-cohomology
of $R\Gamma_{[X]}$, the derived functor of the functor $\Gamma_{[X]}$ of
global sections supported on $X$. But the hyper-cohomology
of $R\Gamma_{[X]}$, applied to the
constant sheaf on $M$, is clearly $H^\bullet(M, M\setminus X)$, and the
isomorphism follows by Poincar\'e duality.

\vskip 4pt

For any complexes $A, B \in D^b(X),$ one defines Ext-groups in the
derived category as shifted $\mbox{Hom}$'s, that is,
$\mbox{{\rm Ext}}^k_{D^b(X)}(A,B):= \mbox{Hom}_{D^b(X)}(A,B[k])\,.$ There is
also an {\it internal $\cal Hom$-complex}, denoted $\cal Hom(A,B)
\in D^b(X),$ such that the Ext-groups above can be expressed as

\begin{equation}
\label{hom-ext}
\mbox{{\rm Ext}}^\bullet_{D^b(X)}(A,B)=
H^\bullet(X,\,\cal Hom(A,B)).
\end{equation}

We now introduce the
Verdier duality functor, $A\mapsto A^\vee$,
which is a contravariant  functor on
the category $D^b(X)$ defined by the formula
$$
A^\vee = {\cal Hom}(A, {\Bbb D}_X).
$$
Note that with this definition we have
$\C_X^\vee= {\Bbb D}_X$.
It is easy to show that for $\cal F\in D^b(X)$, 
\begin{equation}
\label{calc_shift_def}
(\cal F[n])^\vee = (\cal F^\vee)[-n]\quad\mbox{and}\quad (\cal
F^\vee)^\vee = \cal F\,.
\end{equation}

Given an {\it arbitrary} algebraic 
map ${f}:X_1\to X_2$ we have the following four
functors:
\begin{equation}
\label{four_functors_defn}
{f}_*,{f}_!:D^b(X_1)\to D^b(X_2) \quad,\quad {f}^*,{f}^!:D^b(X_2)\to D^b(X_1).
\end{equation}
The functors $(f_*,f^*)$ are defined as the derived functors of
sheaf-theoretic direct and inverse image functors, respectively.  (We
remark that sometimes what we call $f_*$ is written $Rf_*$ in this
context, but as we will {\it never} use the sheaf theoretic
pushforward we will not adopt the derived functor notation.)  
The other pair $(f_!,f^!)$ is defined via Verdier duality:
\begin{equation}
\label{&&&&&}
{f}_!A_1:=(f_*(A_1^\vee))^\vee ,\qquad f^!A_2:=(f^*(A_2^\vee))^\vee\, ,
\end{equation}
for any $A_1\in D^b(X_1)$ and
$A_2\in D^b(X_2)$.
 There is a direct image formula for hyper-cohomology:
\begin{equation}
\label{hypr_coh_forms}
H^\bullet(X_2,\,f_*A_1)= H^\bullet(X_1,\,A_1)
\end{equation}
and two basic inverse image isomorphisms for ``sheaves''  (see [CG] 
for proofs):
$$
f^*\C_{X_2} = 
\C_{X_1},\qquad f^!{\Bbb D}_{X_2}
= {\Bbb D}_{X_1}.
$$

\vskip 4pt

It is further useful to remember that for a map ${f}:X\to Y$ one has 
\begin{equation}
\label{property}
f_! = f_*\, , \,\;\mbox{\it if $f$ is proper;}
\end{equation}
$$
f^! = f^*[2d], \,\;\mbox{\it if $f$ is flat with smooth fibers of complex dimension $d$}.$$

\vskip 4pt

One should mention that, for a closed embedding ${f}:X_1
\hookrightarrow X_2$, the functor $f^!$ coincides with the derived
functor of the ``sections supported on $X_1$'' functor, which was used
earlier in the definition of a dualizing complex.
\medskip

The  functors (\ref{four_functors_defn}) are related by
a base change formula, see [Ver2]. It says that, given a Cartesian square,
$$
\begin{CD}
X\times_Z Y       @>{\tilde f}>>  Y      \\
@V{\tilde g} VV                  @VgVV   \\
X                 @>f>> Z
\end{CD}
$$
for any object $A\in D^b(X)$, we have a canonical isomorphism:
\begin{equation}
\label{base_chg_shf_vers}
g^! f_* A = {\tilde f}_* {\tilde g}^! A \quad\mbox{\it holds in}\quad D^b(Y).
\end{equation}

Let ${i}_\Delta:X\hookrightarrow X\times X$ be the diagonal embedding. We define
two (derived) tensor product functors on $D^b(X)$ by
\begin{equation}
\label{derived_tensors}
A \otimes B = i_\Delta^*(A\boxtimes B),\qquad
A\stackrel{!}{\otimes}B = i_\Delta^!(A\boxtimes B).
\end{equation}
We will be using later the following canonical isomorphism
 in the derived category:
\begin{equation}
\label{can_1}
{\cal Hom} (A,B)= A^\vee\stackrel{!}{\otimes} B,
\end{equation}
which is a sheaf-theoretic version of the well-known 
isomorphism ${\rm Hom}( V, W) \simeq V^* \otimes W$ for finite dimensional
vector spaces.

\vskip 4pt

Let $N$ be a variety and $A_1, A_2, A_3 \in D^b(N)$. For any $p, q \in
{\Bbb Z}$, the composition of morphisms in the category $D^b(N)$ gives a
bilinear product
$$ {\rm Hom}_{D^b(N)}(A_1, A_2[p]) \times {\rm Hom}_{D^b(N)}(A_2[p],
A_3[p+q]) \to {\rm Hom}_{D^b(N)}(A_1, A_3[p+q]).
$$
Using that $\, {\rm Hom}_{D^b(N)}(A_2[p],$ $A_3[p+q]) =
{\rm Hom}_{D^b(N)}(A_2,A_3[q]) = {\rm Ext}^q_{D^b(N)}(A_2,A_3)\,,$ we can rewrite
the composition above as a bilinear product of Ext-groups, called the
Yoneda product,
\begin{equation}
\label{yon}
\mbox{{\rm Ext}}^p_{_{D^b(N)}}(A_1, A_2)\,\otimes\,\mbox{{\rm Ext}}^q_{_{D^b(N)}}(A_2, A_3)
\to \mbox{{\rm Ext}}^{p+q}_{_{D^b(N)}}(A_1, A_3).
\end{equation}

\section{ Perverse sheaves and the Decomposition Theorem}
\label{perverse sheaves}

We will briefly recall some definitions and  list a few
basic results about
the category of perverse sheaves on a complex
algebraic variety. For a detailed treatment the reader is referred
to [BBD].

A locally constant sheaf $\cal L$ will be refered to as a {\it local system}.
Let $Y\,\subset\, X$ be a smooth locally closed subvariety of complex
dimension $d$, and let
$\cal L$ be a local system on $Y$.
The intersection cohomology complex of Deligne-Goresky-MacPherson,
$IC({Y},\cal L)$, is an object of $D^b(X)$ supported on
$\bar Y$, the closure of $Y$, that satisfies the following properties:
\smallskip

 (a) $\cal H^{i} IC({Y},\cal L)=0 \quad\mbox{if}\quad i<-d$,

 (b) $\cal H^{-d} IC({Y},\cal L)|_Y = \cal L,$

 (c) $\dim\mbox{ supp } \cal H^iIC({Y},\cal L)<-i,\quad \mbox{if $i>-d$}$,

 (d) $\dim\mbox{ supp } \cal H^i(IC({Y},\cal L)^\vee )< -i,\quad \mbox{if $i>-d$}$.
\smallskip

An explicit construction of $IC({Y}, \cal L)$ given in [BBD]  yields
the  following result:
\smallskip

\noindent
{\bf Proposition 4.1} {\em Let ${j}:Y\,\hookrightarrow X$ be an
  embedding of a smooth connected locally closed subvariety of complex
  dimension $d>0$ and ${\bar Y}$ the closure of the image.  Then for
  any local system $\cal L$ on $Y$ there exists a unique object
  $IC({Y},\cal L)\in D^b(X)$ such that the above properties {\rm (a) -
    (d)} hold.  Moreover, one has:

\smallskip{\rm
\begin{tabular}[b]{rl}
{(i)}  & {$\!$\em The cohomology sheaves ${\cal H}^i IC(Y,\cal L)
$ vanish unless $-d\leq i < 0$;}\\
{(ii)} & {$\!$\em ${\cal H}^{-d}IC(Y,\cal L) ={\cal H}^0 (j_*\cal L)$;} \\
{(iii)}& {$\!$\em $IC(Y,\cal L^*) = IC(Y,\cal L)^\vee\,,$
where $\cal L^*$ denotes the local system dual to $\cal L$.}
\end{tabular}
$\Box$
\bigskip

If $X$ is a smooth connected variety, $Y=X$ and $\cal L = \C_X$, then
we have $IC(X,\C_X)= \C_X[\dim_{_\C} X]$.  This motivates the
following definition.  Given a smooth variety $X$ with irreducible
components $X_i$ define a complex $\cal C_X$ on $X$ by the equality
$$
\cal C_X |_{X_i} = \C_{X_i}[\dim_\C X_i].
$$
By Proposition 3.1, the complex $\cal C_X$ is self-dual: 
$\;\cal C_X^\vee = \cal C_X\;$. It
will be referred to as the {\it constant perverse sheaf} on $X$,
for it satisfies the conditions of the following definition.

\smallskip

\noindent
{\bf Definition 4.2}
\label{perverse_definition}
A complex $\cal F\in D^b(X)$ is called {\it perverse sheaf } if
\smallskip

 (a) $\dim \mbox{ supp } \cal H^i\cal F\le -i$,

 (b) $\dim \mbox{ supp }\cal H^i(\cal F^\vee)\le -i,$ for any
 $i$.

\smallskip
 
Observe that the dimension estimates involved in the definition of the
intersection complex $IC({Y},\cal L)$ are similar to properties
(c)-(d) in the definition of a perverse sheaf, except that the strict
inequalities are relaxed to non-strict ones.  Hence, any intersection
complex is a perverse sheaf.  If $\phi$ is a local system on an
unspecified locally closed subvariety of $X$ we will sometimes write
$IC_\phi$ for the corresponding intersection cohomology complex, i.e.
if $\phi$ is a local system on $Y$, then by definition $IC_\phi =
IC({Y},\phi).$

\vskip 4pt

\noindent
{\bf Exercise } Let $X=\C^2$ be the plane with coordinates
$(x_1,x_2)$, and $Y =\{(x_1,x_2) \in \C^2 \,|\linebreak x_1\cdot x_2=0\}$ the
``coordinate cross''.  Check whether the complex $\C_Y[1]$, extended
by $0$ to $\C^2\setminus Y$, is a perverse sheaf on $\C^2$.

\vskip 4pt
\noindent
{\bf Theorem 4.3 [BBD]}
{\em
{\rm (i)} The full subcategory of $D^b(X)$ whose objects are 
perverse sheaves on $X$ is an abelian category, $\Perv(X)$.

{\rm (ii)} The simple  objects of $\Perv(X)$  are the  intersection complexes
$IC({Y},\cal L)$ as  $\cal   L$ runs through the   irreducible locally
constant  sheaves on  various    smooth locally  closed   subvarieties
$Y\,\subset\, X$. $\quad\square$ }

\medskip

\noindent
{\bf Corollary 4.4}\\
\indent{\em
{\rm (a)} There are no {\it negative} degree global
Ext-groups between perverse sheaves,
in particular \\
$\,\mbox{{\rm Ext} }^k_{_{D^b(N)}}(IC_\phi,IC_\psi)=0$ for all $k<0$.

{\rm (b)}  For any irreducible locally constant sheaves $\phi$ and $\psi$ we have
$$
 {\rm Hom}_{_{D^b(N)}}(IC_\phi,IC_\psi)= {\rm Hom}_{_{\Perv(X)}}(IC_{\;\phi},
IC_{\;\psi})= \C\cdot\delta_{\phi,\psi}.
$$}
\indent
Let $X^{\circ}$ be a smooth Zariski open subset in a (possibly
\noindent
singular) algebraic variety $X$.
\medskip

\noindent
{\bf Exercises}

(i) If $A \in \Perv(X)$ then $A^\vee \in \Perv(X)$.

(ii) If $\cal L$ is a local system on $X^{\circ}$ then $IC(X, \cal L)$
has neither subobjects nor quotients in $\Perv(X)$ supported on $X
\setminus X^{\circ}$.

(iii) Deduce from (ii) the following

\medskip

\noindent
{\bf Proposition 4.5 (Perverse Continuation Principle)} 
{\it Any morphism $a: \cal L_1 \to
\cal L_2$ of local systems on $X^{\circ}$ can be uniquely extended to
a morphism $IC(a): IC(X, \cal L_1) \to IC(X, \cal L_2)$, and the map
$a \mapsto IC(a)$ gives an isomorphism} 
$${\rm Hom}(\cal L_1, \cal L_2)
\stackrel{\sim}\to {\rm Hom}(IC(X, \cal L_1), IC(X, \cal L_2)). \quad \Box
$$

\medskip 

We will often be concerned with the homology or cohomology of the
fibers $M_x=\mu^{-1}(x)$ of a proper algebraic morphism $\mu:M\to N$,
where $M$ is a smooth and $N$ is an arbitrary complex algebraic
variety. We first consider the simplest case where $\mu$ is a locally
trivial (in the ordinary Hausdorff topology) topological fibration
with connected base $N$. The (co)homology of the fibers then clearly
form a local system on $N$. In the sheaf-theoretic language, one takes
$\mu_*\C_M$, the derived direct image of the constant sheaf on $M$.
Then the cohomology sheaf $\cal H^j(\mu_*\C_M)$ is locally constant
and its stalk at $x \in N$ equals $H^j(M_x)$.  Replacing $\C_M$ by
${\Bbb D}_M$, the dualizing complex, one sees that the stalk at $x$ of
the local system $\cal H^{-j}(\mu_*{\Bbb D}_M)$ is isomorphic to
$H_j(M_x)$.

Recall now that for any
connected,
locally simply connected topological
space $N$, and a choice of base point $x\in N$, there is an
equivalence of categories
\begin{equation}
\label{rep_equiv_cat}
\left\{{\mbox{local systems}}\atop{\mbox{on $X$}}
\right\}
\leftrightarrow
\left\{{\mbox{representations of the}\atop\mbox{fundamental group $\pi_1(N,x)$}}
\right\}
\end{equation}
sending a  local  system to its fiber   at $x$, which   is naturally a
$\pi_1(N,x)$-module via the mono\-dromy action.   In particular, given a
locally trivial topological fibration $\mu: M\to N$  and a point $x\in
N$, there is a natural $\pi_1(N, x)$-action on $H^\bullet(M_x)$ and on
$H_\bullet(M_x)$, respectively.  We will see below (as a very special,
though not at all trivial, case of the Decomposition Theorem) that this
action  is {\it completely reducible},   that is, both $H^\bullet(M_x)$
and $H_\bullet(M_x)$ are direct sums of irreducible representations of
the group $\pi_1(N, x)$.  For  an irreducible representation $\chi$ of
$\pi_1(N,  x)$,    let  $H_\bullet(M_x)_{\chi}=     {\rm Hom}_{\pi_1(N,
  x)}(\chi,\,H_\bullet(M_x, \C))$  be the  $\chi$-isotypic component
of the homology of  the fiber with {\it  complex}  coefficients. (Up to
now we could work with, say, rational homology. But since some
irreducible representations of $\pi_1(N, x)$ may not be defined over
${\Bbb Q}$ we have to take ${\Bbb C}$ as the ground field from now on.)  This
way we get the direct   sum decompositions into isotypic  components
with respect to the fundamental group

\begin{equation}
\label{fund_action_1}
H^\bullet(M_x, \C)= \bigoplus_{\chi\in \widehat{\pi_1(N,x)}}
\chi\otimes H^\bullet(M_x)_{\chi}\quad,\quad
H_\bullet(M_x, \C)= \bigoplus_{\chi\in \widehat{\pi_1(N,x)}}
\chi\otimes H_\bullet(M_x)_{\chi}\,.
\end{equation}
The first decomposition reflects the corresponding direct
sum decomposition of local systems
\begin{equation}
\label{iso_dec_fnd_act}
\cal H^\bullet(\mu_*\C_M) =\bigoplus_{\chi\in \widehat{\pi_1(N,x)}}
\chi\otimes H^\bullet(M_x)_{\chi}\,,
\end{equation}
where now the LHS stands for the cohomology sheaves; $\chi$ is viewed,
by the correspondence (\ref{rep_equiv_cat}), as an irreducible local
system on $N$, and the vector spaces $H^\bullet(M_x)_{\chi}$ play the
role of multiplicities.  Note that there is no need to write a second
formula of this type, corresponding to homology (as opposed to
cohomology), because on the smooth variety $M$ one has ${\Bbb
  D}_M=\C_M[2\dim_{_{\!\C}}M]$, and the second decomposition is
nothing but the one above shifted by $[2\dim_{_{\C}}M]$.

We recall that a morphism $\mu: M \to N$ is called {\it projective} if
it can  be   factored  as  a  composition  of a  closed   embedding $M
\hookrightarrow {\Bbb P}^n \times N$ and the projection ${\Bbb P}^n \times
N \to N$. Any proper algebraic map between quasi-projective
varieties is known to be projective.
  In  the case of a  projective morphism our analysis will be
based  on  the   very  deep   ``Decomposition Theorem'',  which  has  no
elementary  proof  and is deduced (see   [BBD] and references therein)
from  the Weil  conjectures by  reduction  to ground fields of  finite
characteristic.
\medskip

\noindent
{\bf Decomposition Theorem 4.6 [BBD]}
{\em
\label{decomp_thm}
\label{decomp_eqn_bbd}
Let $\mu:M\to N$ be a projective morphism and $X \,\subset\, M$ a smooth locally
closed subvariety. Then we have a finite direct sum decomposition in $D^b(N)$
\[
\mu_* IC(X,\C_X)= \bigoplus_{(i,Y,\chi)} L_{Y,\chi}(i)\otimes 
IC({Y},\chi)[i],
\]
where $Y$ runs over locally closed subvarieties of $N$,
$\chi$ is an irreducible local system
on $Y$, $[i]$ stands for the shift in the derived category and 
$L_{Y,\chi}(i)$ are certain finite dimensional vector spaces. $\square$
}
\medskip

Now let $M$ be a smooth complex algebraic variety, $\mu:M\to N$ a
projective morphism, and $N=\bigsqcup N_\alpha$ an algebraic
stratification such that, for each $\beta$, the restriction map
$\mu:\mu^{-1}(N_\beta)\to N_\beta$ is a locally trivial {\it
  topological} fibration (such a stratification always exists, see
[Ver1]).  Applying the Decomposition Theorem to $\mu_*\cal C_M$ we see
that all the complexes on the RHS of the decomposition have locally
constant cohomology sheaves along each stratum $N_\beta$. Thus, the
decomposition takes the form
\begin{equation}
\label{decom-star-ng}
\mu_*\cal C_M= \bigoplus_{k\in \Z,\phi=(N_\beta,\chi_\beta)}\;
L_\phi(k)\otimes IC_\phi[k],
\end{equation}
where $IC_\phi$ is the intersection cohomology complex associated
with an irreducible local system $\chi_\beta$ on a stratum
$N_\beta$.

\section{Sheaf-theoretic analysis of the convolution algebra}

    Given  a smooth  complex  variety   $M$ and a  proper map
$\mu:M\to N$, where $N$ is not necessarily smooth, following the setup
of the end of section 2 we  put $Z=M\times_{_N}M$.  Then $Z\circ Z=Z $
so that $H_\bullet(Z)$ has a natural associative algebra structure.

This construction can be ``localized'' with respect to the base $N$
using sheaf-theoretic language as follows.  Consider the constant
perverse sheaf $\cal C_M$ and the complex vector space
$\mbox{{\rm Ext}}^{\bullet}_{_{D^b(N)}}(\mu_* \cal C_M,\,\mu_* \cal C_M)$.
The latter space has a natural (non-commutative) graded ${\Bbb
  C}$-algebra structure given by the Yoneda product of Ext-groups, see
(\ref{yon}).  This Ext-algebra construction is ``local'' in the sense
that one may replace the space $N$ here by any open subset $N'\subset
N$ to obtain a similar Ext-algebra on $N'$.

In the sequel we will often be dealing with linear maps between graded
spaces that do not necessarily respect the gradings. It will be
convenient to introduce the following.

\medskip 
\noindent
{\bf  Notation} Given graded vector  spaces $V, W$, we write $V\doteq
W$ for  a linear isomorphism that  does {\it not} necessarily preserve
the   gradings.  We  will also  use  the  notation $\doteq$ to  denote
quasi-isomorphisms that  only  hold up  to  a  shift in  the   derived
category.
\medskip

We  are going to prove  an  algebra isomorphism $\, H_\bullet(Z)\doteq
\mbox{{\rm Ext}}^\bullet_{_{D^b(N)}}(\mu_* \cal C_M,\, \mu_* \cal C_M).\,$ 
This   important  isomorphism  will allow us    to  study  the algebra
structure of $H_\bullet(Z)$  via the  sheaf-theoretic decomposition of
$\mu_* \cal C_M$.

\medskip

\noindent
{\bf Proposition 5.1} {\em There exists a $($not necessarily grading 
preserving$)$ natural algebra isomorphism 
$$
H_\bullet(Z)\doteq \mbox{{\rm Ext}}^\bullet_{_{D^b(N)}}(\mu_* \cal
C_M,\,\mu_*  \cal C_M).$$
} 

\noindent
{\bf Proof}\enspace\, Since $\cal C_M = {\Bbb C}_M[\dim M]$ we may
replace $\cal C_M$ by ${\Bbb C}_M$ in the statement of the proposition
without affecting the Ext-groups.  Further, we have seen in
(\ref{equ2}) that $H_{\bullet}(Z) \doteq H^{\bullet}(Z, {\Bbb D}_Z)$.
Now use the following Cartesian square:
$$
\begin{CD}
Z = M \times_N M @>{\tilde i}>> M \times M \\
@V{\mu_{\Delta}}VV            @V{\mu \times \mu}VV \\
N_{\Delta} @>i>> N \times N 
\end{CD}
$$
to obtain (denoting $\mu_* {\Bbb C}_M$ by $\cal L)$:
$$\begin{array}{rlcl}
H^{\bullet}(Z, {\Bbb D}_Z) 
&\!\!\!\doteq H^{\bullet}(Z, \tilde{i}^! {\Bbb C}_{M  \times M})
  &&\mbox{since $M \times M$ is smooth}\\
&\!\!\!= H^{\bullet}(N_{\Delta},(\mu_{\Delta})_* \tilde{i}^! {\Bbb C}_{M
    \times M}) 
  &&\mbox{by (\ref{hypr_coh_forms})}\\
&\!\!\!=H^{\bullet}(N_{\Delta},i^!(\mu\times\mu)_*{\Bbb C}_{M \times M})
  &&\mbox{by (\ref{base_chg_shf_vers})}\\
&\!\!\!=H^{\bullet}(N_{\Delta}, i^! (\cal L \boxtimes \cal L)) 
  &&\mbox{by definition of $\cal L$}\\
&\!\!\!=H^{\bullet}(N_{\Delta}, i^! ((\mu_* {\Bbb C}_M^{\vee})^{\vee} \boxtimes \cal L)) 
  &&\mbox{since $\mu$ is proper, and Verdier}\\
& &&\mbox{duality commutes with $\mu_*$}\\
&\!\!\!\doteq H^{\bullet}(N_{\Delta}, i^! (\cal L^{\vee} \boxtimes 
\cal L)) 
  &&\mbox{since $M$ is smooth, hence ${\Bbb C}_M^{\vee}\doteq{\Bbb C}_M$}\\
&\!\!\!= H^{\bullet}(N_{\Delta}, \cal L^{\vee} \stackrel{!}\otimes
\cal L) 
  &&\mbox{by definition of $\stackrel{!}\otimes$}\\
&\!\!\!= {\rm Ext}_{D^b(N)}^{\bullet}(\cal L, \cal L)
  &&\mbox{by (\ref{can_1}) and (\ref{hom-ext})}
\end{array}$$

  This shows that the two spaces are isomorphic  as vector spaces over
  ${\Bbb  C}$. The fact  that this  isomorphism  agrees with the algebra
  structures  is  more complicated;  it  is   proved  in [CG,  Theorem
  8.6.7]. $\Box$

\medskip

Assume from now on that the morphism $\mu: M\to N$ is projective and
that $N=\bigsqcup N_\alpha$ is an algebraic stratification such that,
for each $\beta$, the restriction map $\mu:\mu^{-1}(N_\beta)\to
N_\beta$ is a locally trivial {\it topological} fibration.  We will
study the structure of the convolution algebra $H_\bullet(Z)$ by
combining Proposition 5.1 with the known structure of the complex
$\mu_{*} \cal C_{M}$, provided by the Decomposition Theorem, see
(\ref{decom-star-ng}).  In this way we will be able to find a complete
collection of simple $H_\bullet(Z)$-modules.

By Proposition 5.1 and (\ref{decom-star-ng}) we have 

\begin{eqnarray*}
H_\bullet(Z)
&\doteq& \bigoplus_{k\in {\Bbb Z}}\; \mbox{{\rm Ext}}^k_{_{D^b(N)}}(\mu_{*}  
\cal C_{M},\mu_{*} \cal C_{M})\\
&=&      \bigoplus_{i,j,k\in {\Bbb Z}, \phi,\psi}\mbox{{\rm Hom}}_{_\C}(L_\phi(i),L_\psi(j))\otimes \mbox{{\rm Ext}}^k_{_{D^b(N)}}(IC_\phi[i], IC_\psi[j])\\
&=&      \bigoplus_{i,j,k\in {\Bbb Z}, \phi,\psi}\,\mbox{{\rm Hom}}_{_\C}(L_\phi(i),L_\psi(j))\otimes \mbox{{\rm Ext}}^{k+j-i}_{_{D^b(N)}}(IC_\phi,IC_\psi).
\end{eqnarray*}
Since the summation runs  over all $i,j,k\in {\Bbb Z}$,
the expression in the last line will not
be affected if $k+j-i$ is  replaced by $k$. Thus, we obtain
$$
H_\bullet(Z)\doteq \bigoplus_{i,j,k\in {\Bbb Z},\phi,\psi}\;
\mbox{{\rm Hom}}_{_\C}(L_\phi(i),L_\psi(j))\otimes 
\mbox{{\rm Ext}}_{_{D^b(N)}}^k(IC_\phi,IC_\psi).
$$
Introduce the notation 
$L_\phi= \bigoplus_{i\in{\Bbb Z}}\, L_\phi(i)$. 
Using the vanishing of
$\mbox{{\rm Ext}}^k_{_{D^b(N)}}(IC_\phi,IC_\psi)=0$ 
for all
$k<0$ by Corollary 4.4, one finds
\begin{equation}
\label{H(Z)-structure}
H_\bullet(Z)\doteq \bigoplus_{k\geq 0, \phi,\psi} \;
\mbox{{\rm Hom}}_{_\C}(L_\phi,L_\psi)\otimes 
\mbox{{\rm Ext}}_{_{D^b(N)}}^k(IC_\phi,IC_\psi).
\end{equation}
Observe that the RHS of this formula has an algebra structure,
essentially via the Yoneda product. Moreover, it is clear that
decomposition with respect to $k$, the degree of the Ext-group, puts a
grading on this algebra, which is compatible with the product
structure.

Recall further that
$\mbox{{\rm Hom}}(IC_\phi,IC_\psi)=0$ unless $\phi=\psi$. This yields
\begin{equation}
\label{major_iso_gin} 
H_{\bullet}(Z)= \Bigl(
\bigoplus_{\phi} \;\mbox{{\rm End}}_{\,{\Bbb C}} L_\phi\Bigr)\;\oplus\;
\Bigl(
\bigoplus_{\phi,\psi,k>0} \mbox{{\rm Hom}}_{_\C}(L_\phi,L_\psi)\otimes 
\mbox{{\rm Ext}}_{_{D^b(N)}}^k(IC_\phi,IC_\psi)\Bigr)\;.
\end{equation}
The first sum in this expression is a direct
sum of the matrix algebras $\mbox{End}\; L_\phi$, hence is a 
semisimple subalgebra (as any direct some of matrix algebras). The
second sum is concentrated
in degrees $k > 0$, hence is a nilpotent ideal
$H_{\bullet}(Z)_+\,\subset\,
H_{\bullet}(Z)$. This nilpotent
ideal is the {\it radical} of our algebra, since
\[H_{\bullet}(Z)/H_{\bullet}(Z)_+
\simeq\bigoplus_{\phi} \mbox{End}\;(L_\phi)\,\]
is  a semisimple algebra.
Now, for each $\psi$, the composition
\begin{equation}
\label{conv11}
H_{\bullet}(Z) \twoheadrightarrow H_{\bullet}(Z)/H_{\bullet}(Z)_+
\; = \; \bigoplus_{\phi} \;\mbox{End}\; L_\phi \stackrel{\pi} 
\twoheadrightarrow
\mbox{End}\;{L_\psi}
\end{equation}
(where $\pi$ is projection to the $\psi$-summand) yields an
irreducible representation of the algebra $H_{\bullet}(Z)$ on the 
vector space $L_\psi$.  Since $H_{\bullet}(Z)_+$ is the radical
of our algebra, and all simple modules of the semisimple algebra 
$\bigoplus_\phi\,\mbox{End}\; L_\phi$ are of the form $L_\psi$,
one obtains in this way the following result.

\medskip

\noindent
{\bf  Theorem 5.2}   {\em   The non-zero members  of  the  collection
$\{L_\phi\}$ $($arising from $(\ref{decom-star-ng}))$  form a complete set
\noindent
of the isomorphism classes of simple $H_\bullet(Z)$-modules  }.

\medskip

\bigskip
\noindent
{\bf Special case: Semi-small maps}
\smallskip

In this subsection we fix a {\it smooth} complex algebraic variety $M$
with connected components $M_1, \ldots, M_r$ and assume that 
$\mu: M \to N$ is {\it projective}. 
Given  $x\in N_\alpha$, we put
$M_x=\mu^{-1}(x)$ and $M_{x,k} := M_x \cap M_k\,,\,k=1,\ldots,r.$

\smallskip

\noindent
 {\bf Notation}\enspace
We introduce the following integers:
$$
m_k=\dim_{\raisebox{-0,5mm}{$\scriptstyle\C$}}M_k\enspace\;,\;\enspace
n_\alpha=\dim_{\raisebox{-0,5mm}{$\scriptstyle\C$}}N_\alpha\;\enspace,
\;\enspace
d_{\alpha,k}=\dim_{\raisebox{-0,5mm}{$\scriptstyle\C$}}M_{x,k}\;,\enspace\mbox{for}\;\;x\in N_\alpha
$$
If $M$ is connected we simply write $m= \dim_{\raisebox{-0,5mm}{$\scriptstyle\C$}}M$, and
simplify $d_{\alpha,k}$ to $d_{\alpha}$.  Given a stratum of $N$ and a
local system $\chi$ on this stratum, we will write
$\phi=(N_\phi,\chi_\phi)$ for such a pair, and in this case we use the
notation $n_\phi$ and $d_{\phi,k}$ for the corresponding dimensions.

\smallskip
The following notion is introduced in [GM], cf. also [BM].

\smallskip
\noindent
{\bf Definition 5.3} The morphism $\mu$ is called {\it semi-small}
with respect to the stratification $N=\bigsqcup N_\alpha$
if, for any component $M_k$  we have 
$ n_\alpha + 2d_{\alpha,k}  \le m_k$
for all $\alpha $ such that $N_{\alpha} \subset \mu(M_k)$. If we
always have $ n_\alpha + 2d_{\alpha,k}  = m_k$ we say that $\mu$ is {\it
  strictly semi-small}; and if we have $ n_\alpha + 2d_{\alpha,k} <
m_k$ for all $N_{\alpha}$ that are not dense in an irreducible
component of $N$ we say that $\mu$ is {\it small}.

\vskip 4pt
The results below copy, to a large extent, the results we have already
obtained  before, but in the semi-small case 
 all formulas become ``cleaner", since most
shifts in the derived category disappear.
The following theorem may be regarded as an especially nice
version of the Decomposition Theorem 
and is one of the main reasons to single out the semi-small maps.
Denote $M_i \times_N M_j $ by $Z_{ij}$. Set $H (Z)=
\bigoplus_{ij}H_{m_i + m_j}(Z_{ij})$, where 
$m_i = \dim_{\; {\Bbb C}}M_i$. Thus, $H(Z)$ is the ``middle-dimension''
subalgebra of $H_*(Z)$.

 Given $x \in N_{\alpha}$, let $M_x = \mu^{-1}(x)$. Put
$H (M_x)= \bigoplus_k H_{2d_{\alpha,k}}(M_{x, k})$, the ``top''
homology of $M_x$. 

\medskip
\pagebreak[3]

\noindent
{\bf Theorem 5.4} {\it
{\rm (i)} Let $\cal C_M$ be the constant perverse sheaf on  $M$. If $\mu$ is
  semi-small then $\mu_*   \cal  C_M$ is    perverse  and we  have   a
  decomposition without shifts:
\begin{equation}
\label{sans_shift}
\mu_* \cal C_M = \bigoplus_{\phi=(N_{\phi}, \chi_{\phi})} L_{\phi}
\otimes IC_{\phi}.
\end{equation}
Furthermore, $H (Z)$ is a subalgebra of $H_{\bullet}(Z)$ and one has
algebra isomorphisms:
$$ H(Z) = {\rm Hom} (\mu_* \cal C_M, \mu_* \cal C_M) = \bigoplus_{\phi} 
\mbox{\rm End}\;_{\C}(L_{\phi}). 
$$ 

{\rm (ii)} For any stratum $N_{\alpha}$,
the family of spaces $\{ H (M_x), x \in N_{\alpha}\}$ forms a local 
system on $N_{\alpha}$. If
$L_{\phi}$ is the multiplicity space in $(\ref{sans_shift})$ such
that $N_{\phi} \ni x $ and $\chi_{\phi}$ is the representation of
$\pi_1( N_{\phi}, x)$ associated with $\phi$ then
$$
L_{\phi} = H(M_x)_{\phi}= 
{\rm Hom}_{\pi_1 (N_{\phi}, x)} ( H(M_x), \chi_{\phi})
$$
In   other       words, each   multiplicity    space   $L_{\phi}$   in
$(\ref{sans_shift})$ can  be  obtained by taking  $\chi_{\phi}$-isotypic
component of  the  local system   on $N_{\phi}$  formed by top  degree
Borel-Moore homology of the fibers.

{\rm (iii)} If $\mu$ is small and $N$ is irreducible then $\mu_* \cal
C_M = IC(\mu_* \C_{|_{N_0}})$, where $N_0$ is the dense stratum $($that
is, the decomposition in {\rm (i)} contains only summands coming from
irreducible local systems on the open stratum $N_0)$.}
\smallskip

\noindent
{\bf Proof }

(i) By the Decomposition Theorem it suffices to show that $\mu_* \cal C_M$
is a perverse sheaf, and for this we may assume without loss of
generality that $M$ is connected of complex dimension $m$ and $N =
\mu(M)$. 

First check condition (a) in Definition 4.2. 
Fix any $x \in N$ and write $i_x: \{x \} \hookrightarrow N$ for
the embedding. Then one has
$$
H^ji_x^*(\mu_* \cal C_M)= H^ji_x^*(\mu_* \C_M[m])= H^{j+m}(M_x)
$$ 
Hence if $x \in N_{\alpha}$ then we have a chain of implications
\begin{eqnarray*}
 H^ji^*_x(\mu_* \cal C_M) \neq 0 
&\!\!\Rightarrow&\!\! j+m \leq 2 \dim_{\C}M_x \leq m - \dim_{\C} N_{\alpha} 
           \mbox{ (by definition of semi-smallness) }\\
&\!\!\Rightarrow&\!\! \dim_{\C} N_{\alpha} \leq - j 
\end{eqnarray*}

\noindent and condition (a) of Definition 4.2 follows. Condition (b)
follows automatically form (a) due to self-duality of $\mu_* \cal
C_M$.

  To prove the second part, notice first that $H(Z)$ is a
  subalgebra due to the dimension property of Section 2. We can repeat
  the proof of Proposition 5.1 (this time minding the superscripts)
  to get $H_k(Z_{ij})= {\rm Ext}^{m_i + m_j - k}(\mu_* \cal C_{M_i}, 
\mu_* \cal C_{M_j})$, and this proves $H (Z)= {\rm Hom}(\mu_* \cal C_M,
\mu_* \cal C_M)$ which implies our assertion in view of Corollary 4.4 (b).

(ii) We can assume that $M$ is connected of complex dimension $m$
(since all the objects involved are direct sums over connected
components of $M$).We use the obvious Cartesian square:
$$
\begin{CD}
M_x @>{\tilde i}>> M \\
@V{\mu}VV   @V{\mu}VV  \\ 
\{x\} @>i>> N
\end{CD}
$$

Then
$$\begin{array}{rlcl}
H_k(M_x) 
&\!\!\!= H^{-k}(M_x, {\Bbb D}_{M_x})
  &&\mbox{by (\ref{equ2})}\\
&\!\!\!=H^{-k}(\{x\}, \mu_* {\Bbb D}_{M_x}) = H^{-k}(\{x\}, \mu_*
\tilde{i}^! \cal C_M[m])
  &&\mbox{since $M$ is smooth,}\\
& &&\mbox{\hfill hence ${\Bbb D}_M= \cal C_M[m]$}\\
&\!\!\!= H^{-k}(\{x\}, i^!_x \mu_* \cal C_M[m])
  &&\mbox{by base change}\\
&\!\!\!=H^{m-k}(\{x\}, i^!_x \mu_* \cal C_M).
\end{array}$$

Using this computation, definition of a semi-small map and
(\ref{sans_shift}), we find:
\begin{eqnarray*}
H(M_x) 
&=& H_{2d_{\alpha}}(M_x) = H^{m-2d_{\alpha}}(\{x\}, i^!_x \mu_*
\cal C_M)= H^{n_{\alpha}}(\{x\}, i^!_x \mu_* \cal C_M) \\
&=&H^{n_{\alpha}}\Big(\{x\}, i^!_x \bigoplus_{\phi} L_{\phi} \otimes
IC_{\phi}\Big) = 
\bigoplus_{\phi} L_{\phi} \otimes H^{n_{\alpha}}(\{x\}, i^!_x IC_{\phi})
\end{eqnarray*}

If the closure of $N_{\phi}$ does not contain $N_{\alpha}$ then $x$ is
not contained in the support of $IC_{\phi}$ and $i^!_x IC_{\phi} = 0$.
If $N_{\alpha} \subset \overline{N_{\phi}}$ then we use $i^!_x
IC_{\phi}= (i^*_x( IC_{\phi}^{\vee}))^{\vee}$.  By Proposition 4.1 (i)
we obtain
\begin{eqnarray*}
H^{n_{\alpha}}(\{x\}, i^!_x IC_{\phi})  
&=&H^{n_{\alpha}}(\{x\}, (i^*_x IC ( \cal L_{\phi}^*))^{\vee})  \\
&=& \Big(H^{- n_{\alpha}}(\{x\}, i^*_x IC ( \cal L_{\phi}^*)) \Big)^* 
= \Big( \cal H^{- n_{\alpha}}_x IC ( \cal L_{\phi}^*) \Big)^*\, .
\end{eqnarray*}

If $N_{\alpha} \neq N_{\phi}$ then $ n_{\alpha}= \dim_{\; {\Bbb C}}
N_{\alpha} < \dim_{\; {\Bbb C}} N_{\phi}$. Denote $\dim_{\; {\Bbb C}}
N_{\phi}$ by $d$ and apply property (c) of an $IC$-complex (cf. the
beginning of Section 4) to $i = -n_{\alpha}$. We find
$$
\quad \mbox{dim(supp } \cal H^{- n_{\alpha}}IC ( \cal L_{\phi}^*)) 
< n_{\alpha}$$ 
and since $IC(\cal L_{\phi}^*)$ is locally constant along
$N_{\alpha}$ this means that $\cal H^{- n_{\alpha}}_x IC ( \cal L_{\phi}^*)=0$ 
for all $x \in N_{\alpha}$.

  If $N_{\alpha} = N_{\phi}$ then by property (b) of an $IC$-complex
  we have $ \cal H^{- n_{\alpha}}_x IC ( \cal L_{\phi}^*)
  =  (\cal L_{\phi}^*)_x  = (\cal L_{\phi}^*)_x$, so
$H(M_x)= \bigoplus_{\phi} L_{\phi} \otimes (\cal L_{\phi})_x$, where the
sum runs over all $\phi$ satisfying $N_{\alpha} = N_{\phi}$. Now
claim (ii) follows from Schur's Lemma.
\medskip

(iii) Suppose the decomposition of (i) has a component $L_{\phi}
\otimes IC_{\phi}$ coming from a local system $\cal L_{\phi}$ on a
stratum $N_{\alpha}$ of dimension $n_{\alpha}$ and
$\overline{N_{\alpha}} \neq N $. Then by property (b) of the
intersection homology complex $IC( \cal L_{\phi})$ (cf. the beginning
of Section 4), for any point $x \in N_{\alpha}$, we get
$$
 H^{-n_{\alpha}}i_x^*(L_{\phi} \otimes IC_{\phi}) = ( L_{\phi} \otimes
\cal L_{\phi})_x \neq 0.
$$

 As in the proof of (i) we have:
$$
H^{-n_{\alpha}}i^*_x(\mu_* \cal C_{M_k}) \neq 0 \Rightarrow 
-n_{\alpha} + m_k \leq 2 d_{\alpha, k}\, ,
$$

\noindent
 while the definition of smallness requires ``$>$'' instead of
 the inequality ``$\leq$'' that we just obtained. Contradiction. $\Box$

\setcounter{section}{5}
\section{Representations of Weyl groups}

Fix a   complex semisimple connected Lie  group  $G$ with Lie  algebra
${\frak  g}$, often viewed  as a $G$-module via the  adjoint action.  We
introduce a few standard  objects  associated with a semisimple  group
(see [Bo3] for more details about  the structure of algebraic groups).
Let $B$ be a Borel subgroup, i.e.  a maximal solvable subgroup of $G$,
see [Bo2], and let $T$ be  a maximal torus  contained in $B$.  Let $U$
be the unipotent  radical of $B$  so that $B=T\cdot U$;  in particular
$B$ is connected.  Let $\frak b$, resp.  $\frak  h$, $\frak n$, denote
the Lie algebra of  $B$, resp. $T$,  $U$, so that  $\frak b = \frak  h
\oplus \frak n$. We  also consider the  normalizer $N_G(T)$ of $T$.
The  quotient $W := N_G(T)/T$ is called the {\it Weyl group} of $G$.
It  is known [Bo2], [Se] that
$W$  is a finite  group  generated by  reflections, if viewed as  a  subgroup of
$GL({\frak h})$.  
The main
result of this section is a geometric description of the group algebra
of $W$ as well as a classification of all its irreducible representations.

Let $\cal B$ be the set of all Borel subalgebras in ${\frak g}$.  
By definition, $\cal B$ is the closed subvariety of the
Grassmannian of $(\dim {\frak b})$-dimensional subspaces in ${\frak g}$ formed
by all solvable Lie subalgebras. Hence, $\cal B$ is a projective
variety called {\it flag variety}.  
Recall that all Borel subalgebras are conjugate under the
adjoint action of $G$ and that  $G_{\frak b}$, the isotropy
subgroup of $\frak b$ in $G$, is equal to $B$ by (cf. [Bo2]).  
Thus, the assignment $g \mapsto \mbox{Ad}\,g({\frak b})$ gives a bijection
$$
G/B \simeq \cal B.
$$
Furthermore, the LHS has the natural structure of a smooth algebraic 
$G$-variety (cf. [Bo2]), and the above bijection becomes a $G$-equivariant
isomorphism of algebraic varieties.

Recall that an element $x\in{\frak g}$ is  called {\it nilpotent} if the
operator $\mbox{ad}\; x:{\frak g}\to{\frak g}$  is nilpotent.  This agrees
with the   usual notion  of  nilpotency when   ${\frak g}={\frak  s} {\frak
l}_n(\C)$.  Let  $\cal N$ denote the set  of all nilpotent  elements of
${\frak g}$.  Clearly  $\cal  N$  is a   closed $\mbox{Ad  }   G$-stable
subvariety  of ${\frak g}$. The set  $\cal N$ is also $\C^*$-stable with
respect to dilations, i.e.  $\cal N$ is a cone-variety.

\medskip Set 
$ \tilde{\cal N} := \{(x,{\frak b})\in\cal N\times \cal B
\mid x\in{\frak b}\}.  $ The fiber over a Borel subalgebra ${\frak b}\in
\cal B$ of the second projection $\pi: \tilde{\cal N} \to \cal B$ is
formed by the nilpotent elements of $\frak b$. But it is clear that
the operator $\mbox{ad} x\,,\,x\in {\frak b}$, is nilpotent if and only
if $x$ has no Cartan component in a decomposition ${\frak b}={\frak
h}\oplus{\frak n}$, where ${\frak n}:=[{\frak b},{\frak b}]$ is the nil-radical
of ${\frak b}$. Thus, an element of ${\frak b}$ is nilpotent if and only
if it belongs to ${\frak n}$. It follows that the projection $\pi$ makes
$\tilde{\cal N}$ a vector bundle over $\cal B$ with fiber ${\frak n}$.
Furthermore, since any nilpotent element of ${{\frak g}}$ is $G$-conjugate
into ${\frak n}$, we get a $G$-equivariant vector bundle isomorphism
$$
\tilde{\cal N} \simeq G\times_{_B} {\frak
n}\,\stackrel{\pi}{\longrightarrow}
\,G/B=\cal B,
$$
where $B$, the Borel subgroup of $G$ corresponding to a fixed Borel
subalgebra $\frak b$, acts on the second factor $\frak n = [\frak b,
\frak b]$ by the adjoint action. In particular, $\tilde{\cal N}$ is
a smooth variety, while $\cal N$ itself is always singular at the origin.

\medskip

Identify ${\frak g}\simeq{\frak g}^*$ via the $G$-equivariant isomorphism
given by an invariant  bilinear form on ${\frak g}$, e.g. the Killing form
$(x, y)= \mbox{Tr}(\mbox{ad} x \cdot \mbox{ad} y)\,,$
cf. [Hum], [Se].

\medskip

\noindent
{\bf Lemma 6.1 (cf. e.g. [BoB])}
{\it
There is a natural $G$-equivariant vector bundle isomorphism
$$
\tilde{\cal N} \simeq T^*\cal B\;(= 
 \mbox{\it cotangent bundle on } {\cal B}).
$$}
{\bf Proof }
Let $e= 1\!\cdot\! B / B\,\in G/B$ be the base point.  
We have $T_{e} (G/B) = {\frak g}/\frak b$
and $T^*_{e} (G/B)= ({\frak g}/{\frak b})^* = {\frak b}^{\perp} \subset\,
{\frak g}^*$.
It follows  that, for any $g\in G$,
$$
T^*_{g\cdot e}(G/B)= \mbox{Ad}\,g\,({\frak b}^{\perp}). 
$$
This shows that the vector bundles $T^*(G/B)$ and $G\times_{_B}
{\frak b}^{\perp}$ have the same fibers at each point of $G/B$, hence
are equal as sets. To prove that they are ismorphic as manofolds, one
can refine the argument as follows.

Consider the trivial bundle ${\frak g}_{_{G/B}}= {\frak  g} \times G/B$ on
$G/B$ with fiber   ${\frak g}$.  The  infinitesimal  ${\frak g}$-action on
$G/B$ gives rise to a vector bundle morphism ${\frak g}_{_{G/B}} \,\to\,
T(G/B)$. It is clear that the kernel of this morphism is the subbundle
$\underline{\frak b} \,\subset\,  {\frak  g}_{_{G/B}}$ whose fiber  at a
point  $x\in G/B$ is the isotropy  Lie algebra  ${\frak b}_x\, \subset\,
{\frak  g}$  at $x$. This   gives  an isomorphism $T(G/B)\,\simeq\,{\frak
g}_{_{G/B}}/\underline{\frak  b}$.      Hence,        $T(G/B)\,\simeq\,
G\times_{_B}  ({\frak  g}/\frak b)$. Taking the dual on each side we get
$T^*(G/B) \simeq G \times_B ({\frak g}/ {\frak b})^*$.

 Note that $({\frak g} / {\frak b})^* \simeq {\frak b} ^{\perp} =$
annihilator in ${\frak g}^*$ of the vector subsepace ${\frak b} \subset
{\frak g}$. Under the isomorphism ${\frak g}\simeq{\frak g}^*$,
the annihilator ${\frak b}^{\perp}\,\subset\, {\frak g}^*$ gets identified with
the annihilator of $\frak b$ in ${\frak g}$ with respect to
the invariant form. The latter is equal to $\frak n$, the nil-radical of
$\frak b$.
Thus, $T^*\cal B= G\times_{_B} {\frak n} = \tilde{\cal N}$. $\Box$

\medskip
Define the map $\mu: \tilde{\cal N} \to \cal N$ to be the restriction
to  $\tilde{\cal N} \subset \cal N \times \cal B$ of the first
projection $\cal N \times \cal B \to \cal N$.

\medskip

\noindent
{\bf Theorem 6.2} {\it The map $\mu: T^*\cal B = \tilde{\cal N} \to \cal
N$ is proper and surjective. Moreover, $\cal N$ is irreducible and $\mu$ is a
resolution of singularities for $\cal N$.}
\medskip

\noindent
{\bf Proof } First of all, $\mu$ is surjective since any nilpotent
element of ${\frak g}$ is known to be contained in the nil-radical of a Borel
subalgebra (cf. [Hum]). The surjectivity of $\mu$ implies 
\smallskip

(i) irreducibility of $\cal N$ (since $T^* \cal B$ is irreducible), 

(ii) the dimension bound: $\dim \cal N \leq \dim T^*\cal B=2 \dim \cal B
= 2 \dim \frak n$.
\smallskip

To prove $\dim \cal N \geq 2 \dim \frak n$, recall first that $\cal N$
can be  also defined  as  a common  zero    set of all   $G$-invariant
polynomials on ${\frak g}$ without  constant term (cf.  [CG, Proposition
3.2.5]).  By the Chevalley  Restriction Theorem we have $\C[{\frak g}]^G
\simeq \C[{\frak h}]^W$ (cf.  for example   [CG, Theorem 3.1.38]), hence 
one  has
exactly $\mbox{rk}\,  {\frak g}$  algebraically independent  $G$-invariant
polynomials and therefore $\dim \cal N \geq \dim {\frak g} - \mbox{rk}\;
{\frak g} = 2 \dim {\frak n}$. Thus $T^*\cal B$ and $\cal N$ have the same
dimension.

One can  prove that  the set  of all {\it  regular} nilpotent elements
(that is,  elements  for which the  dimension of  their centralizer is
equal to $\mbox{rk}\; {\frak  g}$) is a single  congugacy class which is
Zariski-open and dense in $\cal N$ (cf.  [CG, Proposition 3.2.10]). By
$G$-equivariance  and   the  dimension equality   the  preimage of any
regular nilpotent element   is a finite  set. To  prove that  $\mu$ is
generically one to one, by $G$-equivariance it  is enough to show that
some particular regular  nilpotent element has  just one point  in its
preimage. If $e_1,  \ldots, e_l \in  {\frak n} $   are the root  vectors
corresponding to positive simple roots  with respect to $\frak b$, then
one  shows ([CG, 3.2])   that $n = e_1 +   \ldots + e_l$  is a regular
nilpotent   which is  contained in  a  unique  Borel subalgebra,
the one containing $e_1 ,   \ldots , e_l$.  This
implies $\#(\mu^{-1}(n))=1$  and therefore a generic point of  $\cal N$
has  exactly one preimage in $T^*\cal  B$, hence $\mu$ is a resolution
of singularities. $\Box$

\medskip
\noindent
{\bf Definition 6.3} The map $\mu:\tilde{\cal N} = T^*\cal B\to \cal
N$ is called the {\it Springer resolution.}

\medskip
\noindent
{\bf Remark} The map $\mu$ is also the moment map with respect to the
canonical Hamiltonian $G$-action on $T^* \cal B$.

\medskip

\noindent
{\bf Theorem 6.4} {\em The Springer resolution is strictly semi-small\/}
 (cf. Definition 5.3).
\smallskip

\noindent
{\bf Proof } See [CG, 3.3]. $\Box$
\smallskip

Now    we can apply    the  machinery of   Section  5  to the Springer
resolution.  Set $Z = \tilde{\cal  N} \times_{\cal  N} \tilde{\cal N}$,
and call  it the {\it Steinberg  variety}.  Consider the convolution
algebra $H(Z)$ associated with $Z$. If $x \in \cal  N$ then it follows
from the definitions that the set
$M_x = \mu^{-1} (x)$ is  formed by pairs $\{(x, {\frak
b}) | {\frak b}\ni x \}$ where $\frak b$ runs  over the subset $\cal B_x
\subset  \cal B$ of $x$-invariant points  of  $\cal B$ (any element $x
\in {\frak g}$ induces a vector field $\xi_x$ on $\cal B$ and $\cal B_x$
is the subvariety of zeros of this vector field).
 
Denote by $G(x)$ the centralizer of $x \in \cal N$ in $G$ and by
$G(x)^{\circ}$ its identity component.  Since $\mu$ is
$G$-equivariant, $G(x)$ acts on $\cal B_x$ and this induces a
$G(x)$-action of $H_*(\cal B_x)$.  Moreover, since the automorphism of
$\cal B_x$ induced by any $g \in G(x)^{\circ}$ is homotopic to the
identity, we conclude that the group $G(x)^{\circ}$ acts trivially on
$H_*(\cal B_x)$.  Hence we obtain an action of the finite group
$A(x):= G(x) / G(X)^{\circ}$ on the homology groups $H_*(\cal B_x)$.
Write $A(x)^\vee$ for the set of isomorphism classes of irreducible
representations of $A(x)$ occurring in the top homology $H_{d(x)}(\cal
B_x)$, $d(x)= \dim_{{\Bbb R}} \cal B_x$, with non-zero multiplicity.

\smallskip

 Our main result concerning representations of Weyl group is the
 following.
\medskip

\noindent
{\bf Theorem (Geometric Construction of $W$)}

(i) {\em There is an algebra isomorphism $ H(Z) \simeq \C[W]$.

{\rm (ii)} The collection $\{H(\cal B_x)_{\phi} \}$, where $(x, \phi)$ runs over
$G$-conjugacy classes of pairs $x \in \cal N$, $\phi \in A(x)^\vee$,
 is a complete set of irreducible representations of $W$.}

\medskip

 The proof of (i) will be given in Section 9. 
 Then (ii) follows from part (i) and Theorem 5.4.
Claim (ii) of the theorem is known as the ``Springer 
construction of Weyl group representations''.

\bigskip

\noindent
{\bf Geometric construction of $W$ and Chern-Mather classes}

\smallskip
There is an interesting connection between the geometric
  construction of the Weyl group $W$ given above and 
 a general construction of Chern classes for singular
  varieties, see [Sab], outlined below.

  Let $X$ be a smooth complex variety of dimension $n$. A ${\Bbb
  C}^*$-stable subvariety of $T^*X$ will be referred to as a {\it
  cone-subvariety}. Let $L_i(T^*X)$ be the group of algebraic cycles
  generated by {\it isotropic} cone-subvarieties in $T^*X$ of
  dimension $\leq i$. (Recall that $T^*X$ has a canonical symplectic
  2-form $\omega$, and a subvariety $\Lambda \subset T^*X$ is called
  isotropic if the pull-back of $\omega$ to $\Lambda$ vanishes. An
  isotropic subvariety of pure dimension $n = {1 \over 2} \dim T^*X$
  is called ``Lagrangian''.) Define the group of {\it Lagrangian
  cone-cycles} as $L(X): = L_n(T^*X)/L_{n-1}(T^*X)$.

\medskip
\noindent
{\bf Example}\enspace
 If $Y \subset X$ is a smooth subvariety, then its
conormal bundle, $T^*_YX$, is a Lagrangian cone-subvariety in $T^*X$.

In general, given a closed (possibly singular) subvariety $Y \subset
X$, write $Y^{reg}$ for the smooth locus of $Y$. Then $Y^{reg}$ is
dense in $Y$, and $T^*_{Y^{reg}}X$ is a locally closed Lagrangian
cone-subvariety in $T^*X$. Let $\Lambda_Y:= \overline{T^*_{Y^{reg}}X}$
be its closure. Then $\Lambda_Y$ is a Lagrangian cone-cycle, hence an
element of $L(X)$. 

This example is in effect typical since one has:

\smallskip

\noindent
{\bf Lemma 6.5 (see [CG, Lemma 1.3.27])} {\it Any irreducible closed
  Lagrangian cone-sub\-va\-rie\-ty in $T^*X$ is of the form $\Lambda_Y$ for
  an appropriate locally closed smooth subvariety $Y \subset X$. }
   $\Box$

\medskip

\noindent
{\bf Corollary 6.6} {\it The group $L(X)$ is spanned by classes of
the form $\Lambda_Y$, $Y \subset X$.} $\Box$

\medskip

Let $Gr_k(TX)$ be the Grassmann bundle on $X$ formed by all
$k$-dimensional subspaces in the tangent bundle $TX$. Given an
irreducible $k$-dimensional (possibly singular) closed subvariety $Y
\subset X$, the tangent spaces to $Y$ at the regular points of $Y$
give rise to a section $\tau:Y^{reg} \to Gr_k(TX)$. The closure
$\widehat{Y}:= \overline{\tau(Y^{reg})} \subset Gr_k(TX)$ of the image
of $\tau$ is called the
{\it Nash resolution} of $Y$. The natural projection $Gr_k(TX) \to X$
restricts to a proper map $p: \widehat{Y} \to Y$, which is an
isomorphism over $Y^{reg}$. The variety
$\widehat{Y}$ is not necessarily smooth, but it carries a natural rank
$k$ vector bundle $\widehat{T}_Y$, the restriction of the tautological
rank $k$ vector bundle on $Gr_k(TX)$. The bundle $\widehat{T}_Y$ plays
the role of the ``tangent bundle of the singular variety $Y$'' because
we have $\widehat{T}_Y|_{p^{-1}(Y^{reg})}= p^*(TY^{reg})$.

 One defines the Chern-Mather class $c_M(Y) \in H_*(Y)$ of the
 singular variety $Y$ as follows. Write $c(\widehat{T}_Y) \in
 H^*(\widehat{Y})$ for the total Chern class of the vector bundle
 $\widehat{T}_Y$. Let $c(\widehat{T}_Y) \cap [\widehat{Y}] \in
 H_*(\widehat{Y})$ be the corresponding class in Borel-Moore
 homology. We set 
$$
c_M(Y) = p_*(c(\widehat{T}_Y) \cap
 [\widehat{Y}]).$$ 
The class in $H_*(Y)$ thus defined is independent
 of the choice of an ambient smooth variety $X$, see [Mac].

 We use Chern-Mather classes and Lemma 6.5 to define a homomorphism
 $c_M : L(X) \to H_*(X)$ by the formula:
\begin{equation}
\label{mama}
c_M: \Lambda_Y \mapsto i_*(c_M(Y)), \mbox{ where } i: Y
\hookrightarrow X.
\end{equation} 
A totally different, but {\it a posteriori} equivalent construction of
homomorphism (\ref{mama}), based on K-theory, is given in [Gi4].

 One can extend the above construction to the bivariant framework
 [FM]. Thus, given two smooth varieties $X_i$,
 $i=1, 2$, one defines a homomorphism
$$
c^{\mbox{\footnotesize biv}}: L(X_1 \times X_2) \to H_*(X_1 \times X_2)
$$
by assigning to $\Lambda_Y$, $Y \subset X_1 \times X_2$, the
relative Chern-Mather class of the fibers of the projection $Y \to
X_2$. The main reason for introducing bivariant Chern-Mather classes
is that they behave nicely with respect to the convolution. In the
special case of a ``push-forward'', see Example (ii) above formula (13);
this has been shown by MacPherson [Mac]. In the general case, given
smooth varieties $X_i$, $i=1, 2,3$, one defines (under appropriate
``properness'' assumptions like in (\ref{some_number})) a convolution
map $L(X_1 \times X_2) \times L(X_2 \times X_3) \to L(X_1 \times
X_3)$. In [Gi2] we proved the following

\medskip

\noindent
{\bf Theorem 6.7} {\it The map} $c^{\mbox{\footnotesize biv}}$ {\it
  commutes with convolution,  i.e. the following diagram commutes:}
$$
\begin{CD}
L(X_1 \times X_2) \otimes L(X_2 \times X_3) @>{\mbox{\footnotesize
    convolution}}>> L(X_1 \times X_3) \\
@V{c^{\mbox{\tiny biv}} \boxtimes c^{\mbox{\tiny biv}}}VV @V{c^{\mbox{\tiny biv}}}VV \\
H_*(X_1 \times X_2) \otimes H_*(X_2 \times X_3) @>{\mbox{\footnotesize
    convolution}}>> H_*(X_1 \times X_3) \quad \Box
\end{CD}
$$

\medskip

Now,   let $X_1 =  X_2 =  X_3 =  \cal B$  be the flag   manifold for a
semisimple group $G$. Since the Steinberg  variety $Z \subset T^*(\cal
B \times \cal B)$ is a  Lagrangian cone-subvariety, see [CG, Corollary
3.3.4], we may view $H(Z)$, the top homology of  $Z$, as a subgroup in
$L(\cal B \times  \cal B)$. Recall  that the convolution product makes
both  $H(Z)$ and $H_*(\cal  B \times \cal  B)$ an associative algebra. 
Thus Theorem 6.7 yields:

\medskip

\noindent
{\bf Corollary   6.8} {\it  The map}  $c^{\mbox{\footnotesize biv}}:
  H(Z)  \to H_*(\cal B  \times \cal B)$ {\it is an  algebra homomorphism.} 
$\Box$

\medskip

Further, the theorem on the geometric construction   of Weyl groups gives an  algebra
isomorphism ${\Bbb C}[W] \simeq H(Z)$.  
Thus, the  Chern-Mather homomorphism of    Corollary 6.8 may  be
thought of as a homomorphism $c^{\mbox{\footnotesize biv}}: {\Bbb C} [W]
\to H_*(\cal B \times \cal B)$.
We will describe the latter map quite explicitly as follows.

Let $G$ act diagonally on $\cal B \times \cal B$. Choose a basepoint
$e \in \cal B$ fixed by the maximal torus $T$. Assign to $w \in W =
N_G(T)/T$ the $G$-diagonal orbit through the point $(\tilde{w} \cdot
e, e)$, where $\tilde{w}$ is a representative of $w$ in $N_G(T)$. This
assignment gives a canonical bijection between $W$ and the set of
$G$-diagonal orbits in $\cal B \times \cal B$. We write $O_w$ for the
orbit corresponding to $w \in W$, and let $\overline{O}_w$ denote its
closure, and $[O_w] \in H_*(\cal B \times \cal B)$ the fundamental
class of $\overline{O}_w$. In particular, for $w=e= $({\rm unit of}
$W$), we have $O_e=\Delta$, the diagonal in $\cal B \times \cal B$.
Recall further, that the Weyl group $W$ acts naturally on $H_*(\cal
B)$. Hence there is a $W \times W$-action on $H_*(\cal B \times \cal
B)$, and for any $w \in W$ we may form the class $(e \boxtimes
w)(\Delta)$.

\medskip

\noindent
{\bf Proposition 6.9} {\it For any simple reflection $s \in W$ we
  have}
$$
c^{\mbox{\footnotesize biv}}(s) = (e \boxtimes s )(\Delta) + [O_s]. \quad \Box
$$

\medskip Recall (cf. [Se]) that the Weyl group is generated by simple
reflections $s_1, \dots, s_l$, where $l=rk\,{\frak g}$, and that to each
element $w\in W$ we can associate its length $l(w)$, equal to the
number of factors in any minimal decomposition $w = s_{i_1}\cdot
\ldots \cdot s_{i_{l(w)}}$ into a product of simple reflections.

Put $n = \dim_{{\Bbb C}} \cal B$. Since
$c^{\mbox{\footnotesize biv}}$ is an algebra map and the simple
reflection generate $W$, we deduce from the proposition:

\medskip

\noindent
{\bf Corollary 6.10} {\it For any $w \in W$ we have}
$$
c^{\mbox{\footnotesize  biv}}(w)   = (e    \boxtimes w)(\Delta) +
\sum_{i=1}^{l(w)}    c^{\mbox{\footnotesize       biv}}_i  (w)\quad,  \quad\quad
c_i^{\mbox{\footnotesize biv}}(w)  \in  H_{2(n+i)}(\cal B  \times \cal
B).
$$
{\it Furthermore} $c^{\mbox{\footnotesize biv}}_{l(w)}(w) = [O_w]$. $\Box$

\medskip

\section{Springer theory for $\cal U({\frak s}{\frak l}_n)$}

We are going to demonstrate in this section that, as a special case of
the 
general  machinery developed above,  one can construct representations
of ${\frak  s}{\frak l}_n(\C)$  and,  maybe, other semisimple algebras  as
well, cf.  [Na]. Many  of the objects we  use here
for studying the  ${\frak
s}{\frak  l}_n(\C)$-case are analogous to  the objects in  the Weyl group
case studied in the previous section.

We fix  an  integer $n\ge 1$  corresponding to  ${\frak s}{\frak l}_n(\C)$
whose representations we want to study.  We also  fix an integer $d\ge
1$ {\it bearing no relation to $n$}.

An  $n$-step 
partial flag $F$ in the vector space $\C^d$ is a sequence of subspaces
 $0=F_0\subset F_1 \subset \dots \subset F_n=\C^d$, where the
inclusions are not necessarily proper.
Write $\cal F$ for the
set of all $n$-step partial flags in $\C^d$.
In the current situation, $\cal F$ will play the role that the
flag variety $\cal B$ played in representations of Weyl groups.
The space $\cal F$ is a smooth compact manifold with connected components
parametrized by all partitions
$$
\mbox{\bf d} = (d_1 + d_2 +\dots + d_n =d),\quad d_i\in {\Bbb Z}_{\ge 0}.
$$
We emphasize that each $d_i$ here may take {\it any} value $0\le d_i \le d$, 
zero in particular.
To the partition $\mbox{\bf d}=(d_1+\dots + d_n)$ we associate
the connected component of $\cal F$ consisting  of flags
\begin{equation}
\label{components_of_f}
\cal F_{\mbox{\bf d}} =\{ F=(0=F_0\subset \dots \subset F_n=\C^d)\mid
  \dim F_i/F_{i-1} = d_i\}.
\end{equation}

Next we introduce an analogue of the  nilpotent variety in the current
situation  to  be  the set  $N=\{x:\C^d\to  \C^d\mid  \; x  \mbox{  is
  linear}, x^n=0\}$.

We are going to define an analogue of the Springer resolution.
Write $M$ for the set of pairs $M=\{(x,F)\in N\times\cal F\mid
x(F_i)\;\subset\; F_{i-1},  i=1,2,\dots,n\}$.  
Note that the requirement $x\in N$ in this formula is
superfluous because $x(F_i)\subset F_{i-1}$,  for all $i$, necessarily implies
that $x^n=0$. The first and second projections give rise
to a natural diagram:
$$
\begin{CD}
 M @>{\pi}>> \cal F \\
@V{\mu}VV \\
N
\end{CD}
$$

The natural action of $\mbox{GL}_d(\C)$ on $\C^d$ gives rise to 
$\mbox{GL}_d(\C)$-actions on $\cal F, N$ and $M$
by conjugation. 
The projections clearly commute with the $\mbox{GL}_d(\C)$-action. 
\smallskip

We have the following description of the cotangent bundle on $\cal F$;
its proof is entirely analogous to  the proof in the  case of the flag
variety (see Lemma 6.1).
\medskip

\noindent
{\bf    Proposition     7.1}     {\em    There   is       a   natural
  $\mbox{GL}_d(\C)$-equivariant vector bundle isomorphism
$$
\;M\simeq T^*\cal F\;
$$
making the map $\pi$ above into the canonical projection $T^*\cal
F\to\cal F$.
$\square$}
\medskip

The decomposition of $\cal F$ into connected components $\cal F_{\bold
  d}$ gives rise to a decomposition of $M$ into connected components
according to $n$-step partitions of $d$:
$$
M= \bigsqcup_{\bold d}\;\, M_{\bold d}, \qquad M_{\bold d}=T^*\cal
F_{\bold d}.
$$ 

The variety $N$ of nilpotent endomorphisms is naturally stratified by
$GL_d ({\Bbb C})$-conjugacy classes, $N = \bigsqcup_{\alpha}
N_{\alpha}$. For any point $x \in N$ the fiber $\mu^{-1}(x)$ consists
of pairs $(x, \cal F)$ satifying $x(F_i) \subset F_{i-1}$, for all
$i$, and may be identified with a subvariety of $\cal F$ that we
denote by $\cal F_x$.
\medskip

\noindent
{\bf Lemma 7.2 [Spa]} {\em  For any $x \in N$, and any $n$-step
  partition $\bold d$, the set $\cal F_x \cap M_{\bold d}$ is a 
connected variety of pure dimension $($that is, each irreducible
component has the same dimension$)$ and
$$
\dim {\Bbb O}_x + 2\cdot
\dim (\cal F_x\cap \cal F_{\bold d})= 2\cdot \dim \cal F_{\bold d};
$$
\noindent
where ${\Bbb O}_x$ denotes the $GL_d({\Bbb C})$ orbit of $x$}. 
\medskip

This   result was  proved   by   Spaltenstein  [Spa] via an   explicit
computation. The connectivity part  of the lemma  can  be proved in  a
more conceptual way  using  Zariski's Main Theorem (see  [Mum]).  This
theorem works  because $N$  (and, more  generally, the  closure of any
nilpotent conjugacy class in ${\frak s}{\frak  l}_d(\C)$) is known to be a
normal variety.

The second claim of the Lemma concerning  dimension just says that the
morphism  $\mu: M   \to N$  is   strictly semi-small.  The  inequality
between the dimensions  $LHS \leq RHS$  (which amounts to saying  that
$\mu$ is {\rm semi-small}) can be  proved by showing that $Z=M \times_N
M$ is a Lagrangian subvariety  of $M \times M$.   Proof of the  strict
equality as well as of the equidimensionality assertion given in [Spa]
exploits some specific features of $SL_n$ in an essential way and will
not be reproduced  here. These assertions  fail  for simple groups  of
types other than $SL_n$. $\Box$

As  before,  we set  $Z= M  \times_N  M$ and  consider the convolution
algebra $H_*(Z)$. By Lemma 7.2, the  map $\mu$ is (strictly) semi-small
and hence  by  Theorem 5.4   (i) the  subspace $H(Z)  \subset  H_*(Z)$
spanned  by the fundamental classes  of  the irreducible components of
$Z$ is a  semisimple  subalgebra.  As at  the  end of Section  2,  the
algebra  $H_*(Z)$ acts on $H_*(\cal   F_x)$ by convolution.  Using the
dimension property  of Section 2 we  deduce that the  subspace $H(\cal
F_x)$  spanned by the classes  of  the irreducible components of $\cal
F_x$ is stable under $H(Z)$-action.

Using  the  general prescription of Theorem  5.4  (iii) we  should now
decompose $H(\cal F_x)$ into isotypic  components with respect to  the
monodromy action and  the multiplicity spaces  will be the irreducible
modules over    $H(Z)$. In our   particular case    the  decomposition
simplifies due to the following lemma:

\medskip
\noindent
{\bf Lemma 7.3}  {\em  The monodromy  action on
  $H(\cal F_x)$ is trivial for any $x \in N$.} $\Box$

\medskip
\noindent
{\bf Corollary 7.4} 
(i) {\em If $x, y \in N$ are $GL_d({\Bbb C})$-conjugate  then the $H(Z)$-modules
$H(\cal F_x)$ and $H(\cal F_y)$ are isomorphic.

{\rm (ii)}
The spaces $\{ H(\cal F_x) \}$, where $x$ runs over the representatives of
$GL_d({\Bbb C})$-conjugacy classes in $N$, 
  form a complete  collection of irreducible $H(Z)$-modules.} $\Box$

\medskip

Now we identify  the convolution algebra  $H(Z)$. Consider the natural
${\frak s} {\frak l}_n$-action on ${\Bbb C}^n$ and the induced action of the
universal enveloping algebra $\cal  U( {\frak s}  {\frak l}_n)$ on  $({\Bbb
C}^n)^{\otimes d}$, the  $d$-th tensor power. Let
$$
I_d= \Ann({\Bbb C}^n)^{\otimes d} \subset  \cal U({\frak s}{\frak l}_n)
$$
be the annihilator of $({\Bbb C}^n)^{\otimes d}$, a two-sided ideal of finite
codimension in $\cal U({\frak s}{\frak l}_n)$.

\bigskip

\noindent
{\bf Theorem (Geometric Construction of $\cal U({\frak s} {\frak l}_n)$)} 
{\em There exists a natural algebra isomorphism} 
$$
\cal U({\frak s} {\frak l}_n)/I_d  \simeq H(Z).
$$

\medskip
 This theorem will be proved in Section 10.
\smallskip

Recall next that the set of finite dimensional irreducible
representations of ${\frak s} {\frak l}_n({\Bbb C})$ is known to be in
bijective correspondence with the set of all dominant weights of
${\frak s} {\frak l}_n({\Bbb C})$. The latter can be viewed as $n$-tuples
of integers $d_1 \geq \ldots \geq d_n$ modulo the ${\Bbb Z}$-action by
simultaneous translation.

 On  the other hand, by  the Geometric Construction of
$\cal  U({\frak s} {\frak
 l}_n)$ , any   simple   $H(Z)$-module gives  rise,
via the projection $\cal U({\frak s} {\frak l}_n)\twoheadrightarrow
 H(Z),$ to  an  irreducible
 representation of  the Lie algebra ${\frak s}  {\frak l}_n(\C)$.  We wish
 to establish  a relationship between $GL_d({\Bbb C})$-conjugacy classes
 in   $N$  parametrizing irreducible   representations  of $H(Z)$, and
 dominant weights       parametrizing        corresponding irreducible
 representations of ${\frak s} {\frak l}_n({\Bbb C})$.

 Let $x\in N$ be a  linear operator in $\C^d$ such  that $x^n=0$.  Put
 formally $x^0 = Id$. Then  there are two distinguished flags attached
 to $x$:
$$
F^{max}(x)=(0=\mbox{ Ker } (x^0)\,\subset\,\mbox{ Ker } (x)\,\subset\, 
\mbox { Ker } (x^2) \,\subset\, \dots \,\subset\, \mbox{ Ker }(x^n) = 
\C^d)\,,
$$
$$
F^{min}(x) = (0=\mbox{ Im }(x^n)\,\subset\, \mbox{ Im }(x^{n-1}) \,
\subset\, \dots \, \subset\, \mbox { Im }(x)\,\subset\, \mbox{ Im }
(x^0)=\C^d)\,.  
$$
Observe that $F^{max}(x),F^{min}(x)\in \cal F_x$.
We assign to each $x\in N$ the $n$-tuple
$$
{\bold d}(x) = (d_1+\dots+d_n=d), \quad\mbox{where}
\quad d_i= \dim \mbox{ Ker }(x^{i})-\dim\mbox{ Ker }(x^{i-1}).
$$
This is the partition associated to the flag $F^{max}(x)$.

\medskip

\noindent
{\bf Lemma 7.5} {\em The $n$-tuple ${\bold d}(x)$ is a dominant weight.}
\medskip

\noindent
{\bf Proof }
For any $i \geq 1$ we have $x(\mbox{Ker } (x^i))
\subset
\mbox{ Ker } (x^{i-1})$. Hence,
the operator $x$ induces, for each $i \geq 1$, a linear map

$$
{ \mbox{ Ker }(x^{i+1}) \over \mbox{ Ker }(x^i)} \hookrightarrow 
{ \mbox{ Ker }(x^i) \over \mbox{ Ker }(x^{i-1})} 
$$

Observe that this map is injective. Whence $d_{i} \geq d_{i+1}$.
The lemma follows. $\; \; \Box$

\medskip
\noindent
{\bf Remark}\enspace
For any flag $F\in \cal F_x$ we have $F^{min}\le F \le F^{max}$
in the sense that $F^{min}_i(x) \;\subset\; F_i \;\subset\; F^{max}_i(x)$,
for each $i=1,2,\ldots, n$.  To see that $F \le F^{max}$ note that, for
any $F=(0=F_0\subset \dots \subset F_n=\C^d)\,\in \cal F_x $ and any
$i=1,2,\ldots, n ,$ one has
$x^i(F_i)\,\subset\,x^{i-1}(F_{i-1})\,$
$\subset\ldots\,\subset x(F_1)=0 .$ Hence, $F_i\,\subset\,{\mbox{ Ker
    } (x^{i})}$.
The other inclusion is proved similarly. $\; \; \Box$

\medskip
Here is the main corollary of geometric construction of $\cal U({\frak s}
{\frak l}_n)$.  It provides a  construction of all irreducible finite 
dimensional representations
of the Lie algebra ${\frak s} {\frak l}_n(\C)$.

\medskip

\pagebreak[3]
\noindent
{\bf Theorem 7.6 (Springer Theorem for $\cal U ({\frak s} {\frak l}_n)$)}
 {\em For any $x\in N$, we have

{\rm (a)} The simple ${\frak s} {\frak l}_n({\Bbb C})$-module
$H(\cal F_x)$  has the highest weight
$${\bold d}(x) = ( d_1 \geq d_2 \geq \dots \geq d_n)\quad,\;
d_i= \dim \ker (x^{i})-
\dim \ker (x^{i-1})$$
In particular, every finite-dimensional irreducible representation of
the
Lie algebra ${\frak s} {\frak l}_n({\Bbb C})$ is of the form $H(\cal F_x)$.

{\rm (b)}  The flags $F^{max}(x)$ and  $F^{min}(x)$ are isolated
points of the fiber $\cal F_x$. The corresponding
fundamental classes $[F^{max}(x)]\in H(\cal F_x)$ and 
$[F^{min}(x)]\in H(\cal F_x)$
are  a highest weight and a lowest weight
vector in $H(\cal F_x)$, respectively.} $\; \; \Box$

\medskip

The fundamental classes of the irreducible components of the fiber
$\cal F_x$ form a distinguished basis in $H(\cal F_x)$.  This
basis is a weight basis with respect to the Cartan subalgebra
of diagonal matrices in ${\frak s} {\frak l}_n(\C)$. It is not known to
the author whether or not this basis coincides with the canonical
basis constructed by Lusztig and Kashiwara, see [KSa].

\medskip
\setcounter{section}{7}

\section{Fourier transform}

In this section we will introduce the main tool used in the proofs of
the geometric constructions of $W$ and $\cal U ({\frak s} {\frak l}_n)$,
the Fourier transform. Given a complex manifold $X$ and a vector
bundle $E \to X$, the Fourier transform sends certain complexes of
sheaves on the total space of $E$ to complexes of sheaves on the total
space of the dual bundle $E^*$.
\noindent
The first step is to define the appropriate categories.

\medskip
  
\noindent
{\bf Monodromic sheaves}
\smallskip

Let $E \to X$ be a complex holomorphic vector bundle. 

\medskip

\noindent
 {\bf Definition 8.1}\enspace%
A sheaf of $\C$-vector spaces on the total
space of the bundle  $E$ is called
 {\em monodromic} if it is locally constant over the orbits of the
 natural $\C^*$-action on $E$. A complex of sheaves is called {\em
   monodromic} if all its cohomology sheaves are monodromic. Let
 $D^b_{mon}(E)$ be the derived category of the category of bounded
 complexes (of $\C$-vector spaces) with monodromic cohomology sheaves.

\medskip

\noindent
 {\bf Definition 8.2}\enspace%
Denote by $\Perv_{mon}(E)$ the full subcategory
 of the category of perverse sheaves, $\Perv(E)$, on the total
space of the bundle $E$ (see Definition 4.2) formed by the 
{\em monodromic} perverse sheaves. 

\medskip

Consider $E$ as a {\em real} vector bundle.  Then the {\it complex}
dual bundle $E^*$ can be identified with the {\em real} dual of $E$
via the pairing $x, \xi\mapsto \langle x, \xi \rangle = Re(\xi(x))\in
{{\Bbb R}}$. Let $I^{\bullet}$ be a complex of injective sheaves on
$E$ which is bounded below, and let $\tau: E \to X$, ${\check\tau}:
E^* \to X$ be the projections.  Given an open subset $U \subset E^*$,
define $U^{\circ} \subset E$, the polar set for $U$, as the set of all
$x \in E$ such that
\smallskip

\begin{tabular}{rl}
(i)& $\tau(x) \in {\check\tau}(U)$, and\\
(ii)& $\langle x, \xi \rangle \ge 0 \mbox{ for all } \xi \in U 
\mbox{ satisfying }  {\check\tau} (\xi)=\tau(x) $.
\end{tabular}
\smallskip

\noindent
Then $U^{\circ}$ is a closed subset of $\tilde U  = \tau^{-1}{\check\tau}(U)$. 
The assignment, see [Br, p.63]:

\begin{equation}
\label{fourier}
  U \mapsto \Gamma_{U^{\circ}} ( \tilde U, I^{\bullet}|_{\tilde U })  
\end{equation}
\noindent
defines a complex of presheaves on $E^*$.
 
\smallskip

 Denote by  $\tilde{\bold F}(I^{\bullet})$ the  sheafification of this
 complex. Using injective  resolutions of mono\-dromic complexes we  can
 extend $\tilde{\bold F}$ to a functor

$$
\tilde{\bold F} : D^b_{mon}(E) \to D^b_{mon}(E^*)
$$  
\noindent
(cf. [Br, Proposition 6.11] for the proof of the fact that the image of a
monodromic sheaf under $\tilde{\bold F}$ is also monodromic). This
definition is rather technical and will never be used here. What is
only important are the properties of the Fourier transform listed in
Proposition 8.3 below.

It is convenient to consider the shifted functor ${\bold F} =
\tilde{\bold F}[r]$, where $r= rk(E)$.

\medskip

\noindent
{\bf  Remark}\enspace%
The  above  definition of    Fourier  transform looks
mysterious; it  is not even   clear why $\bold  F$  should be called a
Fourier transform functor.  In fact, this functor  comes from a functor
on modules   over the  ring  $\cal  D_E$  of algebraic linear
  differential
operators on  $E$  ($\cal  D_E$-modules for   short).   For any  $\cal
D_E$-module $M$  the  sheaf $\tau_*(M)$  is  a sheaf   of modules over
$\tau_*(\cal D_E)$. In this way we can think of a $\cal D_E$-module as a
sheaf on  $X$  having the structure of   a (left)  $\tau_*(\cal D_E)$-
module.

Assuming  $E$ to be trivial (as  it will  be  in applications), we can
construct a natural  isomorphism  $\tau_*(\cal D_E) \simeq  {\check\tau}_*(\cal
D_{E^*})$ which can be written as  ${\bold F}(x_i) = \partial / \partial
\xi_i$ , ${\bold F}(\partial / \partial x_i) = - \xi_i$ in coordinates $
x_1, \ldots, x_r$ (resp.\  $\xi_1, \ldots, \xi_r$) along the fibers of
$E$ (resp.  $E^*$). Recall  that the classical Fourier transform  also
interchanges multiplication by a   coordinate function with  taking  a
partial  derivative.    This isomorphism   establishes  an equivalence
between the category  of $\cal D_E$-modules and $\cal D_{E^*}$-modules
which we also denote by  ${\bold F}$. More generally, 
for  a non-trivial bundle $E$ one
has $\tau_*(\cal   D_E) = (\Lambda^rE) \otimes_{\cal O_X} {\check\tau}_*(\cal
D_{E^*}) \otimes_{\cal O_X} (\Lambda^r E)^{-1} $,  where $r$ is the rank of
$E$ and $\cal O_X$ is the sheaf of regular functions  on $X$ , cf.
[BMV].    Hence if $\cal   M$   is a $\tau_*(\cal  D_E)$-module,  then
$(\Lambda^rE)^{-1}    \otimes_{\cal O_X}  \cal      M$ is a    ${\check\tau}_*(\cal
D_{E^*})$-module.

There   is a  natural  De Rham   functor from  the  category of  $\cal
D$-modules on a complex manifold  $X$ to the category of perverse
sheaves on $X$ which sends  $\cal M$ to  $DR(\cal M) := \Omega_X
\otimes_{\cal D_X} \cal M$, where the sheaf $\Omega_X$ of top degree
holomorphic differential forms on $X$ has a natural {\it right} $\cal
D_X$-module structure.

  One can show that 
the  Fourier  transform  functor ${\bold  F}$  on  perverse sheaves is
obtained by  applying $DR(\cdot)$ to  the  Fourier transform of  $\cal
D$-modules,  that is, the mysterious definition  above was designed so
that one has (cf. [Br, Corollary 7.22])
$$
{\bold F}_{Perv}(DR(\cal M))= DR({\bold F}_{\cal D-\mbox{\scriptsize mod}} \cal M).
$$

One has other  similarities between $\bold  F$  and the classical  Fourier
transform. For example, there is a $\star$-convolution ({\it not}  the convolution
in Borel-Moore homology defined in Section 2) on $D^b_{mon}({{\Bbb C}}^n)$, similar
to the classical convolution  of functions on a group, defined by :
$$
\cal F^{\bullet} \star \cal G^{\bullet} = s_! (\cal F^{\bullet}
\boxtimes \cal G^{\bullet}), \quad \mbox{ where } s: E \times E \to E
\mbox { is the sum map: } s(x, y) = x+y.
$$
\noindent
One has an isomorphism of functors ${\bold F}( \cal F^{\bullet} \star
\cal G^{\bullet}) = {\bold F}(\cal F^{\bullet}) \otimes {\bold F}(\cal
G^{\bullet})$ (cf. [Br, Corollary 6.3]) analogous to the corresponding
classical result, saying that the Fourier transform takes the
convolution of functions into the product of their Fourier transforms.
\bigskip

 We summarize the properties of $\bold F$ in the following
\smallskip

\noindent
{\bf Proposition 8.3}\\
\indent {\em
{\rm (1)} For a monodromic complex $\cal G^{\bullet}$ there exists a natural
isomorphism 
${\bold F} \circ {\bold F} (\cal G^{\bullet}) \simeq a^* \cal G^{\bullet}$
where $a$ is the automorphism of the total space of the vector bundle
$E$ given by multiplication by (-1). 

{\rm (2)} The image under $\bold F$ of a monodromic perverse sheaf  is also a
monodromic perverse    sheaf.

{\rm (3)} $\bold F$ sets up an equivalence of categories
$\Perv_{mon}(E)\stackrel{\sim}{\rightarrow} \Perv_{mon}(E^*)$.

{\rm (4)} Let $ i_V: V \hookrightarrow E $ be a subbundle and $\C_V=
(i_V)_*(\C)$.  Then $ {\bold} F(\C_V) \doteq \C_{V^{\perp}}$, where
$i_{V^{\perp}}: V^{\perp} \hookrightarrow E^*$ is the embedding of the
annihilator $($in $E^*)$ of the subbundle $V$.}
\smallskip

\noindent
{\bf Proof } (1) is proved in [Br, Proposition 6.11] and (2) follows
from (1). (3) is obtained from a similar statement for monodromic
$\cal D$-modules (cf. [Br, Corollary 7.26] see also the remark above).
To prove (4) notice first that if $V=E$, then by definition $\bold F$
is the constant sheaf supported at zero section. The general case
follows from this by functoriality of $\bold F$ with respect to $i_V$
[Br, Theorem 6.1(2)]. $\Box$

\medskip

\noindent
{\bf Direct Image} \enspace%
Given a vector space $E$ and  a complex variety $X$  we form a trivial
vector  bundle $E_X= E  \times X$. Since $E$ can  be viewed as a
vector bundle over a point, we  have two Fourier transforms defined on
$\Perv_{mon}(E_X)$ and on   $\Perv_{mon}(E)$, respectively.

\medskip

\noindent
{\bf Claim 8.4} {\it For a compact algebraic variety $X$ the following diagram commutes:}
\begin{equation}
\label{four_dirim}
\begin{CD}
 \Perv_{mon}(E_X)  @>{{\bold F}_{E_X}}>> \Perv_{mon}(E_X^*)  \\
@V{(pr_E)_*}VV          @VV{(pr_{E^*})_*}V \\
 \Perv_{mon}(E) @>>{{\bold F}_{E}}> \Perv_{mon}(E^*)
\end{CD}
\end{equation}

\noindent
{\bf  Proof } This follows from [Br, Proposition 6.8]. $\Box$

\section{Proof of the geometric construction of $W$}

Recall that we  want to construct  an algebra isomorphism $H(Z) \simeq
\C[W]$ where $Z$, the Steinberg variety,  arises from the Springer map
 $\mu: \tilde{\cal   N}  \to \cal N$   via  the basic  construction of
 \S 2.

\medskip

\noindent
{\bf The Fourier transform reduction}
\smallskip

We have seen in Proposition 5.1 that $H(Z) \subset {\rm Ext}^{\bullet}(\mu_*
{\cal C}_{\tilde{\cal  N}},  \mu_*{\cal C}_{\tilde{\cal  N}}) $.  
Since   $\mu$ is
semi-small, $\mu_*\cal C_{\tilde{\cal N}}$ is  perverse, and therefore
$H(Z)=\mbox{End}_{Perv} (\mu_*  \cal C_{\tilde{\cal N}}),$  by
Theorem 5.4 (i).

Using the  embedding  $\cal N   \hookrightarrow {\frak g}$   we can view
perverse sheaves on $\cal N$ as perverse sheaves on ${\frak g}$. In this
way  $H(Z)$ is  realized  as an endomorphism  algebra  of the perverse
sheaf $\mu_*  \cal C_{\tilde N}$ on ${\frak  g}$.   The group ${\Bbb C}^*$
acts both on $\cal N$, since a multiple of a nilpotent element is also
nilpotent, and on $\tilde{\cal N}$  (along the fibers of  $\tilde{\cal
N}  \to \cal B$).  Since $\mu$  is  $\C^*$-equivariant, all the direct
image  sheaves   involved  are ${\Bbb C}^*$-equivariant,  in  particular
monodromic.  Therefore we   can apply the Fourier  transform  functor
$\bold  F$ to $\mu_* \cal  C_{\tilde{\cal N}}$.   Then ${\bold F} (\mu_*
\cal C_{\tilde{\cal N}})$ is a perverse sheaf on ${\frak g}^*$ and since
$\bold F$ is an equivalence of categories, $H(Z) \simeq {\rm Hom}_{D^b({\frak
g}^*)}({\bold  F}(\mu_*\cal   C_{\tilde{\cal  N}}),   {\bold F}  (\mu_*\cal
C_{\tilde{\cal N}}))$.

  Recall the definition of the Springer variety:
$$
\tilde{\cal N}= \{ (x, {\frak b}) \in {{\frak g}} \times \cal B | x \in
{\frak n}_{\frak b}
\}
$$ 
\noindent(here ${\frak n}_{\frak b}$ denotes the nil-radical of $\frak
b$). It is clear that $\tilde{\cal N}$ is a subbundle of the trivial
bundle ${\frak g}_{_{\cal B}}= {\frak g} \times \cal B \to \cal B$:
$$
\begin{CD}
\tilde{\cal N} @>>>  {\frak g} \times \cal B  \\
@V{\mu}VV                          @V{\mbox{\footnotesize pr}_{{\frak g}}}VV \\
 \cal N        @>>>  {\frak g}
\end{CD}
$$

To compute $\tilde{{\bold F}}(\cal C_{\tilde{\cal N}})$ recall that
the fiber of $\tilde{\cal N}$ at a point ${\frak b} \in \cal B$ is
isomorphic to the nil-radical $\frak n$ of $\frak b$. If we identify
${\frak g} = {\frak g}^*$ via the Killing form, then ${\frak
  b}^{\perp}= {\frak n}$, and hence ${\bold F}_{{\frak g}}(\cal
C_{\frak n}) \doteq \cal C_{\frak b}$ where ${\frak b} \subset {\frak
  g} \simeq {\frak g}^*$.  By Proposition 8.3(4) the Fourier transform
of $\cal C_{\tilde{\cal N}}$ is isomorphic (up to shift) to the
constant perverse sheaf supported on the subbundle $\tilde{{\frak g}}
= \{(x, {\frak b}) \in {\frak g} \times \cal B\;|$ $ x \in {\frak b}
\} \subset {\frak g}_{_{\cal B}}$.

We will denote the restriction of $\mbox{\it pr}_{_{{\frak g}^*}}\!\!:
{\frak g}^* \times \cal B \to {\frak g}^*$ to $\tilde{{\frak g}}$ by
$\tilde{\mu}$.  The main reason for passing from $\mu: \tilde{\cal N}
\to {\frak g}$ to $\tilde{\mu}: \tilde{{\frak g}} \to {\frak g}^*
\simeq {\frak g}$, a sort of dual map, is in replacing the semi-small
map $\mu$ by a small map {\rm (cf. Definition 5.3)}, as follows from
the result below:

\medskip

\noindent
{\bf Proposition 9.1}\enspace%
{\it The map $\tilde{\mu}: \tilde{\frak g} \to {\frak g} $ is small.}
$\Box$ \medskip

Denote by ${\frak g}^{rs}$ the open subset of ${\frak g}$ of all
regular semisimple elements.  Then the restriction of
$\tilde{\mu}_*\cal C_{\tilde{{\frak g}}}$ to ${\frak g}^{rs}$ is a
local system, $\tilde{\mu}_*\cal C_{\tilde{{\frak g}}}|_{{\frak
    g}^{rs}}$, and we write $IC(\tilde{\mu}_*\cal C_{\tilde{{\frak
      g}}}|_{{\frak g}^{rs}})$ for the corresponding IC-complex on
${\frak g}$, as defined in Section 4.
\medskip
 
\noindent
{\bf Corollary 9.2}\enspace
 {\em $\tilde{\mu}_* \cal C_{\tilde{{\frak g}}}  =
  IC(\tilde{\mu}_*\cal C_{\tilde{{\frak g}}}|_{{\frak  g}^{rs}})$}.

\smallskip
\noindent
{\bf Proof } Follows from Proposition 9.1 and Theorem 5.4 (iii). $\Box$

\medskip 

The basic idea of the Fourier transform reduction can be now summarized
 as follows:
\begin{quote}
  The convolution algebra $H(Z)$ is isomorphic to $ {\rm End} (\mu_*
  \cal C_{\tilde{\cal N}})$.  View $\cal C_{\tilde{\cal N}}$ as a
  perverse sheaf on ${\frak g}_{_{\cal B}}$ extended by zero from
  $\tilde{\cal N}$.  Then \pagebreak[3]
\begin{eqnarray}
\label{strochka}
H(Z) &\simeq &{\rm End}\big(\mu_* \cal C_{\tilde{\cal N}} \big) \nonumber\\
&\simeq& {\rm End}\big({\bold F}(\mu_* \cal C_{\tilde{\cal N}})\big)
  \simeq {\rm End} \; \tilde{\mu}_* \cal C_{\tilde{{\frak g}}} \\
&\simeq& {\rm End}\big(IC(\tilde{\mu}_*\cal C_{\tilde{{\frak g}}}|_{{\frak g}^{rs}} ) \big) \simeq
{\rm End}\big(\tilde{\mu}_*\cal C_{\tilde{{\frak g}}}|_{{\frak  g}^{rs}} \big),
\nonumber
\end{eqnarray}

\smallskip
\noindent
where $pr_{{\frak g}^*}: {\frak g}_{_{\cal B}}^*= {\frak g}^* \times
\cal B \to {\frak g}^* $ is the projection and ${\bold F}( \cal
C_{\tilde{\cal N}})$ is a perverse sheaf on ${\frak g}_{_{\cal B}}^*$.
The fourth isomorphism follows from Corollary 9.2 and the last one
from the Perverse
\noindent
Continuation Principle 4.5.
\end{quote}

\bigskip
\noindent
{\bf Analysis of $\tilde{\mu}|_{\tilde{{\frak g}^{rs}}}$}
\smallskip

We denote $\tilde{\mu}^{-1}({\frak g}^{rs})$ by $\tilde{{\frak
    g}}^{rs}$.  By definition of $\tilde{\mu}$ the fiber of
$\tilde{\mu}: \tilde{\frak g}^{rs} \to {\frak g}^{rs}$ at a point $x
\in {\frak g}^{rs}$ is the set of all Borel subalgebras that contain
$x$.

 Recall that for a fixed maximal torus $T \subset G$ the Weyl group
 $W$ is defined as the quotient $N_G(T)/T$.

\medskip

\noindent
{\bf Proposition 9.3}

{\em
  
{\rm (i)} There exists a free action of $W$ on $\tilde{\frak
    g}^{rs}$ such that $\tilde{{\frak g}}^{rs}/W$ is isomorphic to
  ${\frak g}^{rs}$. Moreover, under this isomorphism $\tilde{\mu}$
  corresponds to the quotient map $\tilde{{\frak g}}^{rs} \to
  \tilde{{\frak g}}^{rs}/W$.
  
  {\rm (ii)} The map $\tilde{\mu}:\tilde{{\frak g}}^{rs} \to {\frak
  g}^{rs}$ is a regular $($i.e. Galois$)$ covering with automorphism
  group~$W$.}

\noindent
{\bf Proof } Choose and fix a Cartan subalgebra $\frak h$
(corresponding to a maximal torus $T$) and a Borel subalgebra $\frak
b$ (corresponding to a Borel subgroup $B$) containing $\frak h$.
Let ${\frak h}^{rs}= {\frak g}^{rs} \cap \frak h$ be the set of regular
elements in $\frak h$.
Define a map $\phi: G/T \times {\frak h}^{rs} \to {\frak g}^{rs}$ by
$\phi(g, h)= (\mbox{Ad}(g) {\frak b}, \mbox{Ad}(g) h)$. Since $T$ acts
trivially on $\frak h$ and maps $\frak b$ into itself, $\phi(g, h)$
depends only on the image of $g$ in $G/T$.

We claim that $\phi$ is an isomorphism. To show that $\phi$ is
injective assume that for $g_1, g_2 \in G$ we have $\mbox{Ad}(g_1)
{\frak b} = \mbox{Ad}(g_2) \frak b$ and $\mbox{Ad}(g_1) h =
\mbox{Ad}(g_2) h$. The first equality implies $g_2^{-1}g_1 \in N_G(B)$
and since $h$ is regular semisimple the second equality implies
$g_2^{-1}g_1 \in N_G(T)$. By [Hum] $N_G(B)= B$ and $B \cap N_G(T)= T$,
so $g_2^{-1}g_1 \in T$, hence $g_1$ and $g_2$ represent the same point
in $G/T$.
 
To prove that $\phi$ is surjective assume that a point $({\frak b}',
x') \in \tilde{{\frak g}}^{rs}$ is given (i.e. $x' \in {\frak b}'$ and
$x'$ is regular semisimple).  Since all Borel subalgebras are
conjugate (cf. [Hum]) we have ${\frak b}'= \mbox{Ad}(g')\frak b$ for
some $g' \in G$ and $B'=\mbox{Ad}(g')B$. Then $\mbox{Ad}(g')\frak h$
is a Cartan subalgebra of $\frak b$ hence by [Hum] there exists an
element $u \in B'= \mbox{Ad}(g')B$ such that $\mbox{Ad}(u) x' \in
\mbox{Ad}(g')\frak h$.  Denote $\mbox{Ad}((g')^{-1}u) x' \in {\frak
  h}^{rs}$ by $x$ and $u^{-1}g'$ by $g$. Then $\mbox{Ad}(g){\frak b} =
\mbox{Ad}(u^{-1})\mbox{Ad}(g') {\frak b} = \mbox{Ad}(u^{-1}) {\frak
  b}' = \frak b'$ and $\mbox{Ad}(g)x =
\mbox{Ad}(u^{-1})\mbox{Ad}(g')\mbox{Ad}((g')^{-1}u) x' = x'$, so
$\phi(g,x) = ({\frak b}', x')$ and $\phi$ is surjective.

  The Weyl group $W=N_G(T)/T$ acts on $G/T \times {\frak h}^{rs}$ by 
 $w \cdot (g, h)= (gn^{-1}, \mbox{Ad}(n)h)$ where $n\in N_G(T)$ is
any representative of $w \in W$.

  It follows from the definitions of $\phi$ and $\tilde{\mu}$ that
  $\tau= \tilde{\mu} \circ \phi: G/T \times {\frak h}^{rs} \to {\frak
  g}^{rs}$ is the quotient map for the action of $W$ on $ G/T \times
{\frak h}^{rs}$. 

To prove the second statement of the proposition one has to show (by
general theory of coverings) that $\tilde{{\frak g}}^{rs}$ is
connected. Since $G/T$ is connected it suffices to show that ${\frak
h}^{rs}$ is connected. But ${\frak h}^{rs}$ is a complement of finitely many
complex hyperplanes in a complex vector space ${\frak h}$ (cf. [Hum]),
therefore it is connected. $\Box$

\medskip
\noindent
{\bf Remark} In fact, the action of $W$ on $\tilde{{\frak g}}^{rs}$ does
not depend on the choice of $T$ and $B$ if one views $W$ as an
abstract Weyl group associated with the root system of ${\frak g}$, cf.
[CG, Chapter~3].

\medskip

Next, for the covering $\tilde{{\frak g}}^{rs} \to {\frak g}^{rs}$,
one has a decomposition
$$
\tilde{\mu}_*\cal C_{\tilde{{\frak g}}}|_{{\frak g}^{rs}}=
\bigoplus_{\psi} L_{\psi} \otimes \cal L_{\psi}
$$
of the local system $\tilde{\mu}_*\cal C_{\tilde{{\frak
      g}}}|_{{\frak g}^{rs}}$ into a direct sum of irreducible
pairwise distinct local systems $\cal L_{\psi}$ (with multiplicity
spaces $L_{\psi}$). Therefore we have by the Perverse Continuation
principle
\begin{eqnarray*}
{\rm End}(\mu_* \cal C_{\tilde{{\frak g}}^{rs}}) 
&=& {\rm Hom} \Big( \bigoplus_{\psi} L_{\psi} \otimes \cal L_{\psi}, 
\bigoplus_{\psi'} L_{\psi'} \otimes \cal L_{\psi'} \Big) \\
&=& \bigoplus_{\psi, \psi'}
{\rm Hom} ( L_{\psi}, L_{\psi'}) \otimes {\rm Hom} ( \cal L_{\psi}, \cal
L_{\psi'}) = \bigoplus_{\psi} {\rm End}_{\; {\Bbb C}}( L_{\psi})  .
\end{eqnarray*}

\noindent
{\bf Corollary 9.4} {\em  There exists a natural algebra isomorphism
$$
H(Z) \simeq \bigoplus_{\psi} {\rm End}_{\; {\Bbb C}}( L_{\psi}). $$}

\noindent
{\bf Proof } Follows from (\ref{strochka}) and the computation above.
$\Box$

\medskip

\noindent
{\bf End of proof of the  geometric construction of $W$} \enspace
 By Corollary
9.4 we just have to compute  the multiplicity spaces $L_{\psi}$ above. 
By the   correspondence  between monodromy representations   and local
systems  (see   (\ref{rep_equiv_cat})) it  suffices  to decompose  the
monodromy  representation   of   $\tilde{\mu}_*  {\Bbb  C}_{\tilde{{\frak
    g}}^{rs}}$    into   irreducible representations.  The   latter, by
Proposition 9.3(ii), is nothing but the regular representation of $W$:
$$
{\Bbb C}[W]= \bigoplus_{\chi \in W^{\vee}} L_{\chi} \otimes \chi\, ,
$$
where the summation is over the set $W^\vee$ of isomorphism classes of
irreducible representations of $W$,  and $L_{\chi}$ is isomorphic as a
{\it vector space} to the dual  of $\chi$. Moreover, the decomposition
of {\it vector spaces}:
$$
{\Bbb C}[W]= \bigoplus_{\chi \in W^{\vee}} \mbox{End}_{\, {\Bbb C}}(L_{\chi})
$$
is an isomorphism of {\it algebras}, hence the theorem  follows by
Corollary 9.4. $\Box$

\medskip

\noindent
{\bf Remark} \enspace
 Applying the inverse Fourier transform to the decomposition
$$
\tilde{\mu}_* \cal C_{\tilde{{\frak g}}} = \bigoplus L_{\psi} \otimes
IC_{\psi}
$$
we get $\,\mu_* \cal C_{\tilde{\cal N}} = \bigoplus_{\psi}\,
L_{\psi} \otimes {\bold F}^{-1}(IC_{\psi})$. Since $\bold F$ is an
equivalence of categories, ${\bold F}^{-1}(IC_{\psi})$ is a simple
perverse sheaf, hence of the form $IC(\cal L)$ for some irreducible
local system $\cal L$. Thus we have proved (!) the Decomposition
Theorem in this case.

\bigskip

\noindent
{\bf Digression: the braid group ${\bold B}_W$}
\smallskip

  The proof of Proposition 9.3 allows us to analyze the
  fundamental group $\pi_1(\tilde{{\frak g}}^{rs})$. Firstly, one has an
  exact sequence 
$$
1 \to \pi_1(G/T \times {\frak h}^{rs}) \to \pi_1({\frak g}^{rs})
\to W \to 1.
$$

  The homogeneous space $G/T$ is naturally a fibration over $G/B$ with
  fibers isomorphic (non-canonically) to $B/T$. Since $\cal B = G/B$
  is simply-connected and $B/T$ is contractible, the exact sequence
  above turns into

\begin{equation}
\label{braid}
1 \to \pi_1({\frak h}^{rs}) \to \pi_1({\frak g}^{rs}) \to W \to 1.
\end{equation}

The  map $\pi_1({\frak h}^{rs})  \to \pi_1({\frak g}^{rs})$  is induced by
the restriction $\tau|_{e \times {\frak h}^{rs}}: {\frak h}^{rs} \to {\frak
g}^{rs}$  where $\tau: G/T  \times {\frak  h}^{rs}  \to {\frak g}^{rs}$  is
defined at the end of the  proof of Proposition 9.3.  This restriction
coincides with  the  natural embedding  ${\frak h}^{rs}  \hookrightarrow
{\frak g}^{rs}$.

\medskip

\noindent
{\bf Definition 9.5} The fundamental group $\pi_1({\frak g}^{rs})$ is
called the {\it braid group} of ${\frak g}$.

\medskip

 Recall    (cf.  [Se]) that the   Weyl  group  is  generated by simple
 reflections $s_1, \ldots,  s_l$ where $l  = \mbox{rk}\,  {\frak g}$ and
 $\{s_i\}_{i=1, \ldots,\, l}$ satisfy the following relations:
\vspace{-2mm}
\begin{enumerate}
\item[(i)] $s_i^2= e$,\vspace{-3mm}
\item[(ii)] For any $i \neq j$ one has $s_is_js_i \dots = s_j
  s_i s_j \dots$ \quad with $m_{ij}$ terms on each side, where
  $m_{ij}$ is defined in the standard way from the corresponding
  Dynkin diagram (cf. [Se]).
\end{enumerate}

The relations (ii) are called Coxeter relations.  \smallskip

There are two different (equivalent) definitions of the braid group
associated to $W$:

\medskip

{\it $1^{st}$ description:} the braid group is an abstract group,
${\bold B}_W$, generated by the elements $T_1, \ldots, T_l$, subject
to the Coxeter relations (ii) associated with $W$ (but with relations
(i) omitted).

{\it $2^{nd}$ description:} the braid group is the group generated by the
elements $T_w, w \in W$, satisfying $T_w T_{w'}= T_{ww'}$ whenever
$l(w)+l(w')= l(ww')$. 

In particular, $W$ is embedded {\it as a set} (not as a subgroup) into
${\bold B}_W$  via $ w  \mapsto T_w$.   One also  has a surjective group
homomorphism ${\bold B}_W  \twoheadrightarrow W$.  Its kernel is  called
the ``colored braid group''.

\medskip

 The second description is ``better'' in some categorical sense
 explained in [D].

\bigskip

\section{ Proof of the geometric construction of $\cal U ({\frak s} {\frak l}_n)$}

The proof\footnote{The idea of using the Fourier transform to prove
  the geometric construction of $\cal U ({\frak s} {\frak l}_n)$ is
  due to A.Braverman and D.  Gaitsgory.} of the geometric construction
of $\cal U ({\frak s} {\frak l}_n)$ follows the same pattern as the
proof of geometric construction of $W$. First notice that we have a
diagram
$$
\begin{CD}
M @>>> {\frak g} {\frak l}_d \times \cal F\\
@V{\mu}VV  @V{\mu}VV \\
N  @>>> {\frak g} {\frak l}_d
\end{CD}
$$

The convolution algebra $H(Z)$ is isomorphic  to $ {\rm Hom}(\mu_* {\Bbb C}_M,
\mu_* {\Bbb C}_M)$ (by Lemma 7.2 and Theorem 5.4).

We can view ${\Bbb C}_M$ as a sheaf on the total space of the bundle
${\frak g} {\frak l}_d \times \cal F \to \cal F $ supported on the
subbundle $M \hookrightarrow {\frak g} {\frak l}_d \times \cal F $ and
$\mu_* {\Bbb C}_M$ as a sheaf on ${\frak g} {\frak l}_d$. Applying the
Fourier transform as in the proof of the geometric construction of $W$,
we obtain by Proposition 8.3:
$$
H(Z) \simeq {\rm Hom} ( \mbox{\it pr}_* {\Bbb C}_{M^{\perp}},
\mbox{\it pr}_* {\Bbb C}_{M^{\perp}}),
$$ 
where $\mbox{\it pr}$ is the projection $({\frak g} {\frak l}_d)^* \times
\cal F \to ({\frak g} {\frak l}_d)^*$ and $M^{\perp}$ is the subbundle of
the trivial bundle $({\frak g} {\frak l}_d)^* \times \cal F \to \cal F$
 annihilating $M$.
 
 Let $\tilde{\mu}: M^{\perp} \to {\frak g} {\frak l}_d$ be the restriction
 of $pr: {\frak g} {\frak l}_d^* \times \cal F \to {\frak g} {\frak l}_d^*$ to
 $M^{\perp}$.   If we identify  $({\frak g} {\frak  l}_d)^*$ with ${\frak g}
 {\frak l}_d$ via the form $\langle A, B \rangle = Tr(AB)$ then, for any
 flag  $F \in \cal  F$, the   fiber  of the  vector bundle  $M^{\perp}
 \hookrightarrow \cal F$  over $F$ is the subspace  $\{ x  \in {\frak g}
 {\frak l}_d  \; |  \; x(F_i) \subset   F_i \quad \forall i=  1, \ldots,
 n\}$.      Thus for  $x      \in {\frak g}  {\frak   l}_d$   the  fiber
 $\tilde{\mu}^{-1}(x)$ over $x$  is equal to  the set of all flags $F=
 \big(  0=F_0 \subset \dots \subset F_n=  {\Bbb C}^d \big)$ such that $
 x(F_i) \subset F_{i}$, for all $i$. We see that
$$
\tilde{\mu}^{-1}(x)=\{ \mbox{$n$-step partial flags in } {\Bbb C}^d \mbox{
  fixed by } x \}.
$$  
 Notice that this fiber is non-empty for {\it any} $x \in
{\frak g} {\frak l}_d$ as is clear from  the Jordan normal form of
$x$. Again as in Section 9 one proves:

\medskip
\noindent
{\bf Proposition 10.1} {\it The map $\tilde{\mu}: M^{\perp} \to {\frak g}
  {\frak l}_d$ is small {\rm (cf. Definition 5.3)}.}$\Box$

\medskip 

First we analyze the fibers of $\tilde{\mu}$ over a regular semisimple
element $x \in {\frak g} {\frak l}^{rs}_d$, that is, a diagonalizable $d
\times d$ matrix with pairwise distinct eigenvalues.
Choose the basis  $(e_1, \ldots, e_d)$ of  ${\Bbb C}^d$ in which $x$  is
  diagonal. Then any subspace  $F_i$ satisfying $x(F_i)\subset F_i$ is
  spanned    by  a subset    of  our basis.    Hence   we can describe
  $\tilde{\mu}^{-1}(x)$ as the set of maps $\phi:  \{1 \ldots d \} \to
  \{ 1  \ldots n \}   $, where $\phi(i)$  is defined  for each  $i= 1,
  \ldots,  d$   as  the {\it   minimal}  number  such  that  $e_i  \in
  F_{\phi(i)}$.

  Therefore the stalk $\big( \tilde{\mu}_* {\Bbb C}_{M^{\perp}}
  \big)_x$ of the local system $\tilde{\mu}_* {\Bbb C}_{M^{\perp}}
  |_{{\frak g} {\frak l}_d^{rs}}$ is isomorphic to the complex vector
  space with base $\{ \phi \;| \;\phi \in \mbox{Maps}(\{1,
  \ldots, d\}, \{1, \ldots, n\}) \}$. This space can be
  identified with $({\Bbb C}^n)^{\otimes d}$ by choosing a basis $f_1,
  \ldots, f_n$ of ${\Bbb C}^n$ and mapping $\phi$ to $
  f_{\phi(1)}\otimes \ldots \otimes f_{\phi(d)} \in ({\Bbb
    C}^n)^{\otimes d}$.

\medskip

\noindent
{\bf   Lemma  10.2} {\em For   $x  \in {\frak  g}  {\frak  l}_d^{rs}$, the
monodromy action  of $\pi_1(  {\frak g}  {\frak  l}_d^{rs}, x)$  on $({\Bbb
C}^n)^{\otimes d}$ factors through   the natural representation of  the
symmetric group $S_d$ on $({\Bbb C}^n)^{\otimes d}$.}
\smallskip

\noindent
{\bf Proof } First notice that Proposition 9.3 is in fact valid for
any reductive Lie algebra ${\frak g}$ (not necessarily semisimple). In
particular we can take ${\frak g} = {\frak g} {\frak l}_d$. By
(\ref{braid}) the monodromy action factors through $W=S_d$ if and only
if the restriction of the local system on ${\frak h}^{rs}$ is trivial.
But in our case, ${\frak h}^{rs}$ is the space of diagonal $(d \times
d)$-matrices with pairwise distinct eigenvalues.  Hence the fibers of
$\tilde{\mu}$ are canonically identified with each other, and the
local system $\tilde{\mu}_* {\Bbb C}$ is constant on ${\frak h}^{rs}$.

 To show that the action of $S_d$ on $({\Bbb C}^n)^{\otimes d}$
 coincides with the natural one, notice that a lift of an element 
  $w \in W$ to a loop in ${\frak g}^{rs}$ is the image (under $\tau: G/T
 \times {\frak h}^{rs} \to {\frak g}^{rs}$ defined in the proof of
 Proposition 9.3) of
  some path in $\psi: [0, 1] \to (G/T) \times {\frak h}^{rs}$ connecting
  $(n^{-1}, w
  \cdot h)$ with $(e, h)$ (where $n \in N_G(T)$ is any preimage
  of $w \in W$). Let $\psi(t)=(g(t), h(t))$. Then the fiber of
  $\tilde{\mu}$ over $\tau(g(t), h(t)$ as a subset of $\cal F$ is
  isomorphic of $g(t) \cdot \tilde{\mu}^{-1}(h)$. Since
  $\tilde{\mu}^{-1}(w \cdot h)$ is canonically identified with
  $\tilde{\mu}^{-1}(h)$ and $g(1)= n^{-1}$, our assertion follows. $\Box$

\medskip

 As in  Corollary 9.4 we have a decomposition of local systems 
on ${\frak g} {\frak l}_d^{rs}$,
$$
\tilde{\mu}_* {\Bbb C}_{M^{\perp}} = \bigoplus_{\psi} L_{\psi} \otimes
\cal L_{\psi},
$$
into a direct sum of irreducible local systems with multiplicities
(the vector spaces $L_{\psi}$). Using the
connection between monodromy and local systems (\ref{rep_equiv_cat})
and Lemma 10.2, we can reformulate the conclusion of Corollary 9.4 as
an algebra isomorphism
$$
H(Z)\simeq {\rm End}_{S_d}(({\Bbb C}^n)^{\otimes d}).
$$

\medskip 
Recall that $I_d= \Ann({\Bbb C}^n)^{\otimes d} \subset  \cal U({\frak
  s}{\frak l}_n)$, see \S7.
The following Lemma is a classical result of H.Weyl which is at
the origin of Schur-Weyl duality.
\medskip

\noindent
{\bf Lemma 10.3} {\it The image of natural homomorphism $\cal U({\frak
    s} {\frak l}_n) \to {\rm End}_{\; {\Bbb C}}(({\Bbb C}^n)^{\otimes
    d})$ commutes with the $S_d$-action and induces an algebra
  isomorphism
$$
{\rm End}_{S_d}(({\Bbb C}^n)^{\otimes d}) \simeq \cal U({\frak s}
{\frak l}_n)/ I_d. \quad \Box
$$
}
\indent Lemma 10.3 completes the proof of the geometric construction
of $\cal U({\frak s} {\frak l}_n)$. $\Box$

\medskip
\section{$q$-Deformations: Hecke algebras and a quantum group}

In this section we will state generalizations of the geometric
constructions of $W$ and $\cal U({\frak s} {\frak l}_n)$ that will
give geometric interpretations to ``quantized'' versions of these
algebras.

First, we want to recall the notations of Section 6 and introduce some
more.  Let $G$ be a  complex semisimple simply-connected Lie  group.
Choose and fix a Borel subgroup $B \subset G$ and let $T$ be a maximal
torus contained in $B$. The  Lie algebra $\frak h$  of $T$ acts on the
Lie algebra ${\frak g}$ of $G$ via the adjoint representation. Denote by
$\Delta \subset {\frak  h}^*$ the set of all   roots of $G$.   For  any
root $\alpha  \in \Delta$  one has  a  reflection $s_{\alpha} \in  W$.
Write $Q\subset \frak h^*$ and  $Q^\vee\subset  \frak h$ for the  root
and coroot lattices, respectively, and  let $X^{*}(T)={\rm Hom}_{alg}(T, \C^*)$ and
$X_*(T)  =  {\rm Hom}_{alg}(\C^*,   T)$   be  the corresponding    weight and
coweight lattices  (where ${\rm Hom}_{alg}$    stands  for
`algebraic  group  homomorphisms') . We have
\[Q\subset X^{*}(T) \quad\mbox{and}\quad
Q^\vee\subset X_*(T). \]

The group $\Pi= X^{*}(T)/Q$ is finite and is known to be isomorphic to
the fundamental group of  $G$.  We also  denote by $\theta^\vee  $ the
maximal coroot of $G$.

\medskip

\noindent
{\bf Definition 11.1} 

(a) The {\it affine Weyl  group $W^a$} of $G$ is  defined as the group
of affine transformations  of $\frak h^*$ generated  by $W$ and $s_0$,
an additional  reflection with respect  to  the affine  hyperplane $\{
h\in {\frak h}^*: \langle\theta^\vee, h\rangle+ 1=0\}$.  Thus $W^a$ is a
Coxeter group with generators $s_i, i=0, ..., l$ (see Definition 9.6).
It is known that $W^a$ is a semidirect product  of $W$ and the co-root
lattice $Q$.

(b) The {\it  extended  affine Weyl  group $\widetilde W$}  of  $G$ is
defined as a semidirect  product of $W$  and $X^{*}(T)$ and is not, in
general, a Coxeter group.  It is clear  that $W^a\subset \widetilde W$
is a normal  subgroup and we have $\widetilde  W/W^a = \Pi$.  For $\pi
\in \Pi$ and $i  \in \{0, \ldots, l\}$,  the transformation
$\pi s_i \pi^{-1}:  \frak h^* \to \frak  h^*$  is a  simple reflection
 again,  which we denote by  $s_{\pi(i)}$.  Thus, $\widetilde W \simeq
 \Pi \ltimes W^a$ with the  commutation relations $\pi s_i \pi^{-1}  =
s_{\pi(i)} $.
 
(c) For $w \in W^a$ we define the {\it length} function $l(w)$ as we
did above Corollary 6.10, and we extend it to $\widetilde W$ by
requiring $l(\pi \cdot w)= l(w) $ for all $\pi \in \Pi, w \in W^a$.

\medskip

The {\it affine Hecke algebra $\bold H$} of $G$ can be defined in two
different (but equivalent) ways. They are analogous to the two
different definitions of the braid group.  \medskip

\noindent
{\bf First definition 11.2} [KL2] The algebra $\bold H$ is the free
${\Bbb Z}[{\bold q}, {\bold q}^{-1}]$-algebra with basis $T_w, \, w\in
\widetilde W$ and multiplication given by the rules:
\begin{equation}
\label{sq_rel}
(T_w+ 1)(T_w - {\bold q}) = 0, \quad w\in \{s_0, ..., s_l\},
\end{equation}
\begin{equation}
\label{len_rel}
T_wT_y = T_{wy}, \quad \mbox{if}\quad l(wy) = l(w) + l(y).
\end{equation}

Denote by  ${\bold H}^a \subset {\bold H}$  the  subspace spanned by
the $T_w$, $w\in W^a$ only. This is clearly a subalgebra, and ${\bold
H}   \simeq {\bold H}^a[\Pi]$ is  the  twisted group  algebra for the
$\Pi$-action   on ${\bold H}^a$.  In  other  words,  $\bold H$  is
generated by the sets $\{T_w, w\in  W^a\}$ and $\{T_\pi, \pi\in \Pi\}$
with the relations   (\ref{sq_rel}), (\ref{len_rel}) for the  $T_w$'s,
and the relations
$$T_\pi T_{\pi'} = T_{\pi+\pi'}; \,\,\,
T_\pi T_{s_i}T_\pi = T_{s_{\pi(i)}}, \quad
i = 0, \ldots, l .$$

\noindent
{\bf  Second  definition 11.3}  The algebra $\bold  H$ is the ${\Bbb
Z}[ {\bold q}, {\bold q}^{-1}]$-algebra with  generators
$\{T_w,  w\in W\}$ and $\{Y_\lambda, \lambda   \in X^{*}(T)\}$ subject
 to the relations:
\smallskip

(i) The $T_w, w\in W$, satisfy (\ref{sq_rel}) and (\ref{len_rel}).

(ii) $Y_\lambda Y_\mu = Y_{\lambda+\mu}.$

(iii) $ T_{s_i} Y_{\lambda} - Y_{s_i(\lambda)}T_{s_i} = 
(1- {\bold q}){Y_{s_i(\lambda)} - Y_{\lambda} \over 1 - Y_{-\alpha_i} }$,
\quad $i = 1, \dots, n.$

\smallskip

It is known that the elements $T_w Y_\lambda$, $w\in W$, $\lambda\in X^{*}(T)$,
form a ${\C}$-basis of $\bold H$.

\smallskip

\noindent
{\bf Remarks} 

(1) The elements $\{Y_{\lambda}, \lambda \in X^*(T) \}$ span a
commutative subalgebra in $\bold H$ isomorphic to ${\Bbb Z}[{\bold q},
{\bold q}^{-1}] \otimes_{{\Bbb Z}} {\Bbb Z} [ X^*(T) ]$, the group
algebra of the lattice $X^*(T)$ over ${\Bbb Z}[{\bold q}, {\bold
  q}^{-1}]$. The latter is also isomorphic to ${\bold R}(T)[{\bold q},
{\bold q}^{-1}]$.

(2) We have $\widetilde  W\simeq W \ltimes X^{*}(T)$. Accordingly, one
 has a presentation of ${\bold H}$ as $H_W \otimes_{{\Bbb Z}} {\bold
 R}(T)[ {\bold q}, {\bold q}^{-1}]$, where $H_W$ is the subalgebra generated
 by the $T_w$ (and the $\otimes$ above is only tensor product as ${\Bbb
 Z}$-modules, not algebras).
 
 (3) Given a Coxeter group $C$, like $W$ or $W^a$, or its close cousin
 like  $\widetilde{W}$,  write  $H_C$   for  the corresponding   Hecke
 algebra.  Thus  ${\bold H}  =   H_{\widetilde{W}}$, and  $ H_W$  is the
 ``finite''  Hecke algebra. Recall further  that  associated with $C$ is
 the corresponding Braid group ${\bold  B}_C$, see Definition 9.5. It is
 clear    from (\ref{len_rel}) that the    Hecke  algebra $H_C$ is the
 quotient of the   group algebra 
  ${\Bbb Z}   [ {\bold q},  {\bold  q}^{-1}][{\bold B}_C]$
 modulo quadratic relations of type (\ref{sq_rel}).

(4) Hecke algebras arise naturally in Lie theory in at least three
different contexts. First of all, the effect of taking the quotient of
${\Bbb Z} [ {\bold q}, {\bold q}^{-1}][{\bold B}_C]
$ modulo (\ref{sq_rel}) is to get an
algebra of the same ``size'' as the group algebra of the group
$C$. More formally, the Hecke algebra $H_C$ is flat over ${\Bbb Z}[
{\bold q}, {\bold q}^{-1}]$ and its specialization at ${\bold q} = 1$ is the
group algebra ${\Bbb Z} [ C]$. Thus, $H_C$ may be thought of as a
``q-deformation'' of the group algebra.

Recall further  that  the Braid group  ${\bold  B}_W$ has  a topological
interpretation as the fundamental   group  ${\bold B}_W  =   \pi_1({\frak
g}^{rs})$. Thus, the Hecke algebra $H_W$ may be viewed as a quotient of
${\Bbb  Z}[{\bold q}, {\bold  q}^{-1}][\pi_1({\frak  g}^{rs})]$ modulo certain
quadratic  relations. Replacing here the  Lie algebra ${\frak g}$ by the
group $G$ we  get a similar interpretation  of the  {\it affine} Hecke
algebra. Specifically,  write  $\overline{W}:= W  \ltimes X_*(T)$  for  the
semidirect   product  of   $W$ with the      coweight lattice, and  let
$H_{\overline{W}}$ be the corresponding Hecke  algebra (note that the group
$\widetilde{W}$  was defined as   the semidirect product  of $W$  with
$X^*(T)$, the    dual  lattice. Thus,  replacing  $\widetilde{W}$   by
$\overline{W}$  amounts  to  replacing  the  root system  by   its dual, or
equivalently,  replacing the group   $G$ by  its {\it Langlands  dual}
$\ ^LG$). Then  one proves, repeating  the argument in  Proposition 9.3,
that
$$
\pi_1(G^{rs}) \simeq {\bold B}_{\overline{W}}.
$$
The group ${\bold B}_{\overline{W}}$ appears on the RHS instead of
${\bold B}_W$ because, in the group case, the Cartan subalgebra $\frak
h$ gets replaced by the corresponding maximal torus $T$ which has
non-trivial fundamental group $\pi_1(T) \simeq X_*(T)$. This
fundamental group is responsible for the extra generators in ${\bold
  B}_{\overline{W}}$ as compared to ${\bold B}_W$. Thus we may view
$H_{\overline{W}}$ as a quotient of the group algebra of
$\pi_1(G^{rs})$ modulo quadratic relations.

Most fundamentally, Hecke algebras arise as convolution
algebras. Specifically, let $G$ be a split simply-connected reductive
group over ${\Bbb Z}$ with Borel subgroup $B$. Fix a finite field ${\Bbb
F}_q$. Write $G({\Bbb F}_q)$, $B({\Bbb F}_q)$ for the corresponding finite
groups of ${\Bbb F}_q$-rational points. Then ${\Bbb C} [B({\Bbb F}_q)
\backslash G({\Bbb F}_q) /  B({\Bbb F}_q)]$, the algebra (under
convolution) of ${\Bbb C}$-valued $B({\Bbb F}_q)$-biinvariant functions
on $G({\Bbb F}_q)$ is known to be isomorphic to the Hecke algebra $H_W$
specialized at ${\bold q} =q$, i.e.
\begin{equation}
\label{konec}
{\Bbb C} [B({\Bbb F}_q) \backslash G({\Bbb F}_q) /  B({\Bbb F}_q)] \simeq H_W
|_{{\bold q} = q}.
\end{equation}
Note that the convolution algebra on the left is the algebra of
intertwiners of the induced module $\Ind_{B({\Bbb F}_q)}^{G({\Bbb
    F}_q)} {\bold 1}$. Thus, decomposition of the induced module into
irreducible $G({\Bbb F}_q)$-modules is governed by representation
theory of the Hecke algebra $H_W|_{{\bold q} = q}.$

One has a similar interpretation of affine Hecke algebras in terms of
$p$-adic groups. Specifically, let ${\Bbb Q}_p$ be a $p$-adic field
with the ring of integers ${\Bbb Z}_p$ and the residue class field ${\Bbb
F}_p = {\Bbb Z}_p / p {\Bbb Z}_p $. Then the ring maps on the left (below)
induce the following group homomorphisms on the right (below):
$$
{\Bbb F}_p \stackrel{\pi}\twoheadleftarrow {\Bbb Z}_p \hookrightarrow {\Bbb
Q}_p, \quad \quad \quad \quad G({\Bbb F}_p)
\stackrel{\pi}\twoheadleftarrow 
G({\Bbb Z}_p) \hookrightarrow G({\Bbb Q}_p).
$$
The preimage of $B({\Bbb F}_p)$ under the projection $\pi$ is a compact
subgroup $I \subset G({\Bbb Z}_p)$, called an Iwahori subgroup. We may
consider the algebra ${\Bbb C} [I \backslash G({\Bbb Q}_p) / I ]$ of ${\Bbb
C}$-valued $I$-biinvariant functions on $G({\Bbb Q}_p)$ with compact
support. Similarly to (\ref{konec}) one establishes an algebra isomorphism
$$
{\Bbb C} [ I \backslash G({\Bbb F}_q) /  I ] \simeq H_{\overline{W}} 
|_{{\bold q} = p}.
$$

Now let  $\rho:G({\Bbb Q}_p)\to \mbox{{\rm End}}(V)$ be an
admissible representation of 
$G({\Bbb Q}_p)$. For any $I$-bi-invariant compactly supported function
$f$ on $G({\Bbb Q}_p)$, the formula:
$$
\rho(f): v \mapsto \int_{G({\Bbb Q}_p)} f(g) \cdot \rho(g)v\, dg,\quad v \in V^I
$$
defines a $\C[I\backslash G({\Bbb Q}_p)/I]$-module structure on the vector
space $V^I$ of $I$-fixed vectors.
Moreover, the space $V^I$ turns out to be finite-dimensional, and the
assignment $V \mapsto V^I$ is known to provide an equivalence between the category
of admissible $G({\Bbb Q}_p)$-modules 
generated by $I$-fixed vectors
and the category
of
 finite dimensional $H_{\overline{W}} |_{{\bold q} = p}$-modules.

\vspace{1cm}

The interpretations of Hecke algebras given above show the
importance of having a classification of their finite dimensional
irreducible representations. However, none of the above interpretations
helps in finding such a classification. For this, one needs a totally
different geometric interpretation that we are now going to explain.
\medskip

Let $\mu: \widetilde{\cal N}= T^*\cal B \to \cal N$ be the Springer
resolution and $Z = \widetilde{\cal N} \times_{\widetilde{\cal N}}
\widetilde{\cal N}$ the Steinberg variety.  Since $\mu: T^* \cal B
\to \cal N$ is $G$-equivariant, $Z$ is a $G$-variety and since $Z
\circ Z = Z$, the K-group $K^G(Z)$ has the structure of an associative
convolution algebra.
  
  Let $Z_{\Delta} \subset T^* \cal B \times T^* \cal B$ be the
  diagonal copy of $T^* \cal B$. The variety $Z_{\Delta}$ gets
  identified with $T^*_{\cal B_{\Delta}}$, the conormal bundle to the
  diagonal $\cal B_{\Delta} \subset \cal B \times \cal B$. This yields
  the following canonical isomorphisms of ${\bold R}(G)$-algebras
\begin{equation}
\label{k_diag}
K^G(Z_{\Delta})= K^G(T^*_{\cal B_{\Delta}}) \simeq K^G(\cal
B_{\Delta}) \simeq K^G(G/B) \simeq {\bold R}(T),
\end{equation}
where the second isomorphism is the Thom isomorphism (cf. [CG, Lemma 5.4.9]),
and the last one is the induction isomorphism (cf. [CG, Lemma 6.1.6]).

The following result is a $G$-equivariant extension of the geometric
construction of ${\Bbb Z}[W]$ given in Section 6.

\medskip

\noindent
{\bf Theorem 11.4 (see [CG, Theorem 7.2.2])} {\em There is a natural
  algebra isomorphism $K^G(Z) \simeq {\Bbb Z}[\widetilde{W}]$ making
  the following diagram commute}
$$
\begin{array}{cccc}
K^G(Z_{\Delta}) & \hookrightarrow & K^G(Z) & \\
\downarrow \wr & & \downarrow \wr & \\
{\bold R}(T) & \hookrightarrow & {\Bbb Z}[\widetilde{W}]. & \quad \Box
\end{array}
$$

\bigskip

\noindent
{\bf Affine Hecke algebras}
\smallskip

Notice that our picture has an extra symmetry: the group $\C^*$ acts
on $\cal N$, and also on $T^*( \cal B \times \cal B)$, along the
fibers,
by the formula ${\Bbb C}^*\ni z: x \mapsto z^{-1}\!\cdot\!x\,.$
Then, the Steinberg variety $Z = T^* \cal B \times_{\cal N} T^*\cal B$
is a $G \times {\Bbb C}^*$ stable subvariety of $T^*( \cal B \times \cal
B)$ and  we can   consider   $K^{G \times  {\Bbb  C}^*}$-theory  of  $Z$.

Note that any irreducible representation of ${\Bbb C}^*$ has the form $z
\mapsto z^m$ for some $m \in {\Bbb Z}$.
 Therefore we have the natural ring isomorphism
$$
{\bold R}(T) \simeq {\Bbb Z}[{\bold q}, {\bold q}^{-1}],
$$
where ${\bold q}$  is the tautological  representation ${\bold q}:{\Bbb C}^*
\to {\Bbb C}^*$ given by the identity map. One can prove the following 
``$G \times {\Bbb C}^*$-counterpart'' of (\ref{k_diag})
$$
K^{G \times {\Bbb C}^*}(Z_{\Delta}) \simeq {\bold R} (T \times {\Bbb C}^*)
\simeq {\bold R}(T)[{\bold q}, {\bold q}^{-1}].
$$ 
\medskip

\noindent
{\bf Theorem 11.5 (see [CG, Theorem 7.2.5])} {\it There is a natural
  algebra isomorphism $K^{G \times \C^*}(Z) \simeq {\bold H}$ making
  the following diagram commute}
$$
\begin{array}{cccc}
K^{G\times{\Bbb C}^*}(Z_{\Delta}) & \hookrightarrow & K^{G\times{\Bbb C}^*}(Z) & \\
\downarrow \wr & & \downarrow \wr & \\
{\bold R}(T)[{\bold q}, {\bold q}^{-1}] & \hookrightarrow & {\bold H}. & \quad \Box
\end{array}
$$

\bigskip
\noindent
{\bf Remarks} 

(1) For any $G \times {\Bbb C}^*$-variety $M$, the group $K^{G \times
  {\Bbb C}^*}(M)$ is a module over $K^{G \times {\Bbb C}^*}(pt)=
{\bold R}(G \times {\Bbb C}^*)$, the representation ring of $G \times
{\Bbb C}^*$. Restriction to $T\times {\Bbb C}^*$ gives a ring
isomorphism ${\bold R}(G\times {\Bbb C}^*) \simeq {\bold R}(T)^W
[{\bold q},{\bold q}^{-1}]$ where the RHS stands for $W$-invariants in
${\bold R}(T)[{\bold q}, {\bold q}^{-1}]$ (here ${\bold q}$ comes from
the representation ring of ${\Bbb C}^*$).  This ring ${\bold
  R}(T)^W[{\bold q}, {\bold q}^{-1}]$ gets identified, via the second
definition of ${\bold H}$, with a subalgebra of ${\bold H}$.  One can
prove [CG, Proposition 7.1.14] that it coincides with $Z({\bold H})$,
the center of ${\bold H}$.

(2) Recall that the restriction to the Steinberg variety $Z \subset
T^* \cal B \times T^* \cal B$ of either of the two projections $T^*
\cal B \times T^* \cal B \to T^* \cal B$ is proper.  Thus, the
convolution product yields a $K^{G \times {\Bbb C}^*}(Z)$-module
structure on $K^{G \times {\Bbb C}^*}(T^* \cal B)$.  Recall that $K^{G
  \times {\Bbb C}^*}(T^* \cal B) \simeq {\bold R}(T)[{\bold q}, {\bold
  q}^{-1}]$.  Hence, for any $s_{\alpha} \in W$ convolution action of
the class in $K^{G \times {\Bbb C}^*}(Z)$, corresponding to the
element $T_{s_{\alpha}}$ via Theorem 11.5, gives an operator
$\widehat{T}_{s_{\alpha}} \in \mbox{{\rm End}}_{{\Bbb Z}[{\bold q},
  {\bold q}^{-1}]} {\bold R}(T)[{\bold q}, {\bold q}^{-1}] $ (cf.
Definition 11.2).  One can find the following explicit formula for the
action of $\widehat{T}_{s_{\alpha}}$
\begin{equation}
\label{dem_lus}
\widehat{T}_{s_{\alpha}}: Y_{\lambda} \mapsto {Y_{\lambda} -
  Y_{s_{\alpha}(\lambda)} \over Y_{\alpha} - 1} - {\bold q}
{Y_{\lambda} - Y_{s_{\alpha}(\lambda)+ \alpha} \over Y_{\alpha} - 1}.
\end{equation}
This   formula,   discovered by  Lusztig,    was  a  starting  point of
the K-theoretic approach to Hecke algebras.

(3)   Theorem  11.5 implies   Deligne-Langlands-Lusztig conjecture for
${\bold H}$, see [CG], [KL1].

(4) A Fourier  transform argument  does not  work for Theorems  11.2 and
11.5.

\medskip One can give a description of finite-dimensional irreducible
complex representations of ${\bold H}$ similar to that of the Springer
Theorem for $W$.  Firstly, by Remark (1) above, we have a canonical
algebra isomorphism $Z({\bold H}) \simeq {\bold R}(G \times {\Bbb
  C}^*)$.  On any irreducible representation the center $Z({\bold H})$
of ${\bold H}$ acts via an algebra homomorphism $Z({\bold H}) \to
{\Bbb C}$, due to Schur's lemma.  Any such homomorphism may be
identified [CG, 8.1] with the evaluation homomorphism sending a
character $z \in {\bold R}(G \times {\Bbb C}^*)$ to $z(a)$, the value
of $z$ at a semisimple element $a=(s, q) \in G \times {\Bbb C}^*$.  In
particular, the indeterminate ${\bold q}$ specializes to a complex
number $q \in {\Bbb C}^*$.

Given a semisimple element $a=(s,q)\in G\times {\Bbb C}^*$, let ${\Bbb
  C}_a$ be the 1-dimensional complex vector space $\C$ viewed as a
$Z({\bold H})$, equivalently ${\bold R}(G\times \C^*)$-module, via the
action
$$
{\bold R}(G\times {\Bbb C}^*)\times {\Bbb C}_a\to \C_a,\qquad
(z,x)\mapsto z(a)\cdot x,
$$
where $z \mapsto z(a)$ is the corresponding evaluation homomorphism at $a$.

\medskip

\noindent
{\bf Definition  11.6 } The tensor    product ${\bold H}_a:={\Bbb  C}_a
\otimes_{_{Z({\bold H})}}{\bold H}$ is called the Hecke algebra {\it
  specialized} at $a$.

\medskip

Thus, for any simple ${\bold H}$-module  $M$ there exists a semisimple
element $a=(s,t)\in G\times\C^*$, such that the action of ${\bold H}$
factors through an action of  the specialized Hecke algebra  ${\bold
H}_a$. From  now  on we will  fix $a$  and  consider representations of
${\bold H}_a$.

Recall [CG, 6.2] that $a$ acts on  $\cal N$ by the  formula $(s, q): x
\to q^{-1} \cdot sxs^{-1}$  (note  the {\it inverse}  power of  $q$)
 and this action agrees with the action of $a$ on $\tilde{\cal N}= T^*
\cal B $, given by a similar formula.  Denote by  $\tilde{\cal N}^a\,,\,\cal 
 N^a$ and $Z^a$ the   corresponding $a$-fixed point subvarieties.  The
variety $\tilde{\cal N}^a$ is smooth  due to [CG, Lemma 5.11.1],  since
$\tilde{\cal   N}$ is smooth.      Observe   further  that  we    have
$Z^a=\tilde{\cal N}^a\times_{{\cal N}^a} \tilde{\cal N}^a$.  Therefore
$Z^a$ may be viewed as a subvariety in $\tilde{\cal N}^a \times
\tilde{\cal N}^a$ such that
$$
Z^a\circ Z^a= Z^a.
$$
Our general construction of Section 2 (Corollary 2.1)
 makes the Borel-Moore homology
$H_\bullet(Z^a)$ an associative algebra via convolution.

\medskip

\noindent
{\bf Proposition 11.7 [CG, Proposition 8.1.5]} {\it Let $a=(s,q)\in
  G\times{\Bbb C}^*$ be a semisimple element.  Then there is a natural
  algebra isomorphism
$$
{\bold H}_a\simeq H_*(Z^a, {\Bbb C}). \quad \quad \Box
$$} 

\noindent
{\bf Remark on the proof of Proposition 11.7} Let $\cal A$ be the
closed subgroup of $G\times\C^*$ generated by $a$, that is, the
algebraic closure in $G\times\C^*$ of the cyclic group $\{a^n\,,\,
n\in \Z\}\,.$ Clearly $\cal A$ is an abelian reductive subgroup of
$G\times \C^*$, and we have $Z^a=Z^{\cal A}$.  The isomorphism of
Proposition 11.7 is constructed as a composite of the following chain
of algebra isomorphisms (cf. [CG, 8.1.6])
\begin{eqnarray*}
\C_a\otimes_{Z({\bold H})} {\bold H} 
&\stackrel{\sim}\to&
   \C_a\otimes_{_{{\bold R}(G\times\C^*)}}K^{G\times\C^*}(Z)
   \; \simeq \;  \C_a\otimes_{{\bold R}(\cal A)} K^{\cal A}(Z) \\
&\stackrel{r_a}{ \simeq} &
   \C_a\otimes K^{\cal A}(Z^{\cal A})\;\stackrel{\mbox{\small
ev}}\simeq \;  K_\C(Z^{\cal A})\;\stackrel{RR} \simeq \; 
H_\bullet(Z^{\cal  A}, \C)= H_\bullet(Z^a, \C).
\end{eqnarray*}

The first isomorphism  here  is given by  Theorem 11.5  and the first
remark after it, the second  is given by  the restriction property for
$\cal A \subset G \times {\Bbb  C}^*$ (cf.  [CG,  section 4.2(6)]).  The
third map is  given by the algebra  homomorphism $r_a$ of localization
(cf. [CG,  Theorem 5.11.10]).  The fourth  map $\mbox{ev}:\C_a\otimes_{{\bold
R}(\cal A)} K^{\cal A}(Z^{\cal  A})\simeq \C_a\otimes_{{\bold R}(\cal A)}
\big({\bold R}({\cal   A})\otimes  K(Z^{\cal A})\big)\stackrel{\sim}\to
 K_\C(Z^{\cal A})$ is the  evaluation map sending $1\otimes f  \otimes
[\cal F]$ to $f(a)\otimes [\cal F]$ where $f \in {\bold R}({\cal A})$ is
viewed as a character function on ${\cal A}$ and $[\cal F] \in K^{\cal
A}(Z^{\cal  A})$.  The last isomorphism is  the map $RR$  given by the
bivariant Riemann-Roch theorem (cf. [CG, Theorem 5.11.11]). $\Box$

\medskip

We will construct,  for each semisimple  $a=(s,q)\in  G\times \C^*$, a
complete collection of simple ${\bold H}_a$-modules, which will yield
a complete collection of simple ${\bold  H}$-modules as $a$ runs over
all semisimple conjugacy classes in $G \times {\Bbb C}^*$.

To   that  end, consider the map    $\,\mu: \tilde{\cal  N}^a \to \cal
N^a\,,$ the restriction of the Springer resolution  to the fixed point
varieties. Explicitly, we have
$$    
\cal     N^a=\{x\in\cal    N     \mid   sxs^{-1}     =    q\cdot
x\}\quad,\quad\tilde{\cal   N}^a= \{(x,{\frak b})\in\cal   N^a\times\cal
B^a\mid x\in \frak b\}.
$$

Let  $x\in \cal   N^a$.  The  fiber $\mu^{-1}(x)\,\subset\,\tilde{\cal
  N}^a$ may be identified  via  the projection $\cal N^a\to\cal  B\,,$
$\,(x,\frak b)\mapsto \frak   b\,$,  with the  subvariety  $\cal  B^s_x
\,\subset\,  \cal B$ of all Borel  subalgebras simultaneously fixed by
$s$ and $x$.
\medskip

\noindent
{\bf Remark} {\it The variety $\cal B^s_x$ is non-empty.}
\smallskip

\noindent
{\bf Proof }  Recall that $a = (s,q)$   and the relation  $s x s^{-1}=
q\cdot  x$ holds.   Let  $u=\mbox{exp}(z\cdot x)\in G\,,\,z\in{\Bbb C}$.
Then  $u$    is    a  unipotent   element     of   $G$   and   clearly
$sus^{-1}=\mbox{exp}(z\cdot q\cdot x)$.   We see that the elements $s$
and  $\mbox{exp}(z\cdot  x)\,,\,z\in{\Bbb  C}\,,$  generate a   solvable
subgroup of $G$.  Hence there  exists a Borel subgroup $B$  containing
this solvable subgroup.  It  follows that its Lie  algebra is in $\cal
B^s_x$.  $\Box$

\medskip

By our general construction in the end of Section 2, the Borel-Moore
homology of the fibers of the map $\,\mu: \tilde{\cal N}^a\to \cal
N^a\,$ have a natural $H_*(Z^a)$-module structure via convolution.
Hence, for any $x \in \cal N^a$, we get an $H_*(Z^a)$-action on
$H_*(\cal B_x^s)$.  Further, let $G(s,x)$ be the simultaneous
centralizer in $G$ of $s$ and $x$, and let $G(s,x)^{\circ}$ be the
connected component of the identity.  As at the end of Section 6, we
define $A(s, x)= G(s,x)/ G(s,x)^{\circ}$ and notice that the action of
$G(s, x)$ on $\cal B_x^s$ induces an action of $A(s, x)$ on $H_*(\cal
B_x^s)$. We write $A(s, x)^\vee$ for the set of all irreducible
representations of $A(s, x)$ that occur with non-zero multiplicity in
the homology of $\cal B_x^s$.

\medskip

\noindent
{\bf Remark} Notice that in Section 6 we used only the {\it top}
homology to define $A^\vee$ and to state the analogue of Proposition
11.7. However, here we use all homology groups since Proposition 11.7
arises from the K-theoretic Theorem 11.5, and the K-groups are not
graded by dimension.

\medskip

The variety $\cal N^a$ is stable under the adjoint action of the group
$G(s)$, the centralizer of $s$ in $G$. Moreover, $\cal N^a$ is known
[CG, Proposition 8.1.17] to be a finite union of $G(s)$-orbits.

 Now we can apply the direct image decomposition (\ref{decom-star-ng})
 to the morphism $\mu : \tilde{\cal N}^a \to \cal N^a$, to get
\begin{equation}
\label{prikol}
\mu_* \cal C_{\tilde{\cal N}^a} = \bigoplus_{\phi} L_{\phi}(k) \otimes
IC_{\phi}[k],
\end{equation}
where $\phi$ runs over the set of pairs: $(G(s)$-{\it orbit in} $\cal
N^a, G(s)$-{\it equivariant irreducible local system on this orbit}).
Choosing a base point, we can write each $G(s)$-orbit in the form
$G(s)\cdot x\,,\, x\in \cal N^a$; then giving an equivariant
irreducible local system on $G(s) \cdot x$ amounts to giving an
irreducible representation $\chi$ of the group $A(s, x)$. Thus,
Theorem 5.2 says that the multiplicity spaces $L_{\phi}
= \bigoplus_{k \in {\Bbb Z}}\, L_{\phi}(k)$, where $\phi= (G(s) \cdot
x, \chi)$, form a complete collection of simple $H_*(Z^a, {\Bbb
  C})$-modules, to be denoted $L_{a, x,\chi} :=L_{\phi}$. We will now
specify more precisely the range of the parameters $(a, x,\chi)$
labelling the isomorphism classes of {\it non-zero} modules $L_{a,
  x,\chi}$ occurring in this parametrization.
\medskip

We say that two pairs $(x, \chi)$ and $(x', \chi')$ are $G(s)$-{\it
  conjugate} if there is a $g \in G(s)$ such that $x'=gxg^{-1}$ and
conjugation by $g$ intertwines the $A(s, x)$-module $\chi$ with the $A(s,
x')$-module $\chi'$.
Write ${\bold M}$ for the set of $G$-conjugacy classes of the triple
data
$$
{\bold M} = \{a = (s, q) \in G \times {\Bbb C}^*, x \in \cal N^a, \chi \in
A(s, x)^\vee \; | \; s \ is \  semisimple \}/ \mbox{Ad}G
$$ 
The main result on representations of ${\bold H}$ is now deduced from
Theorem 5.2 and reads as follows:

\medskip

\noindent
{\bf Theorem 11.8 (cf.\ [CG, 8.1.13-16] and [KL1])} {\it Assume that $q
  \in {\Bbb C}^*$ is not a root of unity. Then

{\rm (i)} Two modules $L_{a, x, \chi}$ and $L_{a, x', \chi'}$ are isomorphic
if and only if the pairs $(x, \chi)$ and $(x', \chi')$ are
$G(s)$-conjugate to each other.

{\rm (ii)} $L_{a, x, \chi}$ is non-zero simple $H_*(Z^a)$-module, for any
$(a, x, \chi) \in {\bold M} $.

{\rm (iii)} The collection $\{L_{a, x, \chi}\}_{(a, x, \chi) \in {\bold
    M}}$ is a complete set of irreducible ${\bold H}$-modules such
that ${\bold q}$ acts by $q$.} $\Box$ \medskip

\noindent
{\bf Remark}

(i) Theorem  11.8 fails if $q$  is a root of  unity.  In fact, in this
case  some of   the  modules  $L_{a,  x,  \chi}$   may  be  zero.  The
classification of Theorem  11.8(iii) was first  obtained  in [KL1]  in a
different way, see [CG, Introduction] for more historical remarks.

(ii) Note that the morphism $\mu: \tilde{\cal N}^a \to \cal N^a$ is
{\it not} semi-small and, as a result (cf. Proposition 5.1), the
algebra ${\bold H}_a$ is not semisimple.

\bigskip

\noindent
{\bf Quantized loop algebra of ${\frak g} {\frak l}_n$}
\smallskip

Recall briefly the setup of Section 7. We have a variety $N$ of
endomorphisms $x: \C^d \to {\Bbb C}^d$ that satisfy $x^n=0$, a smooth
variety $\cal F$ of $n$-step partial flags in $\C^d$, and a
(semi-small) map $\mu: M= T^* \cal F \to N$. Denote as usual $Z = M
\times_N M$.
  
The group $GL_d \times {\Bbb C}^*$ acts on $M$ and $N$. Since the map
$\mu: M \to N$ is $GL_d \times {\Bbb C}^*$-equivariant, the variety
$Z$ is a $GL_d\times {\Bbb C}^*$-stable subvariety of $M \times M$ and
we can consider $K^{GL_d \times {\Bbb C}^*}(Z)$.
 
To describe an algebraic object arising from equivariant K-theory of
$Z$, consider the {\it loop Lie algebra} $\mbox{L}{\frak g} = {\frak
  g} {\frak l}_n \otimes_{{\Bbb C}} {\Bbb C}[u, u^{-1}]$ which may be
thought of as the space of polynomial maps ${\Bbb C}^* \to {\frak g}
{\frak l}_n$.  One can define (cf. [Dr1], [CP]) a deformation $\cal
U_{{\bold q}}(\mbox{L}{\frak g} )$ of the universal enveloping algebra
$\cal U (\mbox{L}{\frak g} )$. We will call this deformation {\it the
  quantized loop algebra} of ${\frak g} {\frak l}_n$.

We assume, throughout, that $q$ is {\it not a root of unity}.  The
following theorem is announced in [GV] and proved in [V]:

\medskip

\noindent
{\bf Theorem 11.9} {\em There exists a surjective algebra
 homomorphism}
$$
\cal U_{{\bold q}}(\mbox{L}{\frak g} ) \twoheadrightarrow 
{\Bbb C}  \otimes_{{\Bbb Z}} K^{GL_d \times \C^*}(Z). \quad \Box
$$

\medskip

Recall that the algebra $\cal U_{{\bold q}}(\mbox{L}{\frak g} )$ has
generators $E_{i, r}, F_{i, r}, K_i$, $i = 1, \ldots, n$, $r \in {\Bbb
  Z}$ and $H_{i, r}$, $i= 1, \ldots, n$, $r \in {\Bbb Z} \setminus
\{0\}$ subject to a certain explicit set of relations (cf. [CP]).
Recall further that there is a tensor product decomposition as a
vector space 
$$\cal U_{{\bold q}}(\mbox{L}{\frak g} )= \cal U^+ \otimes_{{\Bbb C} [
  {\bold q}, {\bold q}^{-1}]} \cal U^0 \otimes_{{\Bbb C} [ {\bold q},
  {\bold q}^{-1}]} \cal U^-,$$
where $\cal U^+$, $\cal U^0$, $\cal
U^-$ are the subalgebras generated by $\{E_{i, r}\}$, $\{K_i, H_{i,
  r}\}$ and $\{F_{i,r}\}$, respectively.

 $\cal  U^0$ is  a   large commutative  subalgebra   of $\cal U_{{\bold
 q}}(\mbox{L}{\frak g}   )$.   A     finite-dimensional     irreducible
 representation  $V$ of $\cal U_{{\bold q}}(\mbox{L}  {\frak g} )$ is said
 to be {\it of  type} {\bf 1} if  $K_1, \ldots, K_n$ act semisimply on
 $V$ with eigenvalues which are half-integer powers of $q$, the
 specialization of ${\bold q}$-action in $V$.

 We say that $v \in V$ is a {\it pseudo-highest weight vector} if it
 is annihilated by $\cal U^+$ and is a weight vector for $\cal U^0$.
 We write $k_i(V), h_{i, r}(V) $, $i= 1, \ldots, n$, $r \in {\Bbb Z}
 \setminus \{0\}$, for the corresponding eigenvalues of the elements
 $K_i, H_{i, r}$, and call the collection $\{k_i(V), h_{i, r}(V)\} $
 {\it quasi-highest weight}. The following result about irreducible
 finite-dimensional representations of $\cal U_{{\bold
     q}}(\mbox{L}{\frak g})$ is essentially due to Drinfeld [Dr2], see
 [CP, Theorems 12.2.3 and 12.2.6].

\medskip

\noindent
{\bf Proposition 11.10 }
 {\it 
   
{\rm (i)} Any finite-dimensional simple $\cal U_{{\bold
    q}}(\mbox{L} {\frak g})$-module $V$ of type ${\bold 1}$ has a
unique quasi-highest weight vector.

{\rm (ii)} Two simple modules of type ${\bold 1}$ are isomorphic iff their
quasi-highest weights  are the
same.

{\rm (iii)} A collection $\{k_i, h_{i, r} \}$ is the quasi-highest weight of
a finite-dimensional irreducible representation  $V$ iff there exist
  unique monic polynomials $P_{i, V} \in {\Bbb C}[u]$, $i = 1, \ldots,
  n$, all with non-zero constant term, such that, setting $d_i = deg
(P_{i, V})$, we have :
$$
k_i \; exp\Big((q^{1/2}- q^{-1/2}) \sum_{s=1}^{\infty}h_{i, s} u^{s} \Big)
= q^{{ d_i \over 2}} {P_{i, V}(q^{-1}u) \over P_{i, V}(u)} =
k_i^{-1} exp\Big((q^{-1/2}- q^{1/2}) \sum_{s=1}^{\infty}h_{i, -s} u^{-s} \Big)
$$
in the sense that the left- and right-hand sides are the Laurent
expansions of the middle term about $0$ and $\infty$,
respectively. 

Moreover, the polynomials $P_{i, V}$ define the
representation $V$ uniquely, and every n-tuple $(P_i)_{i= 1, \ldots,
  n}$ of monic polynomials with non-zero constant term arises from a
finite-dimensional irreducible $\cal U_{{\bold q}}(\mbox{L}{\frak g})$-module
of type ${\bold 1}$ in this way.} $\Box$

\medskip

The collection $\{P_{i, V}, i = 1, \ldots, n \}$ is called the {\it
  Drinfeld polynomials} associated to the irreducible representation $V$.

For any pair $a = (s, q)$, $ s \in GL_d({\Bbb C})$, $q \in {\Bbb
  C}^*$, we consider the fixed point variety $Z^a$ and construct as in
Proposition 11.6 a surjection ${\Bbb C} \otimes K^{GL_d({\Bbb C})
  \times {\Bbb C}^*}(Z) \twoheadrightarrow H_*(Z^a)$. As in the case
of affine Hecke algebras we obtain from Theorem 5.2 a complete set of
irreducible finite dimensional $H_*(Z^a)$-modules $L_{a, x}$ labeled
by pairs $(a, x)$, where $a=(s, q)$ is a semisimple element of $GL_d
\times {\Bbb C}^*$ and $x \in N$ (as in Lemma 7.3, the monodromy
action is trivial, hence $L_{a, x}$ does not depend on the third
parameter $\chi$).
 
Observe that every simple $K^{GL_d({\Bbb C}) \times {\Bbb
    C}^*}(Z)$-module can be pulled back via the surjection $\cal
U_{{\bold q}}(\mbox{L}{\frak g} ) \twoheadrightarrow {\Bbb C}
\otimes_{{\Bbb Z}} K^{GL_d \times \C^*}(Z)$ of Theorem 11.8 to give a
simple $\cal U_{{\bold q}}(L {\frak g})$-module.  Now we are ready to
identify the modules $L_{a, x}$ as representations of $\cal U_{{\bold
    q}}(\mbox{L} {\frak g})$.
 
Fix $(a, x)$ such that $a=(s, q) \in GL_d \times {\Bbb C}^*$, where
$q$ is not a root of unity, $x \in N$, $s$ is semisimple and $s x
s^{-1}= q x$.  Recall that we have defined before Lemma 7.5 two
$n$-step flags $F^{min}(x)$ and $F^{max}(x)$ in ${\Bbb C}^d$.  Since
$sxs^{-1} = qx$, the flags $F^{min}(x)$ and $ F^{max}(x)$ are
both preserved by $s$ and we can consider, for each $i = 1, \ldots,
n$, the action of $s$ on $F^{max}_i(x)/F^{min}_i(x)$.
 
Let $L_{a, x}$ be the $K^{GL_d \times {\Bbb C}^*}(Z)$-module viewed,
due to the Theorem 11.8, as a $\cal U_{{\bold q}}(\mbox{L}{\frak g})
$-module.  The theorem below is a quantized version of the Springer
theorem for $\cal U({\frak s} {\frak l}_n )$ given in section 7.

\medskip

\noindent
{\bf Theorem 11.11} {\it The $i$-th Drinfeld polynomial, $P_{i,L_{a,
      x}}$, is equal to the characteristic polynomial of the
  $s$-action on $F^{max}_i(x)/F^{min}_i(x), i.e. :$}
$$
P_{i,L_{a, x}}(u)= \mbox{det}(u\cdot Id - s; F^{max}_i(x)/F^{min}_i(x)
)\; \mbox{\it for all } i = 1, 2, \ldots, n. \quad \quad  
$$
{\it In particular, every irreducible finite-dimensional $\cal
  U_{{\bold q}}(\mbox{L}{\frak g}) $-module of type $1$ is of the
  form} $L_{a, x}.\quad\Box$

\section{Equivariant cohomology and degenerate versions}

In this last section we will study degenerate versions of Hecke
algebras and quantized enveloping algebras as the deformation
parameter $q \to 1$. Of course, in the limit $q = 1$, the algebras in
question reduce to their classical counterparts: the Hecke algebra
specializes to the group algebra of the corresponding Weyl group and
the quantized enveloping algebra specializes to the corresponding
classical enveloping algebra. We will consider however another, more
interesting limit, which corresponds in  a sense to taking the ``first
derivative'' with respect to the deformation parameter at $q = 1$,
rather than the value at $q = 1$ itself. We will see that taking
``first derivative at $q=1$'' amounts geometrically to replacing
equivariant K-theory by equivariant cohomology.

\bigskip

\noindent
{\bf The degenerate Hecke algebra}
\smallskip

 Let $\epsilon$ be an indeterminate. Write ${\Bbb
C}[{\frak h} , \epsilon]$ for polynomials in $\epsilon$ with
coefficients in the ring ${\Bbb C} [ {\frak h}]$ of polynomial functions
on $\frak h$.
\smallskip

\noindent
{\bf Definition 12.1} 
The {\it  degenerate affine Hecke algebra
  ${\bold  H}_{deg}$} of $G$ is the unital associative free ${\Bbb
C}[\epsilon]$-algebra  defined by the following properties:
\smallskip

(i)   ${\bold    H}_{deg}   \simeq  {\Bbb   C}   [W]     \otimes_{{\Bbb  C}}
{\Bbb C}[{\frak h}, \epsilon] $ as a ${\Bbb C}$-vector space.

(ii)   the  subspaces    ${\Bbb  C}[W]$ and   $  {\Bbb C} [{\frak  h}, 
\epsilon]$ are subalgebras of ${\bold H}_{deg} $.

(iii) the following relations hold in ${\bold  H}_{deg}$:
$$
 s_i \lambda - s_i(\lambda)s_i = -  \epsilon \cdot \langle
 \alpha^\vee_i , \lambda \rangle,
 \quad \quad i = 1, \ldots, l, \quad \lambda \in  {\frak h}^* \subset
 {\Bbb C} [{\frak h}].
$$
 The algebra ${\bold H}_{deg}$ has a natural grading defined by
$deg(s_i)=0$, $deg(\alpha)=1$, $deg(\epsilon)=1$.

\medskip
\noindent
{\bf Remarks} 

(1) Sometimes in the definition of the degenerate affine Hecke algebra
one takes a  quotient  modulo relation  $\epsilon =  1$. In  fact, all
algebras with $\epsilon$ specialized to a  non-zero complex number are
isomorphic.  We prefer the homogeneous version above  since we want to
relate it to the equivariant cohomology with its natural grading.

(2)  The degenerate affine  Hecke    algebra ${\bold H}_{deg}$ can    be
obtained from   the affine Hecke  algebra  ${\bold H}$  by the following
procedure: we formally make substitution ${\bold q} \mapsto \mbox{exp}(2
\epsilon)$,  $Y_{\lambda} \mapsto   \mbox{exp}(\epsilon     \lambda)$,
$T_{s_\alpha} \mapsto s_{\alpha}$ in the relations defining the affine
Hecke algebra (cf.  Definition 11.2) ,   and then take  the homogeneous
components of  the minimal degree  with respect to  the grading above. 
This   last step is  sometimes  expressed  as ``taking  $\displaystyle
\lim_{\epsilon \to 0}$ ''.

\medskip

\noindent
{\bf Demazure - Lusztig type operators}
\smallskip

Recall that algebraic group homomorphisms $\lambda \in X^*(T)$ clearly
form  a ${\Bbb Z}$-basis of  the  representation ring ${\bold R}(T)$, that
is,  ${\bold R}(T)$  may  be identified with  the  group algebra  of the
lattice $X^*(T)$.  The Weyl group $W$  acts naturally  on ${\bold R}(T)$
and on ${\Bbb C} [ {\frak h}]$;  we write $P \mapsto  w(P)$ for the action
of $w \in W$.

 To each simple reflection $s_{\alpha} \in W$ we have associated in
 (\ref{dem_lus}) a ${\Bbb Z} [ {\bold q}, {\bold q}^{-1}]$-linear map
 $\widehat{T}_{\alpha} : {\bold R}(T) [{\bold q}, {\bold q}^{-1}] \to 
{\bold R}(T) [{\bold q}, {\bold q}^{-1}]$ given by the Demazure-Lusztig formula.

 Similarly, we define a ${\Bbb C} [\epsilon]$-linear map $S_{\alpha}:
 {\Bbb C}  [{\frak h}, \epsilon] \to {\Bbb C}  [{\frak h}, \epsilon]$ by the
 formula, due to [BGG]:
\begin{equation}
\label{nilki}
S_{\alpha}: P \mapsto s_{\alpha}(P) + \epsilon \; {s_{\alpha}(P) -P \over
\alpha}, \quad \quad P \in {\Bbb C} [ {\frak h}].
\end{equation}


\noindent
{\bf Proposition 12.2 (see [Dr2], [Lu3], [CG, Theorem 7.2.16])}

{\it {\rm (i)} The ${\Bbb Z} [{\bold q}, {\bold q}^{-1}]$-subalgebra
  of ${\rm End}_{{\Bbb Z} [{\bold q}, {\bold q}^{-1}]} {\bold
    R}(T)[{\bold q}, {\bold q}^{-1}]$ generated by the operators
  $\{\widehat{T}_{\alpha} , \newline \alpha$ simple root $\}$ and by
  all the multiplication operators $P \mapsto f \cdot P$, $f \in
  {\bold R}(T)$, is isomorphic to the affine Hecke algebra ${\bold
    H}$.
  
  {\rm (ii)} The ${\Bbb C} [\epsilon]$-subalgebra of ${\rm End}_{{\Bbb
      C} [\epsilon]} {\Bbb C} [{\frak h}, \epsilon]$ generated by the
  operators $\{S_{\alpha}, \alpha$ simple root $\}$ and by all the
  multiplication operators $P \mapsto f \cdot P$, $f \in {\Bbb C} [
  {\frak h}]$, is isomorphic to the degenerate affine Hecke algebra
  ${\bold H}_{deg}$.} $\Box$

\medskip

Using this proposition we will frequently identify ${\bold H}$, ${\bold
H}_{deg}$ with the corresponding algebras of operators. 

Let ${\Bbb C} [T] = {\Bbb C} \otimes_{{\Bbb Z}} {\bold R}(T)$ be the
algebra of regular functions on $T$.  Pull-back via the exponential
map $exp: {\frak h} \to T$ gives an embedding $\exp^*: {\Bbb C} [T]
\hookrightarrow {\Bbb C} [[{\frak h}]]$. Similarly, one gets an
embedding $exp^*: {\Bbb C} [T \times {\Bbb C}^*] = {\Bbb C} [T][{\bold
  q}, {\bold q}^{-1}] \hookrightarrow {\Bbb C} [[{\frak h},
\epsilon]]$ (here $\epsilon$ is viewed as a base vector in $Lie \;
{\Bbb C}^*$ and the map $exp^*$ takes ${\bold q}$ to $exp(\epsilon) =
\sum {\epsilon^k \over k!}$). Using the latter embedding and
Proposition 12.2 we may (and will) view ${\bold H}$ as a subalgebra of
${\rm End}_{{\Bbb C}[[\epsilon]]} {\Bbb C}[[{\frak h}, \epsilon]]$ and
form the ${\Bbb C}[[\epsilon]]$-algebra $\widehat{{\bold H}}:= {\Bbb
  C}[[\epsilon]] \otimes_{{\Bbb Z}[{\bold q}, {\bold q}^{-1}]}
{\bold H}$.

 Let $I$ denote the augmentation ideal in ${\Bbb C}[[{\frak h},
 \epsilon]]$. The powers $I^n, n = 1, 2, \ldots$, form the decreasing
 $I$-adic filtration on ${\Bbb C}[[{\frak h}, \epsilon]]$. We define a
 {\it decreasing} ${\Bbb Z}$-filtration $F^{\bullet}$ on ${\rm End}_{{\Bbb
 C}[[\epsilon]]} {\Bbb C}[[{\frak h}, \epsilon]]$ as follows
$$
F^j = \{ u \; | \; u(I^{n}) \subset I^{n+j}, \; \mbox{ for all } n
> - j \}, \quad j \in {\Bbb Z}.
$$

We will write $F^{\bullet} \widehat{{\bold H}}$ for the induced
filtration on the subalgebra $\widehat{{\bold H}} \subset {\rm End}_{{\Bbb
    C}[[\epsilon]]} {\Bbb C}[[{\frak h}, \epsilon]]$. Comparing the
defining relations in ${\bold H}$ and ${\bold H}_{deg}$, one deduces
from Proposition 12.2 the following
\medskip

\noindent
{\bf Proposition 12.3} {\it There is a graded algebra isomorphism}
$$
{\bold H}_{deg} \simeq \mbox{gr}_F \widehat{{\bold H}},
$$
{\it where the RHS denotes the associated graded algebra of $\widehat{{\bold
H}}$.} $\Box$

\medskip
Note that although filtration terms $F^j{\rm End}_{{\Bbb C} [[ \epsilon]]}
{\Bbb C} [[{\frak h}, \epsilon]]$ are non-zero for both positive and
negative values of $j$, formula (\ref{nilki}) shows that the
induced filtration $F^j\widehat{{\bold H}}$ has non-zero terms only for
$j \geq 0$.

\bigskip
\pagebreak
\noindent
{\bf Nil-Hecke algebra}\\[3pt]%
\nopagebreak%
\noindent
Let {\bf Nil} be the ${\Bbb C}$-algebra with base $\{r_w, w \in W\}$
and the following multiplication rules \smallskip

(i) $r_w \cdot r_{w'} = r_{ww'}$ whenever $l(w) + l(w') = l(ww')$,

(ii) $r_w \cdot r_{w'} = 0 $ if $l(w) + l(w') > l(ww')$.
\smallskip

\noindent
Note that $r_e$ is the unit of the algebra {\bf Nil}.

\medskip

\noindent
{\bf Definition 12.4} The {\it nil-Hecke algebra} ${\bold H}_{nil}$ is the
unital associative algebra defined by the following properties:
\smallskip

(1) ${\bold H}_{nil} =\,${\bf Nil}$\,\otimes_{{\Bbb C}} {\Bbb C}
[{\frak h}]$ as a ${\Bbb C}$-vector space;

(2) The natural embeddings {\bf Nil} $\hookrightarrow {\bold H}_{nil}$
and ${\Bbb C} [{\frak h}] \hookrightarrow {\bold H}_{nil}$ are algebra
homomorphisms;

(3) For any simple root $\alpha \in {\frak h}^*$ and any linear
function $\lambda \in {\frak h}^* \subset {\Bbb C} [{\frak h}]$ the
following relations hold in ${\bold H}_{nil}$:
$$
r_{s_{\alpha}} \cdot \lambda - s_{\alpha}(\lambda) \cdot  r_{s_{\alpha}} 
= - \langle\lambda, \alpha ^\vee\rangle.
$$

To formulate an analogue of Proposition 12.2, consider for each simple
root $\alpha$ the linear operator $R_{\alpha}: {\Bbb C} [{\frak h}]
\to {\Bbb C} [{\frak h}]$ given by:
\begin{equation}
\label{dima}
R_{\alpha}: P \mapsto {s_{\alpha}(P) - P \over \alpha}, \quad \quad P \in
 {\Bbb C} [{\frak h}].
\end{equation}

\medskip

\noindent
{\bf Proposition 12.5 (cf. [KK])} {\it The subalgebra of ${\Bbb
    C}$-linear endomorphisms of ${\Bbb C}[{\frak h}]$ generated by
  $\{R_{\alpha}, \alpha$ simple root $\}$ and operators of
  multiplication by the elements of ${\Bbb C} [{\frak h}]$, is
  isomorphic to the nil-Hecke algebra ${\bold H}_{nil}$.} $\Box$

\medskip

\noindent
{\bf Remarks}

(1) Comparison of formulas (\ref{dem_lus}) and (\ref{dima}) shows that
the algebra ${\bold H}_{nil}$ may be thought of as the limit of the
affine Hecke algebra ${\bold H}$ at ${\bold q} \to 0$. Note that since the
variable ${\bold q}$ is invertible in ${\bold H}$ one cannot simply set
${\bold q}$ equal to 0.

(2) The algebra ${\bold H}_{nil}$ has a natural ${\Bbb Z}$-grading
defined on generators by
$$
deg\; r_w = - l(w), \quad deg \; \lambda = 1 \; \mbox{ for all }
\lambda \in {\frak h}^* \subset {\Bbb C} [{\frak h}].
$$

\medskip

There is also an interpretation of the nil-Hecke algebra as an
associated graded algebra, analogous to Proposition 12.3. To explain
it, fix a complex number $q \in {\Bbb C}^*$, $q \neq 1$. Specializing
in the Demazure-Lusztig formula (\ref{dem_lus}) and in Proposition
12.2, the variable ${\bold q}$ to $q$ we obtain an action of ${\bold
  H}_q$, the specialized affine Hecke algebra, on the ${\Bbb
  C}$-vector space ${\Bbb C} [T]$. Using the embedding ${\Bbb C} [T]
\hookrightarrow {\Bbb C} [[{\frak h}]]$ induced by the exponential map
$\exp: {\frak h} \to T$ we view ${\bold H}_q$ as a subalgebra of the
${\Bbb C}$-algebra ${\rm End}_{\,{\Bbb C}} {\Bbb C} [[{\frak h}]]$ (a
``specialized'' analogue of Proposition 12.3 with ${\bold q} = q$ and
no ``$\epsilon$'' at all). Following the same pattern as before
Proposition 12.3, we endow ${\Bbb C} [[{\frak h}]]$ with $J$-adic
filtration, where $J$ is the augmentation ideal in ${\Bbb C} [[{\frak
  h}]]$. We further introduce a {\it decreasing} filtration $\cal F^j,
j \in {\Bbb Z}$, on ${\rm End}_{\,{\Bbb C}} {\Bbb C} [[{\frak h}]]$ by the
formula
$$
\cal F^j {\rm End}_{\,{\Bbb C}} {\Bbb C} [[{\frak h}]] = \{ u \; | \; u(J^k)
\subset J^{k+j}, \; \mbox{ for all } k > j \}.
$$
and let $\cal F^j {\bold H}_q : = {\bold H}_q \cap \cal F^j {\rm End}_{\,{\Bbb C}}
{\Bbb C} [[{\frak h}]]$ be the induced filtration on ${\bold H}_q$. Note
that unlike the filtration $F^{\bullet} \widehat{{\bold H}}$ considered
before Proposition 12.3, the filtration $\cal F^j {\bold H}_q$ 
 is non-trivial for both {\it negative} and {\it
  positive} values of $j$.

Comparison of formulas (\ref{dem_lus}), (\ref{dima}) and Proposition
12.5 yield the following result

\medskip

\noindent
{\bf Proposition 12.6} {\it For any $q \neq 1$ there is a graded
  algebra isomorphism }
$$
{\bold H}_{nil} \simeq \mbox{gr}_{\cal F} {\bold H}_q. \quad \Box
$$

\medskip
We now turn to the geometric interpretation.
\bigskip

\noindent
{\bf Equivariant cohomology}  
\smallskip

 Let $M$ be a manifold equipped with smooth action of a Lie group $G$. If the
 action of $G$ is free, then $M/G$ is a manifold, but in general it is
 not. We define a space $M_G$ (the {\it homotopy quotient}) by
 replacing $M$ by another space with the same homotopy type on which
 $G$ acts freely. This is done by introducing a contractible space
 $EG$ on which $G$ acts freely and defining
$$
M_G= (EG \times M)/G.
$$
$M_G$ is well defined up to homotopy equivalence.

 The first and the second projections $EG \leftarrow M \times EG \to
M\;$ induce, by passing to $G$-orbits, a double fibration
$$
\begin{CD}
M_G @>{p}>> M/G \\
@V{\pi}VV \\
BG
\end{CD}
$$
where the {\it classifying space} $BG$ is defined as $BG=EG/G$.
Here we consider $M/G$ a space of orbits which may have
non-separated topology if the action of $G$ is not free (if $G$ is
compact then $M/G$ is always separated, but may have singularities
even for smooth $M$). The projection $\pi: M_G \to BG$ behaves much
better. It is a locally trivial fibration with fiber $M$.

We define the {\it equivariant cohomology} of $M$ by
$H_G^*(M)=H^*(M_G).$ Clearly $H^*_G(M)$ is a module over $ H^*(BG)$
via $J_1$.

If $G$ acts freely on $M$ then the map $p: M_G \to M/G$ becomes a
fibration with contractible fiber $EG$, and hence $$H^*_G(M) = H^*(M
/G). $$

\noindent
{\bf Example } The action of ${\Bbb C}^*$ on ${\Bbb C}^n \setminus
\{0\}$ is free. Although ${\Bbb C}^n$ is not contractible, we may form
$\displaystyle {\Bbb C}^{\infty} = \{ (x_i)_{i \in {\Bbb N}} \; | \;
x_i \in {\Bbb C}, \, x_i = 0$ for almost all $i \}$ and one shows that
${\Bbb C}^{\infty} \setminus \{0\}$ is a contractible space on which
${\Bbb C}^*$ acts freely. Thus, ${\Bbb C}^{\infty} \setminus \{0\}$ is
a model for $E{\Bbb C}^*$, and $B{\Bbb C}^*$ is homotopy equivalent to
${\Bbb C} {\Bbb P}^{\infty}= ({\Bbb C}^{\infty} \setminus \{0\})/{\Bbb
  C}^*$.
\medskip

\noindent
{\bf Remarks} 

(1)
In our applications $G$ will be a complex reductive Lie
group, and for such a group the spaces $EG$ and $BG$ are necessarily
infinite-dimensional. The reader who doesn't feel comfortable with
leaving the category of algebraic varieties may choose the approach of
[EG] where $EG$ and $BG$ are approximated by open subspaces of complex
vector spaces with linear $G$-action.

(2) If $K \subset G$, $G$ reductive, is a maximal compact subgroup,
then $K$ acts freely on $EG$, so we can consider $EG$ as a model for
$EK$. Hence we have a map $BK=EG/K \to BG=EG/G$ which is a homotopy
equivalence. Therefore, for any $G$-variety $M$, the space $M_G$ is
homotopy equivalent to $M_K$, hence $H^*_K(M) \simeq H^*_G(M)$.

\medskip

Let $M$ be a smooth complex algebraic variety with an algebraic action
of a linear reductive group $G$ and let $\cal F$ be a $G$-equivariant
sheaf on $M$. We may construct a sheaf $\cal F_G$ on $M_G$ as follows.
First, pullback $\cal F$ to $M \times EG$ via the first projection,
and regard the pull-back $pr_1^* \cal F$ as a $G$-equivariant sheaf on
$M \times EG$ relative to the {\it diagonal} $G$-action. The action
being free, the equivariant descent property [CG, 5.2.15] (cf. also
the first remark above) ensures that there exists a uniquely defined
sheaf $\cal F_G$ on $M_G$ such that $pr_1^* \cal F = \pi^* \cal F_G$.
 
 If  $\cal F$ is   locally free we define equivariant  characteristic
 classes of  $\cal F$ to be  the  corresponding classes of  $\cal F_G$
 (viewed as a vector bundle on $M_G$) in $H^*(M_G) = H^*_G(M)$. In the
 general case, since $M$   is smooth, we  can  choose a {\it   finite}
 $G$-equivariant resolution of  $\cal F$ by  locally  free sheaves and
 define equivariant   characteristic classes  via  the usual  formulas
 involving the classes of these locally free sheaves.

  For any finite-dimensional representation $V$ of $G$ we may form the
  associated    vector  bundle  $V_G$  on  $BG$  corresponding  to the
  principal $G$-bundle  $EG  \to BG$.  The assignment $V \mapsto V_G$
 gives a  ring homomorphism
  ${\bold  R}(G) \to   K(BG)$. It is   known  that the  homomorphism  is
  injective and $K(BG)= \widehat{{\bold R}(G)}$ is the completion of the
  image  of ${\bold R}(G)$  at the  augmentation  ideal. Thus  if $G$ is
  reductive we have
$$
H^*(BG, {\Bbb C})\  \stackrel{Chern\atop charact.}\simeq \ {\Bbb C}
\otimes K(BG) \; = \; {\Bbb C} \otimes \widehat{{\bold R}(G)} \; = \;
\widehat{{\Bbb C} [G]^G} \; \simeq \; {\Bbb C} [[ {{\frak g}}]]^G,
$$ 
where the last isomorphism is given by pulling back
$\mbox{Ad}\,G$-invariant functions from $G$ to ${\frak g} = Lie \; G$
via the exponential map $exp: {\frak g} = Lie\ G \to G$.

Similarly, for  most of the $G$-spaces  we consider, e.g. for cellular
fibrations   ([CG, 5.5]),  the Chern  character   gives  an isomorphism
$\widehat H_G(M) = \widehat{K^G(M)}$,  where $\widehat  H_G(M)$ stands
for  the completion  of $H^*_G(M)$  at the  augmentation  ideal of the
ground ring $H^*(BG)$, i.e.
$$
\widehat H_G(M) = \prod_{i=0}^{\infty}H^i_G (M). 
$$

\bigskip

Given a smooth $G$-manifold ${\Bbb M}$ and a closed $G$-stable subset $Z
\subset {\Bbb M}$, we define  the $G$-equivariant cohomology of ${\Bbb M}$
with support in $Z$ by the formula
$$
 H^*(Z |{\Bbb M}; G):= H^*_G ({\Bbb M} , {\Bbb M} \setminus Z)= H^*({\Bbb M}_G,
({\Bbb M} \setminus Z)_G)
$$
(compare this definition with (\ref{pd_defn_bm})).

Recall now the setup of Section 6: we have  the Steinberg variety $Z$,
embedded naturally as a $G$-invariant subvariety in  ${\Bbb M} = T^*\cal
B  \times T^* \cal B$.   Since  $Z \circ Z =Z$   we can repeat all our
constructions  with  convolutions  in   the ``relative situation  over
$BG$''  (i.e.  we view $Z_G$ as   a fibration over   $BG$ which has an
embedding  into ${\Bbb M}_G$,   an infinite-dimensional smooth manifold,
commuting with the projection to $BG$,  etc.)  to define the structure
of a convolution algebra on $H^*(Z | {\Bbb M}; G)$.

\medskip
We are going to interpret the degenerate Hecke algebra geometrically
using equivariant cohomology, in the same way as the affine Hecke
algebra ${\bold H}$ was described via equivariant K-theory.
Identify $H^{2\bullet}_{{\Bbb C}^*}(pt)=H^{2\bullet }({\Bbb C} {\Bbb
 P}^{\infty})$ with the ring ${\Bbb C}[\epsilon] $,
so that a generator $\epsilon \in H^2({\Bbb C} {\Bbb P}^{\infty})$ is
viewed as a {\it degree 1} element in the polynomial ring ${\Bbb C} [
\epsilon]$.

\medskip

\noindent
{\bf Theorem 12.7} {\em Let $Z$ be the Steinberg variety of $G$ with
  the natural action of $G \times {\Bbb C}^*$. Then the convolution 
algebra $H^{2 \bullet}
(Z | {\Bbb M}; G \times {\Bbb C}^*)$  is naturally isomorphic to
the degenerate affine Hecke algebra   ${\bold H}_{deg}$.} 

\medskip

\noindent
{\bf Remarks}

(1) We write $H^{2\bullet} (Z | {\Bbb M}; G \times {\Bbb C}^*)$ since
 there are no equivariant odd cohomology groups with support in $Z$, and
the isomorphism of Theorem 12.7 is doubling degrees, i.e. the
 degree of an element in the equivariant cohomology is  twice the
 degree of the corresponding element in ${\bold H}_{deg}$.

(2) Theorems 12.7, 11.5 and 11.4  allow us to express all the 
results on convolution algebras of the
Steinberg variety $Z$ presented in these lectures, in the following diagram:

$$
\begin{CD}
{\bold H} @=  K^{G \times {\Bbb C}^*}(Z) @>{ch^*_{G \times {\Bbb C}^*}}>>
\widehat H (Z | {\Bbb M}; G \times {\Bbb C}^*) @= \widehat{{\bold H}}_{deg}
\\
@VV{{\bold q} =1}V  @VV{\mbox{\scriptsize forgetting}\atop {\Bbb
  C}^*-\mbox{\scriptsize action}}V
@VV{\mbox{\scriptsize forgetting}\atop {\Bbb C}^*-\mbox{\scriptsize
    action}}V @VV{\epsilon =0}V 
\\
{\Bbb Z} [ \widetilde W] @= K^G(Z) @>{ch^*_G}>> \widehat H (Z | {\Bbb M};
G) @= {\Bbb C} [[{\frak h} ]][W]
\\
@VV{\mbox{\scriptsize proj}}V 
@VV{\mbox{\scriptsize forgetting}\atop G-\mbox{\scriptsize action}}V 
@VV{\mbox{\scriptsize forgetting}\atop G-\mbox{\scriptsize action}}V 
@VV{\mbox{\scriptsize proj}}V
\\
\cal A[W] @= K_{{\Bbb C}}(Z) @>{ch^*}>> H_*(Z, {\Bbb C}) @= \cal A [W]
\\
@VV{\mbox{\scriptsize ev}}V 
@VV{\mbox{\scriptsize support}\atop \mbox{\scriptsize cycle}}V   
@VV{\mbox{\scriptsize top}\atop \mbox{\scriptsize homology}}V 
@VVV
\\
{\Bbb C}[W] @= H(Z) @= H(Z) @= {\Bbb C} [W]
\end{CD}
$$ 
\smallskip

\noindent
where going  from top to bottom  leads to forgetting some amount of
structure  and where the  following notation has been used:\\ 
$\bullet \quad \cal  A = {\Bbb C}   [T]/ I_T \stackrel{exp^*}\simeq  {\Bbb C} [[
{\frak  h}]]/ I_{{\frak  h}}  \ \stackrel{Borel\atop isom}\simeq  H^*(\cal
B)$, where \\ 
$\bullet \quad I_T  =$ ideal  generated  by  $W$-invariant functions on  $T$
vanishing at 1, \\  
$\bullet \quad I_{\frak h}=$ ideal generated  by $W$-invariant power  series
on  $\frak  h$ vanishing at 0,\\   
$\bullet \quad \mbox{ev}:\cal A [W]  \to {\Bbb C} [W]$ is  taking value at  $1
\in T$, resp. $0 \in \frak  h$,\\ 
$\bullet \quad exp^*: {\Bbb C} [T]/ I_T \to {\Bbb C} [[ {\frak h}]]/ I_{\frak h} $
is the pullback via $exp:  {\frak h} \to T$, \\   
$\bullet \quad  \mbox{proj}: {\bold R}(T) \to  \cal A$, resp.
 $\mbox{proj}: {\Bbb C} [[ {\frak h}]] \to  \cal A$ is  a natural projection, \\ 
$\bullet \quad   \widehat{{\bold  H}}_{deg}=$   completed  version of   ${\bold
H}_{deg}$, with $S^{\bullet}({\frak h}^*)[\epsilon]$  replaced by ${\Bbb C}
[[{\frak h}, \epsilon]]$.

\bigskip

\noindent
{\bf Sketch of proof of Theorem 12.7} Recall that $\mu: T^*\cal B \to
\cal N$ denotes the Springer resolution. Since  $\cal B = \mu^{-1}(0)$,
we have  a  natural $H^*(Z  |  {\Bbb M}; G  \times  {\Bbb C}^*)$-action on
$$H^*_{G \times  {\Bbb  C}^*}(\cal B) = H^*_{G  \times  {\Bbb C}^*} (G/ B)=
H^*_{B   \times  {\Bbb  C}^*}(pt)   =   H^*_{T \times  {\Bbb C}^*}(pt)   =
{\Bbb C}[{\frak h}, \epsilon].$$

 This action gives an algebra homomorphism 
$$
\rho: H^*(Z| {\Bbb M}; G \times {\Bbb C}^*) \to {\rm End}_{\;{\Bbb C} [\epsilon]}
{\Bbb C} [{\frak h}, \epsilon].
$$
One shows first, by the same argument as has been used in [CG, Claim
7.6.7] 
to prove a similar result in equivariant K-theory, that the map
$\rho$ is injective. The isomorphism of the theorem will then follow
from Proposition 12.5 provided we show that the image of $\rho$
coincides with the subalgebra of ${\rm End}_{\;{\Bbb C} [\epsilon]}
{\Bbb C} [{\frak h}, \epsilon]$ described in Proposition 12.5. To prove
the latter we argue as follows.

View $ H^*(Z| {\Bbb M}; G \times {\Bbb C}^*)$ as a subalgebra of
${\rm End}_{\;{\Bbb C} [\epsilon]} {\Bbb C} [{\frak h}, \epsilon]$ via $\rho$ and 
recall the filtration $F^{\bullet}  {\rm End}_{\;{\Bbb C} [\epsilon]} {\Bbb C}
[{\frak h}, \epsilon]$ introduced after Proposition 12.2. The
convolution product on $H^*(Z| {\Bbb M}; G \times {\Bbb C}^*)$ is
continuous in the topology defined by the filtration $F^{\bullet}$,
hence extends to the completion. This way one makes 
$ \widehat{H}(Z| {\Bbb M}; G \times {\Bbb C}^*):=
\prod_{i=0}^{\infty} H^i(Z| {\Bbb M}; G \times {\Bbb C}^*)$ a filtered
algebra under convolution and defines an action of this algebra on
$$
\widehat{H}_{G \times {\Bbb C}^*}(\cal B) = \prod_{i=0}^{\infty} H^i_{G 
\times {\Bbb C}^*}(\cal B) = {\Bbb C} [{\frak h}, \epsilon].
$$

Thus, the homomorphism  $ \rho: H^*(Z| {\Bbb  M}; G \times {\Bbb C}^*) \to
{\rm End}_{\;{\Bbb C} [\epsilon]}  {\Bbb C} [{\frak  h}, \epsilon] $, given by the
action, can be extended to a continuous algebra homomorphism
$$
\widehat{\rho}: \widehat{H}(Z| {\Bbb M}; G \times {\Bbb C}^*) \to 
{\rm End}_{\;{\Bbb C} [[\epsilon]]\;} {\Bbb C} [[{\frak h}, \epsilon]].
$$

 Further, we use an equivariant version (i.e. relative version for
 fibrations over $BG$) of the Bivariant Riemann-Roch Theorem [CG,
 Theorem 5.11.11] to construct an algebra homomorphism
$$
RR: K^{ G \times {\Bbb C}^*}(Z) \to \widehat{H}(Z| {\Bbb M}; G \times {\Bbb
C}^*).
$$

 We now apply the isomorphism ${\bold H} \simeq K^{ G \times {\Bbb
 C}^*}(Z)$ of Theorem 11.5 so that composing the maps $RR$ and
 $\widehat{\rho}$ yields an algebra homomorphism
$$
\widehat{\rho} \circ RR:\; {\bold H} \simeq K^{ G \times {\Bbb C}^*}(Z) \to
{\rm End}_{\;{\Bbb C} [[\epsilon]]\;} {\Bbb C} [[{\frak h}, \epsilon]].
$$
 By our construction, we see by means of Remark (2) after Theorem 11.5
 that the map $\widehat{\rho} \circ RR$ above, arising from the
 geometric convolution action, is in effect nothing but the embedding
 ${\bold H} \hookrightarrow {\rm End}_{\;{\Bbb C} [[\epsilon]]\;} {\Bbb C} [[{\frak
 h}, \epsilon]]$ used in Proposition 12.3 for the definition of the
 filtration $F^{\bullet}$ on $\widehat{{\bold H}}$. In particular, the
 map
$$
\widehat{\rho \circ RR}:\; \widehat{{\bold H}} = {\Bbb C} [[\epsilon]]
\otimes_{{\Bbb Z} [{\bold q}, {\bold q}^{-1}]} {\bold H} \to
{\rm End}_{\;{\Bbb C} [[\epsilon]]\;} {\Bbb C} [[{\frak h}, \epsilon]]
$$
induced by $\widehat{\rho} \circ RR$ is compatible with the filtration
$F^{\bullet}$ on both sides. Furthermore, the isomorphism $\mbox{gr}_F
\widehat{{\bold H}} \simeq {\bold H}_{deg}$ of Proposition 12.3 implies
that the image of the associated graded map
$$
\mbox{gr}_F(\widehat{\rho \circ RR}):\; \mbox{gr}_F(\widehat{{\bold H}}) \to 
\mbox{gr}_F({\rm End}_{\;{\Bbb C} [[\epsilon]]\;} {\Bbb C} [[{\frak h}, \epsilon]])
$$
can be identified with ${\bold H}_{deg}$ if we identify 
$\mbox{gr}_F({\rm End}_{\;{\Bbb C} [[\epsilon]]\;} {\Bbb C} [[{\frak h}, \epsilon]])$
with ${\rm End}_{\;{\Bbb C} [\epsilon]\;} {\Bbb C} [{\frak h}, \epsilon]$ in a
natural way (the associated graded of ${\Bbb C} [[{\frak h}, \epsilon]]$
with respect to the $I$-adic filtration is equal to ${\Bbb C} [{\frak h}, 
\epsilon]$). Also, one clearly has $\mbox{gr}_F \widehat{H}(Z|{\Bbb M};
G \times {\Bbb C}^*) \simeq H^*(Z|{\Bbb M};
G \times {\Bbb C}^*)$.

 Thus we conclude:
$$
\mbox{Image}(\rho) = \mbox{Image}(\mbox{gr}_F \widehat{\rho})
\stackrel{\alpha}\hookrightarrow \mbox{Image}\
\mbox{gr}_F(\widehat{\rho \circ RR}) = {\bold H}_{deg}.
$$

 To complete the proof of the theorem it suffices to show that the
 inclusion $\alpha$ above is in fact an equality. This follows from
 the fact that the map $RR$ has dense image, which is proved by the
 same argument as was used in [CG, Theorem 6.2.4] to show that the
 non-equivariant Riemann-Roch map ${\Bbb C} \otimes_{{\Bbb Z}} K(Z) \to
 H_*(Z)$ is a bijection. $\Box$

\bigskip

\pagebreak
\noindent
{\bf Geometric realization of the Nil-Hecke algebra}

\smallskip

 The algebra construction of Corollary 2.1 in the special case of $M =
 \cal B$, $N = pt$ (see Example (ii) after Corollary 2.2) makes
 $H_*(\cal B \times \cal B)$ an associative algebra with convolution
 product, and makes $H_*(\cal B)$ a simple $H_*(\cal B \times \cal
 B)$-module. Since $\cal B$ is smooth we have by Poincar\'e duality
 $H_*(\cal B) \simeq H^*(\cal B)$. Applying further the Borel
 isomorphism we may identify $H_*(\cal B)$ with ${\Bbb C} [{\frak
 h}]/I_{{\frak h}}$, where $I_{{\frak h}}$ is the ideal generated by
 $W$-invariant polynomials without constant term.

Recall that $G$-diagonal orbits on $\cal B \times \cal B$
are in canonical bijection with the Weyl group, $W$ (see above
Proposition 6.9); we write
 $[O_w] \in H_*(\cal B \times \cal B)$ for the fundamental
class of the closure, $\overline{O}_w$, of the
orbit corresponding to $w \in W$.
 The following result is essentially due to [BGG].

\medskip
\noindent
{\bf Proposition 12.8} \  {\it

{\rm (i)} The classes $[O_w]$, $w \in W$, span a \,$\#W$-dimensional subalgebra
in the convolution algebra $H_*(\cal B \times \cal B)$. Moreover, the
map $r_w \mapsto [O_w]$, $w \in W$, gives an isomorphism of this
subalgebra with the algebra} {\bf Nil}.

{\it {\rm (ii)} For any simple reflection $s_{\alpha}$, the action of the
  class $[O_{s_{\alpha}}]$ in the 
$H_*(\cal B \times \cal B)$-module $H_*(\cal B)
  \simeq {\Bbb C} [{\frak h}]/ I_{{\frak h}}$ is given by the operator
  $R_{\alpha}$ defined in {\rm (\ref{dima})}.} $\Box$
\medskip

\noindent
{\bf Remarks}

(1) The operators $R_{\alpha}$ in (\ref{dima}) commute with the
multiplication operator by a $W$-invariant polynomial, hence descend
to well-defined operators on ${\Bbb C} [{\frak h}]/ I_{{\frak h}}$.

(2) One can get an analogue of Proposition 12.8 with the whole
nil-Hecke algebra ${\bold H}_{nil}$ instead of {\bf Nil}, and with the
operators $R_{\alpha}$ acting on ${\Bbb C}[{\frak h}]$ instead of 
${\Bbb C} [{\frak h}]/ I_{{\frak h}}$, by replacing the ordinary homology by
$G$-equivariant (co-)homology; see [Ar].

\medskip

 Further, it follows from the formula of Corollary 6.10 that the map
 $c^{\mbox{\footnotesize biv}}$ gives an 
embedding $c^{\mbox{\footnotesize biv}}: {\Bbb C}[W] \hookrightarrow
 H_*(\cal B \times \cal B)$. We can therefore define a filtration on
 ${\Bbb C} [W]$ by 
$$
E_j({\Bbb C}[W]) = \{u \; | \; c^{\mbox{\footnotesize biv}}(u) \in\; \bigoplus_{i \leq j}
H_{2(n+i)}(\cal B \times \cal B) \},
$$
which is analogous to the filtration used in Proposition 12.3.
From Proposition 6.9 and Corollary 6.10 we obtain the following result:

\medskip

\noindent
{\bf Corollary 12.9} (i) {\it For any simple reflection $s_{\alpha}$,
the convolution action of $c^{\mbox{\footnotesize biv}}(s_{\alpha})\in
H_*(\cal B \times \cal B)$ on $H_*(\cal B)={\Bbb C} [{\frak h}]/ I_{{\frak
h}}$
is given by the operator $S_\alpha$, see $(\ref{nilki})$, 
specialized at $\epsilon=1$.}

(ii) {\it The assignment $w\mapsto r_w\,,\, w\in W,$ gives a graded algebra isomorphism}
$$
\mbox{gr}^E{\Bbb C}[W] \simeq \mbox{\bf Nil}. 
$$

\bigskip

\noindent
{\bf Final remark } Most of  the constructions of  this section can be
extended to quantized enveloping algebra instead of Hecke algebra. The
analogue of the degenerate Hecke algebra ${\bold  H}_{deg}$ is known, in
the  context of  quantized enveloping  algebras, as  the Yangian,  see
[CP]. There is an analogue of Theorem 12.7 for Yangians.

\addtolength{\baselineskip}{-.055\baselineskip}

\end{document}